\documentclass[preprint]{imsart} 

\RequirePackage[OT1]{fontenc}
\RequirePackage{amsthm,amsmath}
\RequirePackage[numbers]{natbib}
\RequirePackage[colorlinks,citecolor=blue,urlcolor=blue]{hyperref}
\usepackage{amsmath}
\usepackage{graphicx} 
\usepackage{epstopdf}
\usepackage{float}
\usepackage{fixmath}
\usepackage{bm}
\usepackage{inputenc,t1enc}
\usepackage{bbm}
\usepackage{amsfonts}
\usepackage{latexsym}
\usepackage{mathrsfs}
\usepackage{mathtools,xparse}
\usepackage{enumitem}  
\usepackage{subcaption}
\usepackage{amssymb} 
\usepackage{multirow}

\startlocaldefs
\numberwithin{equation}{section}
\theoremstyle{plain}

\endlocaldefs

\startlocaldefs
\numberwithin{equation}{section}
\theoremstyle{plain}
\newcommand{\nc}{\newcommand}
\nc{\nt}{\newtheorem}
\nt{defn}{Definition}
\nt{lem}{Lemma}
\nt{pr}{Proposition}
\nt{theorem}{Theorem}
\nt{cor}{Corollary}
\nt{remark}{Remark}
\nt{ex}{Example}
\nt{ass}{Assumption}
\nt{step}{Step}
\nt{case}{Case}
\nt{subcase}{Subcase}
\nt{note}{Note}
\nt{assum}{Assumption}
\nc{\bd}{\begin{defn}}
\nc{\ed}{\end{defn}}
\nc{\blem}{\begin{lem}} 
\nc{\elem}{\end{lem}}
\nc{\bpr}{\begin{pr}} 
\nc{\epr}{\end{pr}}
\nc{\bcor}{\begin{cor}} \nc{\ecor}{\end{cor}}
\nc{\bex}{\begin{ex}}  \nc{\eex}{\end{ex}}
\nc{\bass}{\begin{ass}}  \nc{\eass}{\end{ass}}
\nc{\bstep}{\begin{step}}  \nc{\estep}{\end{step}}
\nc{\bcase}{\begin{case}}  \nc{\ecase}{\end{case}}
\nc{\bsubcase}{\begin{subcase}}  \nc{\esubcase}{\end{subcase}}
\nc{\bnote}{\begin{note}}  \nc{\enote}{\end{note}}
\nc{\bassum}{\begin{assum}}  \nc{\eassum}{\end{assum}}
\nc{\prf}{{\bf Proof.} }
\nc{\eop}{\hfill $\Box$ \\ \\}
\nc{\argmin}{\mathrm{argmin}}
\nc{\argmax}{\mathrm{argmax}}
\nc{\sgn}{\mathrm{sgn}}
\nc{\Var}{\mathrm{Var}}
\nc{\Cov}{\mathrm{Cov}}
\nc{\bak}{\!\!\!\!\!}
\nc{\IBD}{\mathrm{IBD}}
\nc{\supp}{\mathrm{supp}}
\nc{\dom}{\mathrm{dom}}
\nc{\R}{{\mathbb R}}
\nc{\peq}{\preceq}
\nc{\wt}{\widetilde}
\nc{\Mult}{\mathrm{Mult}}

\nc{\Prob}[1]{\mathbb{P}_{#1}}

\endlocaldefs

\begin{document}

\begin{frontmatter}

\title{Stacked Grenander and rearrangement estimators of a discrete distribution}
\runtitle{Stacked Grenander and rearrangement estimators}
\runauthor{V. Pastukhov}
\begin{aug}
\author{\fnms{Vladimir} \snm{Pastukhov}\corref{craut1}\ead[label=e1]{vlapas@chalmers.se}}
\address{Department of Computer Science and Engineering,\\
Chalmers University of Technology\\
\printead{e1}}
\end{aug}


\begin{abstract}
In this paper we consider the stacking of isotonic regression and the method of rearrangement with the empirical estimator to estimate a discrete  distribution with an infinite support. The estimators are proved to be strongly consistent with $\sqrt{n}$-rate of convergence. We obtain the asymptotic distributions of the estimators and construct the asymptotically correct conservative global confidence bands. We show that stacked Grenander estimator outperforms the stacked rearrangement estimator. The new estimators behave well even for small sized data sets and provide a trade-off between goodness-of-fit and shape constraints.
\end{abstract}



\tableofcontents
\end{frontmatter}

\section{Introduction}

This work is largely inspired by recent papers in the estimation of discrete distributions with shape constraints. The first paper in this area is \cite{jankowski2009estimation}, where the authors studied the method of rearrangement and maximum likelihood estimator (MLE) of probability mass function (p.m.f.) under monotonicity constraint. The MLE under monotonicity constraint is also known as Grenander estimator. Next, in the paper \cite{durot2013} the authors introduced the least squares estimator of a discrete distribution under the constraint of convexity and, further, its limiting distribution was obtained in \cite{balabdurkol}. Furthermore, the MLE of log-concave p.m.f. was studied in detail in \cite{baljanruf}, and in \cite{jankowskitian} the problem was generalised to the case of multidimensional discrete support. Next, in paper \cite{baljan} the authors introduced the MLE of unimodal p.m.f. with unknown support, proved the consistency and obtained the asymptotic distribution. The problem of least squares estimation of a completely monotone p.m.f. was considered in papers \cite{balfour, balkul}.

In most of the papers listed above the authors considered both the well- and the mis-specified cases and studied the asymptotic properties of the estimators in both cases. In this work we do not have the mis-specified case in a sense that we assume that the true p.m.f. can be non-monotone and our estimators are strongly consistent even if the true p.m.f. is not decreasing. 

The estimators introduced and studied in this paper are in some sense similar to nearly-isotonic regression approach, cf. \cite{tibshirani2011} and \cite{minami2020} for multidimensional case. Nearly-isotonic regression is a convex optimisation problem, which provides intermediate less restrictive solution and the isotonic regression is included in the path of the solutions. 

At the same time, our approach is in some sense opposite to liso (lasso-isotone), cf. \cite{fang2012}, and to bounded isotonic regression, cf. \cite{Luss2017}. The liso is a combination of isotonic regression and lasso penalties, and bounded isotonic regression imposes additional penalisation to the range of the fitted model. 

In this paper we combine Grenander estimator and the method of rearrangement with cross-validation-based model-mix concept, cf. \cite{stone1974}. The estimator is constructed as a convex combination of the empirical estimator and Grenander estimator or the empirical estimator and rearrangement estimator. Following the terminology for regression and classification problems in \cite{breiman1995, leblanktibshirani1996, wolper1992}, we call the resulting estimators as \textit{stacked Grenander estimator} and \textit{stacked rearrangement estimator}, respectively. Therefore, we do not impose the strict monotonic restriction and let the data decide. 

There are several papers where the authors studied a convex combination of the empirical estimator with a prescribed probability vector, cf. \cite{fienberg1972, fienberg1973, stone1974, tribula1958}. In particular, in \cite{stone1974} the authors proposed the combination of the empirical estimator and a constant p.m.f. with a mixture parameter selected by cross-validation. Also, the minimax estimator of a p.m.f. with respect to $\ell_{2}$-loss with a fixed known finite support and sample size $n$ is given by a convex combination of the empirical estimator and the uniform distribution with a mixture parameter equal to $\frac{\sqrt{n}}{n + \sqrt{n}}$, cf. \cite{tribula1958}. Furthermore, in \cite{fienberg1973} the authors provide a geometrical explanation on the gain from stacking the empirical estimator with a fixed probability vector and show that the improvement of the estimation increases as the size of the support becomes larger. 

In the case of continuous support the first paper on the density estimation via stacking is \cite{smyth1999}, where it is shown that the method of stacking performs better than selecting the best model by cross-validation. Next, in \cite{rigollet2007} the authors studied the approach of linear and convex aggregation of density estimators and, in particular, proved that the aggregation of two estimators allows to combine the advantages of both. To the authors’ knowledge the constrained stacked estimators have not been investigated for the case of continuous density.

To the authors' knowledge, the problem of staking the shape constrained estimators has not been studied much even in a regression setup, except for the paper \cite{wrigh1978}. In the paper \cite{wrigh1978} the author used a convex combination of linear regression with isotonic regression to obtain a strictly monotonic solution. Also, it is worth to mention the paper \cite{Hastie2020}, where it was shown that in terms of prediction accuracy the simplified relaxed lasso (which is stacking of least squares estimator and lasso) performs almost equally to the lasso in low signal-to-noise ratio regimes, and nearly as well as the best subset selection in high signal-to-noise ratio scenarios.


The paper is organised as follows. In Section \ref{statement_pr} we state the problem and
introduce notation. The derivation of cross-validation based mixture parameter is given in Section \ref{alp_sel}. Section \ref{est_prpts} is dedicated to the theoretical properties of the estimators such as consistency, rate of convergence and asymptotic distribution. Also, in Section \ref{est_prpts} we construct asymptotic confidence bands. In Section \ref{est_sim} we do simulation study to compare the performance of the estimators with empirical, minimax, rearrangement and Grenander estimators. The article closes with a conclusion and a discussion of possible generalisations in Section \ref{concl}. The ancillary results and the proofs of some statements are given in Appendix. The R code for the simulations is available upon request.

\section{Statement of the problem and notation}\label{statement_pr}
First, let us introduce notation and several definitions. Assume that $z_{1}, z_{2}, \dots, z_{n}$ is a sample of $n$ i.i.d. random variables with values in $\mathbb{N}$ and generated by a p.m.f. $\bm{p}$. For a given data sample let us create the frequency data $\bm{x} = (x_{0}, \dots, x_{t_{n}})$, where $x_{j}= \sum_{i=1}^{n} 1 \{ z_{i} = j \}$ and $t_{n} = \sup\{j: x_{j}>0 \}$ denotes the largest order statistic for the sample.

The empirical estimator of $\bm{p}$ is given by
\begin{equation*}\label{unrMLEp}
	\hat{p}_{n, j} = \frac{x_{j}}{n}, \quad j \in \mathbb{N},
\end{equation*}
and it is strongly consistent, unbiased and asymptotically normal in $\ell_{2}$-space.

The rearrangement estimator studied in \cite{jankowski2009estimation} is defined as 
\begin{equation}\label{estrear}
	\hat{\bm{r}}_{n} = rear(\hat{\bm{p}}_{n}),
\end{equation}
where $rear(\bm{w})$ denotes the reversed-ordered vector. Also, equivalently, the rearrangement estimator can be written as  $\hat{r}_{n,j} = \sup \{u: Q_{n}(u) \leq j \}$, where $Q_{n}(u) = \#\{k: \, \hat{p}_{n, k} \geq u \}$.
 
The MLE of decreasing p.m.f., or Grenander estimator, which we denote by $\hat{\bm{g}}_{n}$, is equivalent to the isotonic regression of the empirical estimator, cf. \cite{barlowstatistical, jankowski2009estimation, robertsonorder}, i.e.
\begin{equation}\label{orMLEp}
\hat{\bm{g}}_{n} = \Pi(\hat{\bm{p}}_{n}|\mathcal{F}^{decr}) :=  \underset{\bm{f} \in \mathcal{F}^{decr}}{\argmin}\sum_{j}[\hat{p}_{n, j} - f_{j}]^{2},
\end{equation}
where $\mathcal{F}^{decr}$ is the monotonic cone in $\ell_{2}$, i.e. $\mathcal{F}^{decr} = \Big\{\bm{f} \in \ell_{2}: \, f_{0} \geq f_{1} \geq  \dots \Big\}$, $\hat{\bm{p}}_{n}$ is the empirical estimator and $\Pi(\hat{\bm{p}}_{n}|\mathcal{F}^{decr})$ denotes the $\ell_{2}$-projection of $\hat{\bm{p}}_{n}$ onto $\mathcal{F}^{decr}$
\footnote{The notion of "isotonic regression" in (\ref{orMLEp}) might be confusing. Though,
for historical reasons, it is a standard notion in the subject of constrained inference, cf. the
monographs \cite{robertsonorder, silvapsen} and also papers \cite{best1990, stout2013}, dedicated to the computational aspects, where
the notation "isotonic regression" is used for the isotonic projection of a general vector.}.

In our work we construct the estimator in the following way: 
\begin{equation}\label{SGest}
\hat{\bm{\phi}}_{n} = \beta \hat{\bm{h}}_{n} + (1 -\beta)\hat{\bm{p}}_{n},
\end{equation}
where
\begin{equation*}\label{}
\hat{\bm{h}}_{n} = \begin{cases}
    \hat{\bm{r}}_{n}, &\text{for the stacked rearrangement estimator}, \\
    \hat{\bm{g}}_{n}, &\text{for the stacked Grenander estimator},
  \end{cases}
\end{equation*}
with the data-driven selection of $\beta$:
\begin{equation*}
\hat{\beta}_{n} = \underset{\beta \in [0, 1]}{\argmin } \, CV(\beta),
\end{equation*}
where $CV(\beta)$ is a cross-validation criterion, which we introduce and study below. 

We associate each component $x_{j}$ of the frequency vector $\bm{x}$ with multinomial indicator $\bm{\delta}^{[j]} \in \mathbb{R}^{t_{n}+1}$, given by  
\begin{equation}\label{indgm}
\bm{\delta}^{[j]} = (0, \dots, 0, 1, 0, \dots, 0)
\end{equation}
for $j = 0, \dots, t_{n}$, cf. \cite{stone1974}. All elements of $\bm{\delta}^{[j]}$ are zeros, except for the one with index $j$.



Next, let $\hat{\bm{p}}^{\backslash [j]}_{n}$ for $j=0, \dots, t_{n}$ denote the leave-one-out version of the empirical estimator $\hat{\bm{\phi}}_{n}$ for the frequency data $\bm{x} = (x_{0}, \dots, x_{t_{n}})$, i.e. for $j$ such that $x_{j} > 0$ let
\begin{equation*}
\hat{\bm{p}}^{\backslash [j]}_{n} = \frac{\bm{x} - \bm{\delta}^{[j]}}{n-1}.
\end{equation*}
Next, for the rearrangement estimator, the leave-one-out version is given by
\begin{equation*}
\hat{\bm{r}}^{\backslash [j]}_{n} = rear(\hat{\bm{p}}^{\backslash [j]}_{n}),
\end{equation*}
and for Grenander estimator:
\begin{equation*}
\hat{\bm{g}}^{\backslash [j]}_{n} = \Pi \big(\hat{\bm{p}}^{\backslash [j]}_{n}|\mathcal{F}^{decr}\big).
\end{equation*}
Therefore, for $j$ such that $x_{j} > 0$ the leave-one-out versions of stacked rearrangement and stacked Grenander estimators for a fixed misture parameter $\beta$ are given by
 \begin{equation}\label{cvestimj}
\hat{\bm{\phi}}^{\backslash [j]}_{n} =  \beta \, \hat{\bm{h}}^{\backslash [j]}_{n} + (1 - \beta)\hat{\bm{p}}^{\backslash [j]}_{n},
\end{equation}
with $\hat{\bm{h}}^{\backslash [j]}_{n} = \hat{\bm{r}}^{\backslash [j]}_{n}$ for the case of stacked rearrangement estimator, and $\hat{\bm{h}}^{\backslash [j]}_{n} = \hat{\bm{g}}^{\backslash [j]}_{n}$ for the case of stacked Grenander estimator, respectively.



For an arbitrary vector $\bm{f}\in\ell_{k}$ we define $\ell_{k}$-norm
\begin{equation*}\label{}
  ||\bm{f}||_{k}=\begin{cases}
    \Big(\sum_{j=0}^{\infty}|f_{j}|^{k}\Big)^{1/k}, & \text{if } \, k \in \mathbb{N} \backslash \{0\}, \\
    \sup_{j\in\mathbb{N}}|f_{j}|, & \text{if} \,  k=\infty,
  \end{cases}
\end{equation*}
and for $\bm{v} \in \ell_{2}$ and $\bm{w} \in \ell_{2}$ let $\langle \bm{v},\bm{w} \rangle = \sum_{j=0}^{\infty}v_{j}w_{j}$ denote the inner product on $\ell_{2}$.

For a random sequence $b_{n} \in \mathbb{R}$ we will use the notation $b_{n} = O_{p}(n^{q})$ if for any $\varepsilon >0$ there exists  a finite $M > 0$ and a finite $N > 0$ such that  
\begin{equation*}
\mathbb{P}[n^{-q}|b_{n}| > M] < \varepsilon,
\end{equation*}
for any $n > N$.

\section{Data-driven selection of the mixture parameter $\beta$}\label{alp_sel}
Let us consider squared $\ell_{2}$-distance between the true p.m.f. $\bm{p}$ and the stacked estimator $\hat{\bm{\phi}}_{n}$:
\begin{equation}\label{Lsqdist}
\begin{aligned}
L_{n}{} &= {}  ||\hat{\bm{\phi}}_{n} - \bm{p}||^{2}_{2}:\equiv L_{n}^{(1)} - 2 L_{n}^{(2)} + L_{n}^{(3)},
\end{aligned}
\end{equation}
where $L_{n}^{(1)} = \sum_{j=0}^{t_{n}}\hat{\phi}_{n,j}^{2}$, $L_{n}^{(2)} = \sum_{j=0}^{t_{n}}\hat{\phi}_{n,j}p_{j}$ and $L_{n}^{(3)} = \sum_{j=0}^{t_{n}}p_{j}^{2}$.

We aim to minimise $L_{n}$. Obviously, $\bm{p}$ is unknown, and we will use the approach introduced in \cite{ouyang2006} to estimate $L_{n}$. First, note that $L_{n}^{(3)}$ is a constant and can be omitted. Next, note that for a given $n$ we have for $L_{2}$ we have 
\begin{equation*}
L_{2} = \sum_{j=0}^{t_{n}}\hat{\phi}_{n,j}p_{j} = \mathbb{E}[\hat{\bm{\phi}}_{n}],
\end{equation*}
and following \cite{ouyang2006} we estimate $L_{n}^{(2)}$ by 
\begin{equation*}
\hat{L}_{n}^{(2)} = \sum_{j=0}^{t_{n}} \hat{p}_{n,j} \hat{\phi}^{\backslash [j]}_{n,j},
\end{equation*}
with $\hat{\bm{\phi}}^{\backslash [j]}_{n}$ defined in (\ref{cvestimj}). 
Therefore, we select the mixture parameter $\beta$ to minimise 
\begin{equation}\label{lsCV}
CV(\beta) = L_{n}^{(1)} - 2\hat{L}_{n}^{(2)},
\end{equation}
i.e.
\begin{equation*}
\hat{\beta}_{n} = \underset{\beta \in [0, 1]}{\argmin } \, CV(\beta).
\end{equation*}
This cross-validation approach for estimation of discrete distributions was first introduced in \cite{ouyang2006} for smoothing kernel estimator and was also used in, for example, \cite{chu2015,chu2017,racine2020}. The mixture parameter $\hat{\beta}_{n}$ is given in the following theorem.
 
\begin{theorem}\label{thmLSCV}
The leave-one-out  least-squares cross-validation mixture parameter $\hat{\beta}_{n}$ is given by
\begin{equation*}\label{}
  \hat{\beta}_{n}=\begin{cases}
    \frac{b_{n}}{a_{n}}, & \text{if } \, a_{n} \neq 0 \text{ and } 0 \leq b_{n} \leq a_{n}, \\
    1, & \text{if } \, 0 < a_{n} \leq b_{n}, \\
    0, & \text{otherwise}, 
  \end{cases}
\end{equation*}
where 
\begin{equation*}
a_{n} = \sum_{j=0}^{t_{n}}(\hat{h}_{n, j} - \hat{p}_{n, j})^{2},
\end{equation*}
and 
\begin{equation*}
b_{n} = \sum_{j=0}^{t_{n}}\hat{p}_{n, j}(\hat{h}_{n, j}^{\backslash[j]} - \hat{p}_{n, j}^{\backslash[j]}) - \sum_{j=0}^{t_{n}}\hat{p}_{n, j}(\hat{h}_{n, j} - \hat{p}_{n, j}),
\end{equation*}
with $\hat{\bm{h}}^{\backslash [j]}_{n} = \hat{\bm{r}}^{\backslash [j]}_{n}$ for the case of stacked rearrangement estimator, and $\hat{\bm{h}}^{\backslash [j]}_{n} = \hat{\bm{g}}^{\backslash [j]}_{n}$ for the case of stacked Grenander estimator, respectively.
\end{theorem}

In the sequel of the paper we always assume that both $\hat{\bm{\phi}}_{n}$ and $\hat{\bm{\phi}}^{\backslash [j]}_{n}$ are constructed with the  leave-one-out  least-squares cross-validation mixture parameter $\hat{\beta}_{n}$.

\section{Theoretical properties of the estimator}\label{est_prpts}
In this section we study theoretical properties of stacked rearrangement and stacked Grenander estimators. First, let us assume that $\bm{p} \in\mathcal{F}^{decr}$, i.e. the underlying p.m.f. is decreasing. Note that from the subadditivity of the norms for $||\hat{\bm{\phi}}_{n} - \bm{p}||_{k}$, with $1 \leq k \leq \infty$, we have
\begin{equation*}
\begin{aligned}
||\hat{\bm{\phi}}_{n} - \bm{p}||_{k} ={} & ||\hat{\beta}_{n} \hat{\bm{h}}_{n} + (1-\hat{\beta}_{n})\hat{\bm{p}}_{n} - \bm{p}||_{k} \leq\\
&\hat{\beta}_{n}||\hat{\bm{h}}_{n} - \bm{p}||_{k} + (1-\hat{\beta}_{n})||\hat{\bm{p}}_{n} - \bm{p}||_{k}.
\end{aligned}
\end{equation*} 
From the error reduction property of the rearrangement and Grenander estimators, i.e. $||\hat{\bm{h}}_{n} - \bm{p}||_{k} \leq ||\hat{\bm{p}}_{n} - \bm{p}||_{k}$,  with $1 \leq k \leq \infty$, cf. Theorem 2.1 in \cite{jankowski2009estimation}, we have
\begin{equation}\label{erpsge}
||\hat{\bm{\phi}}_{n} - \bm{p}||_{k} \leq ||\hat{\bm{p}}_{n} - \bm{p}||_{k}
\end{equation} 
for all  $1 \leq k \leq \infty$. Therefore, in the case of a decreasing true p.m.f. both the stacked rearrangement and stacked Grenander estimators also provide the error reduction. 

Assume that the true p.m.f. is not decreasing. Let $\bm{r} = rear(\bm{p})$ and $\bm{g} = \Pi \big(\bm{p}|\mathcal{F}^{decr}\big)$. Note that $\bm{r} \neq \bm{p}$ nor $\bm{g} \neq \bm{p}$, if $\bm{p} \not\in\mathcal{F}^{decr}$, i.e. the vector $\bm{r}$ is reversed ordered vector $\bm{p}$ and $\bm{g}$ is decreasing vector in $\ell_{2}$ which is closest in $\ell_{2}$-norm to the true p.m.f. $\bm{p}$. 

Then, since the isotonic regression and the rearrangement, viewed as a mapping from $\ell_{2}$ into $\ell_{2}$, are continuous in the case of a finite support, and the empirical estimator is strongly consistent, then
\begin{equation*}
\hat{\bm{r}}_{n}\stackrel{a.s.}{\to} \bm{r}, \, \text{ and } \, \hat{\bm{g}}_{n}\stackrel{a.s.}{\to} \bm{g},
\end{equation*}
pointwise. Note that from the statements $(i)$, $(ii)$ and $(iv)$ of Lemma \ref{propisot} in Appendix it follows that $\hat{\bm{g}}_{n}$ always exists, and it is a probability vector for all $n$. Clearly, the same result holds for the rearrangement estimator $\hat{\bm{r}}_{n}$ for all $n$. The almost sure convergence in $\ell_{k}$-norm, for $1 \leq k \leq \infty$, of $\hat{\bm{r}}_{n}$ and $\hat{\bm{g}}_{n}$ to  $\bm{r}$ and $\bm{g}$, respectively, now follows from Lemma C.2 in the supporting material of \cite{baljanruf}.

\subsection{Consistency} 
First, let us study the leave-one-out versions of the empirical, rearrangement and  Grenander estimators. Recall that
\begin{equation*}
\hat{\bm{p}}^{\backslash [j]}_{n} = \frac{\bm{x} - \bm{\delta}^{[j]}}{n-1}, \, \, \hat{\bm{r}}^{\backslash [j]}_{n} = rear(\hat{\bm{p}} ^{\backslash [j]}_{n})\, \text{ and } \, \hat{\bm{g}}^{\backslash [j]}_{n} = \Pi \big(\hat{\bm{p}} ^{\backslash [j]}_{n}|\mathcal{F}^{decr}\big),
\end{equation*}
for $j$ such that $x_{j} > 0$. 

Let us define vectors $\hat{\bm{\pi}}_{n} \in \ell_{1}$, $\hat{\bm{\rho}}_{n} \in \ell_{1}$, and $\hat{\bm{\gamma}}_{n} \in \ell_{1}$  as

\begin{equation}\label{pigamma}
\begin{aligned}
  \hat{\pi}_{n,j}=\begin{cases}
    \hat{p}^{\backslash [j]}_{n,j}, & \text{if } \, x_{j} > 0, \\
    0, & \text{otherwise},
  \end{cases}\\
      \hat{\rho}_{n,j}=\begin{cases}
    \hat{r}^{\backslash [j]}_{n,j}, & \text{if } \, x_{j} > 0, \\
    0, & \text{otherwise}, 
  \end{cases}\\
  \hat{\gamma}_{n,j}=\begin{cases}
    \hat{g}^{\backslash [j]}_{n,j}, & \text{if } \, x_{j} > 0, \\
    0, & \text{otherwise}.
  \end{cases} 
\end{aligned}
\end{equation}

\begin{lem}\label{pwconv}
The sequences of vectors $\hat{\bm{\pi}}_{n}$,  $\hat{\bm{\rho}}_{n}$ and $\hat{\bm{\gamma}}_{n}$ converge pointwise a.s. to $\bm{p}$, $\bm{r}$, and $\bm{g}$, respectively.
\end{lem}
\prf
The proof is given in Appendix.
\eop

Next, we prove the following important lemma.
\begin{lem}\label{loois}
For the vectors $\hat{\bm{\pi}}_{n}$ we have
\begin{equation*}
\hat{\pi}_{n,j} \leq \hat{p}_{n,j}
\end{equation*}
for all $j$, and for $\hat{\bm{\rho}}_{n}$ and $\hat{\bm{\gamma}}_{n}$ we have
\begin{equation*}
 \hat{\rho}_{n,j} \leq \frac{n}{n-1}\hat{r}_{n,j}  \, \text{ and }  \,  \hat{\gamma}_{n,j} \leq \frac{n}{n-1} \hat{g}_{n,j}
\end{equation*}
for all $j$.
\end{lem}
\prf
The proof is given in Appendix.
\eop

In Lemma C.2 in the supporting material of \cite{baljanruf} it was proved that for probability mass functions the pointwise  convergence and the convergence in $\ell_{k}$ for $1 \leq k \leq \infty$ are all equivalent. Note, in our case the sequences $\hat{\bm{\pi}}_{n}$, $\hat{\bm{\rho}}_{n}$ and $\hat{\bm{\gamma}}_{n}$ are not probability vectors. Nevertheless, as we prove below, all $\hat{\bm{\pi}}_{n}$, $\hat{\bm{\rho}}_{n}$ and $\hat{\bm{\gamma}}_{n}$ converge a.s. to $\bm{p}$, $\bm{r}$ and $\bm{g}$, respectively, in $\ell_{k}$-norm for $1 \leq k \leq \infty$.
\begin{theorem}\label{asconvpigm}
For the vectors $\hat{\bm{\pi}}_{n}$, $\hat{\bm{\rho}}_{n}$ and $\hat{\bm{\gamma}}_{n}$ we have
\begin{equation*}
\hat{\bm{\pi}}_{n}\stackrel{a.s.}{\to} \bm{p},
\end{equation*}
\begin{equation*}
\hat{\bm{\rho}}_{n}\stackrel{a.s.}{\to} \bm{r},
\end{equation*}
and 
\begin{equation*}
\hat{\bm{\gamma}}_{n}\stackrel{a.s.}{\to} \bm{g}
\end{equation*}
in $\ell_{k}$-norm for $1 \leq k \leq \infty$.
\end{theorem}
\prf
The proof starts in a similar way as the one for Lemma C.2 in \cite{baljanruf}. Let us, first, study the case of $\hat{\bm{\pi}}_{n}$. Fix some $\varepsilon>0$. Then, we can choose $K$ such that 
\begin{equation*}
\sum_{j\leq K}p_{j} \geq 1 - \frac{\varepsilon}{4}.
\end{equation*}

Since both $\bm{\pi}_{n}$ and the empirical estimator $\bm{p}_{n}$ converge to $\mathbf{p}$ pointwise, then there exists random $n_{0}$ such that for all $n\geq n_{0}$
\begin{equation*}
\sup_{j\leq K}|\hat{p}_{n,j} - p_{j}| \leq \frac{\varepsilon}{4(K+1)},
\end{equation*} 
\begin{equation*}
\sup_{j\leq K}|\hat{\pi}_{n,j} - p_{j}| \leq \frac{\varepsilon}{4(K+1)},
\end{equation*}
almost surely.

This implies that for all $n\geq n_{0}$ we have $\sum_{j\leq K}\hat{p}_{n,j} \geq 1 - \frac{\varepsilon}{2}$ and $\sum_{j\leq K}|\hat{\pi}_{n,j} - p_{j}| \leq \frac{\varepsilon}{4}$, almost surely.

Next, for any $n$ 
\begin{equation*}
\sum_{j=0}^{\infty}|\hat{\pi}_{n,j} - p_{j}| = \sum_{j\leq K}|\hat{\pi}_{n,j} - p_{j}| + \sum_{j>K}|\hat{\pi}_{n,j} - p_{j}| \leq \sum_{j\leq K}|\hat{\pi}_{n,j} - p_{j}| + \sum_{j>K} \hat{\pi}_{n,j} + \sum_{j>K}p_{j}.
\end{equation*}
Furthermore, $\sum_{j>K}\hat{\pi}_{n,j} \leq  \sum_{j>K}\hat{p}_{n,j}$ since $0<\hat{\pi}_{n,j} \leq \hat{p}_{n,j}$. Then, for all $n>n_{0}$ we have proved that
\begin{equation*}
\sum_{j=0}^{\infty}|\hat{\pi}_{n,j} - p_{j}| \leq \frac{\varepsilon}{4} + \frac{\varepsilon}{2} + \frac{\varepsilon}{4} = \varepsilon,
\end{equation*}
almost surely. This means that for any $\varepsilon>0$ there exists random $n_{0}$, such that for all $n>n_{0}$ 
\begin{equation*}
||\hat{\bm{\pi}}_{n} - \textbf{p}||_{1} \leq  \varepsilon,
\end{equation*}
almost surely.  

Furthermore, since $\ell_{1} \subset \ell_{k}$, for all $k>1$, then a.s. convergence holds in $\ell_{k}$, for all $1\leq k\leq\infty$. 

Let us prove the convergence for $\hat{\bm{\gamma}}_{n}$. First, from Lemma \ref{loois} it follows that 
\begin{equation*}
\frac{n-1}{n}\hat{\gamma}_{n,j} \leq \hat{g}_{n,j}.
\end{equation*}
Then, since both $\frac{n-1}{n}\hat{\bm{\gamma}}_{n}$ and $\hat{\bm{g}}_{n}$ converge to $\bm{g}$ a.s., we can use the same approach as for $\hat{\bm{\pi}}$ above, and prove that  
\begin{equation*}
\frac{n-1}{n}\hat{\bm{\gamma}}_{n} \stackrel{a.s.}{\to} \bm{g},
\end{equation*}
in $\ell_{k}$, for $1\leq k\leq\infty$, which means that 
\begin{equation*}
\hat{\bm{\gamma}}_{n} \stackrel{a.s.}{\to} \bm{g},
\end{equation*} 
in $\ell_{k}$, for $1\leq k\leq\infty$.

Now, using the result of Lemma \ref{loois}, we can prove the result for $\hat{\bm{\rho}}_{n}$ in the same way as we did for $\hat{\bm{\gamma}}_{n}$.
\eop
Now  we can summarize the above results in the following theorem.
\begin{theorem}\label{conssgrestim}
For any underlying distribution $\bm{p}$, both the stacked rearrangement and stacked Grenander estimators are strongly consistent:
\begin{equation*}
\hat{\bm{\phi}}_{n}\stackrel{a.s.}{\to} \bm{p}
\end{equation*}
in $\ell_{k}$-norm for $1 \leq k \leq \infty$.
\end{theorem}
\prf 
Firs, let us assume that $\bm{p}$ is decreasing. Then the result of the theorem
follows from the strong consistency of $\hat{\bm{g}}_{n}$, $\hat{\bm{r}}_{n}$ and $\hat{\bm{p}}_{n}$.

Next, assume that $\bm{p}$ is not decreasing. From Theorem \ref{asconvpigm} it follows that for the case of stacked rearrangement estimator we have
\begin{equation*}
a_{n}\stackrel{a.s.}{\to} ||\bm{r} - \bm{p}||_{2}^{2},
\end{equation*}
and
\begin{equation*}
b_{n}\stackrel{a.s.}{\to}  \langle \bm{p},(\bm{r} - \bm{p}) \rangle - \langle \bm{p},(\bm{r} - \bm{p}) \rangle = 0,
\end{equation*}
and for the case of stacked Grenander estimator we have
\begin{equation*}
a_{n}\stackrel{a.s.}{\to} ||\bm{g} - \bm{p}||_{2}^{2},
\end{equation*}
and
\begin{equation*}
b_{n}\stackrel{a.s.}{\to}  \langle \bm{p},(\bm{g} - \bm{p}) \rangle - \langle \bm{p},(\bm{g} - \bm{p}) \rangle = 0.
\end{equation*}
Therefore,
\begin{equation*}
\hat{\beta}_{n}\stackrel{a.s.}{\to} 0.
\end{equation*}
Next, since
\begin{equation*}
||\hat{\bm{\phi}}_{n} - \bm{p}||_{k} \leq \hat{\beta}_{n}||\hat{\bm{h}}_{n} - \bm{p}||_{k} + (1-\hat{\beta}_{n})||\hat{\bm{p}}_{n} - \bm{p}||_{k}
\end{equation*}  
for all $1\leq k \leq \infty$, it follows
\begin{equation*}
\hat{\bm{\phi}}_{n}\stackrel{a.s.}{\to} \bm{p}
\end{equation*}
in $\ell_{k}$-norm for $1 \leq k \leq \infty$.
\eop

\subsection{Rate of convergence}
In this section we study the rate of convergence of stacked estimator. In the case of bounded support the $\sqrt{n}$-rate of convergence follows from pointwise convergence of the vectors $\hat{\bm{\pi}}_{n}$, $\hat{\bm{\rho}}_{n}$ and $\hat{\bm{\gamma}}_{n}$. In this work we assume that the support can be infinite.
\begin{theorem}\label{rtconvsge}
Stacked rearrangement and Grenander estimators have $\sqrt{n}$-rate of convergence for any underlying p.m.f. $\bm{p}$:
\begin{equation*}
\sqrt{n}||\hat{\bm{\phi}}_{n} - \bm{p}||_{k} = O_{p}(1)
\end{equation*}
for $1 < k \leq \infty$. Next, if $\sum_{j=0}^{\infty}\sqrt{p}_{j} < \infty$, then
\begin{equation*}
\sqrt{n}||\hat{\bm{\phi}}_{n} - \bm{p}||_{1} = O_{p}(1).
\end{equation*}
\end{theorem}

\prf
Assume that $\bm{p}$ is decreasing. Then the result follows from (\ref{erpsge}) and Corollaries 4.1 and 4.2 in \cite{jankowski2009estimation}. 

Next, assume that $\bm{p}$ is not decreasing. Let us, first, prove the case of stacked Grenander estimator. Recall that
\begin{equation*}\label{}
  \beta_{n}=\begin{cases}
    \frac{b_{n}}{a_{n}}, & \text{if } \, a_{n} \neq 0 \text{ and } 0 \leq b_{n} \leq a_{n}, \\
    1, & \text{if } \, 0 < a_{n} \leq b_{n}, \\
    0, & \text{otherwise}, 
  \end{cases}
\end{equation*}
where 
\begin{equation*}
a_{n} = \sum_{j=0}^{t_{n}}(\hat{g}_{n, j} - \hat{p}_{n, j})^{2},
\end{equation*}
and in the notation introduced in \ref{pigamma}, we can write $b_{n}$ as
\begin{equation*}
b_{n} = \sum_{j=0}^{t_{n}}\hat{p}_{n, j}(\hat{\gamma}_{n, j} - \hat{\pi}_{n, j}) - \sum_{j=0}^{t_{n}}\hat{p}_{n, j}(\hat{g}_{n, j} - \hat{p}_{n, j}).
\end{equation*}

First, as we proved in Theorem \ref{conssgrestim} 
\begin{equation}\label{anconv}
a_{n}\stackrel{a.s.}{\to} ||\bm{g} - \bm{p}||_{2}^{2}> 0.
\end{equation} 

Second, note that from Lemma \ref{loois} it follows that for all $n$ we have
\begin{equation*}
\begin{aligned}
b_{n}  ={} & \sum_{j=0}^{t_{n}}\hat{p}_{n, j}(\hat{\gamma}_{n, j} - \hat{g}_{n, j}) + \sum_{j=0}^{t_{n}}\hat{p}_{n, j}(\hat{p}_{n, j} - \hat{\pi}_{n, j}) 
\leq\\
&\frac{n}{n-1}\sum_{j=0}^{t_{n}}\hat{p}_{n, j}\hat{g}_{n, j} - \sum_{j=0}^{t_{n}}\hat{p}_{n, j} \hat{g}_{n, j} + \sum_{j=0}^{t_{n}}\hat{p}_{n, j}(\hat{p}_{n, j} - \hat{\pi}_{n, j}).
\end{aligned}
\end{equation*} 
Next, 
\begin{equation*}
\begin{aligned}
\frac{n}{n-1}\sum_{j=0}^{t_{n}}\hat{p}_{n, j}\hat{g}_{n, j} - \sum_{j=0}^{t_{n}}\hat{p}_{n, j} \hat{g}_{n, j}  = \frac{\sum_{j=0}^{t_{n}}\hat{p}_{n, j}\hat{g}_{n, j}}{n-1}.
\end{aligned}
\end{equation*} 
Recall that 
\begin{equation*}\label{}
  \hat{\pi}_{n,j} =\begin{cases}
    \frac{x_{j} - 1}{n-1} = \frac{n}{n-1}\hat{p}_{n,j} - \frac{1}{n-1}, & \text{if } \, x_{j}\neq 0, \\
    0, & \text{otherwise},
  \end{cases}
\end{equation*}
which leads to
\begin{equation*}
\sum_{j=0}^{t_{n}}\hat{p}_{n, j}(\hat{p}_{n, j} - \hat{\pi}_{n, j}) = \sum_{j=0}^{t_{n}}\hat{p}_{n, j}(\hat{p}_{n, j} - \hat{\pi}_{n, j}) = \frac{1-\sum_{j=0}^{t_{n}}\hat{p}_{n, j}^{2}}{n-1}.
\end{equation*} 
Therefore, the upper bound for $b_{n}$ is given by
\begin{equation*}\label{}
b_{n} \leq \frac{\sum_{j=0}^{t_{n}}\hat{p}_{n, j}\hat{g}_{n, j}}{n-1} + \frac{1-\sum_{j=0}^{t_{n}}\hat{p}_{n, j}^{2}}{n-1},
\end{equation*}
and, consequently,
\begin{equation}\label{bnub}
\sqrt{n}b_{n} \stackrel{a.s.}{\to} 0,
\end{equation}
since both sequences $\sum_{j=0}^{t_{n}}\hat{p}_{n, j}\hat{g}_{n, j}$ and $\hat{p}_{n, j}^{2}$ are bounded. 

Next, since $\beta_{n} \geq 0$, from (\ref{anconv}) and (\ref{bnub}) it follows that
\begin{equation}\label{convlbd}
\sqrt{n} \hat{\beta}_{n} \stackrel{a.s.}{\to} 0. 
\end{equation}
Then, from (\ref{convlbd}) for any $\bm{p}$ and all $1\leq k \leq \infty$ the following holds
\begin{equation*}
\hat{\beta}_{n}\sqrt{n}||\hat{\bm{g}}_{n} - \bm{p}||_{k} \stackrel{a.s.}{\to} 0, 
\end{equation*} 
for all $1\leq k \leq \infty$. Further, as it follows from Corollary 4.2 in \cite{jankowski2009estimation}, if $\sum_{j=0}^{\infty}\sqrt{p}_{j} < \infty$,
then 
\begin{equation*}
\sqrt{n}||\hat{\bm{p}}_{n} - \bm{p}||_{1} = O_{p}(1).
\end{equation*} 

Therefore, for all $2 \leq k \leq \infty$ and all $\bm{p}$ we have
\begin{equation*}
(1-\hat{\beta}_{n})\sqrt{n}||\hat{\bm{p}}_{n} - \bm{p}||_{k} = O_{p}(1),
\end{equation*} 
and,  if $\sum_{j=0}^{\infty}\sqrt{p}_{j} < \infty$, then we have
\begin{equation*}
(1-\hat{\beta}_{n})\sqrt{n}||\hat{\bm{p}}_{n} - \bm{p}||_{1} = O_{p}(1).
\end{equation*} 

Finally, recall that
\begin{equation*}
\sqrt{n}||\hat{\bm{\phi}}_{n} - \bm{p}||_{k} \leq \hat{\beta}_{n}\sqrt{n}||\hat{\bm{g}}_{n} - \bm{p}||_{k} + (1-\hat{\beta}_{n})\sqrt{n}||\hat{\bm{p}}_{n} - \bm{p}||_{k},
\end{equation*} 
which finishes the prove of theorem for the case of Grenander estimator.

Similarly, using the results of Lemma \ref{loois}, for the case of stacked rearrangement estimator we can show that 
\begin{equation*}\label{}
b_{n} \leq \frac{\sum_{j=0}^{t_{n}}\hat{p}_{n, j}\hat{r}_{n, j}}{n-1} + \frac{1-\sum_{j=0}^{t_{n}}\hat{p}_{n, j}^{2}}{n-1},
\end{equation*}
for all $n$. Then, the rest of the proof is the same as for Grenander estimator with $\hat{\bm{g}}_{n}$ and $\bm{g}$ suitably changed to $\hat{\bm{r}}_{n}$ and $\bm{r}$, respectively.
\eop

\subsection{Asymptotic distribution and global confidence band}
In this section we study the asymptotic distribution of stacked rearrangement and Grenander estimators and discuss calculation of global confidence band for $\bm{p}$. The limit distribution of rearrangement and Grenanader estimators were obtained in \cite{jankowski2009estimation}. The asymptotic distribution of stacked Grenander estimator for the case when true p.m.f. $\bm{p}$ is either not decreasing with a countable support or strictly decreasing with a finite support is given in the next theorem.
\begin{theorem}\label{asymdist}
Assume that $\bm{p}$ is either not decreasing with a countable support or strictly decreasing with a finite support. Then stacked rearrangement and Grenander estimators are asymptotically normal
\begin{equation*}
\sqrt{n}(\hat{\bm{\phi}}_{n} -\bm{p}) \stackrel{d}{\to} \bm{Y}_{\bm{0},C},
\end{equation*}
in $\ell_{2}$, where $\bm{Y}_{\bm{0},C}$ is a Gaussian process in $\ell_{2}$ with mean zero and the covariance operator $C$ such that $\langle C \bm{e}_{i}, \bm{e}_{i'}\rangle = p_{i}\delta_{i,i'} - p_{i}p_{i'}$, with $\bm{e}_{i} \in \ell_{2}$ the orthonormal basis in $\ell_{2}$ such that in a vector $\bm{e}_{i}$ all elements are equal to zero but the one with the index $i$ is equal to $1$, and $\delta_{i,j}=1$, if $i=j$ and $0$ otherwise, cf. \cite{jankowski2009estimation}.
\end{theorem}
\prf
The proof is given in Appendix.
\eop

\begin{figure}[!htbp] 
  \begin{subfigure}{3.9cm}
    \centering\includegraphics[scale=0.24]{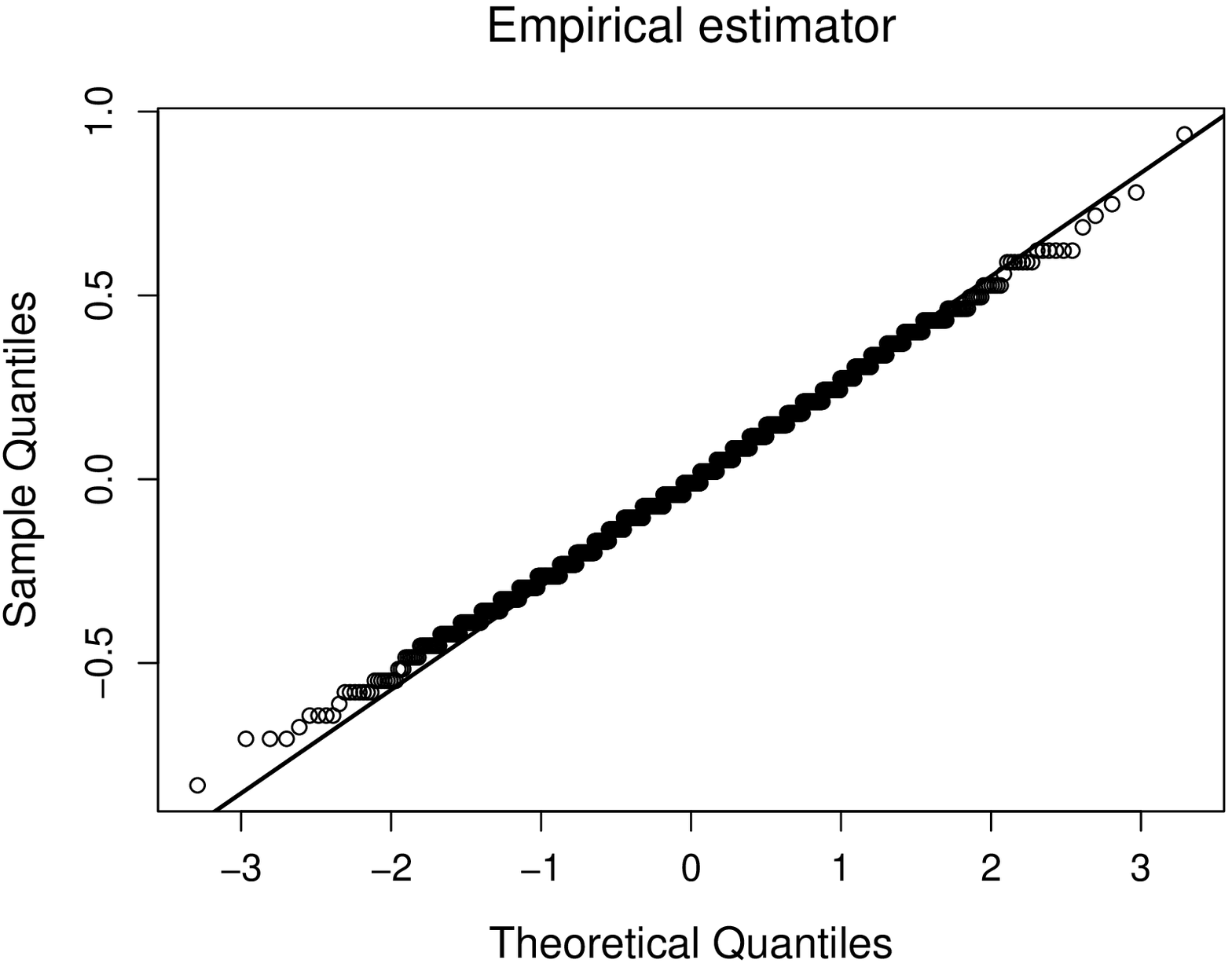}
  \end{subfigure}
  \begin{subfigure}{3.9cm}
    \centering\includegraphics[scale=0.24]{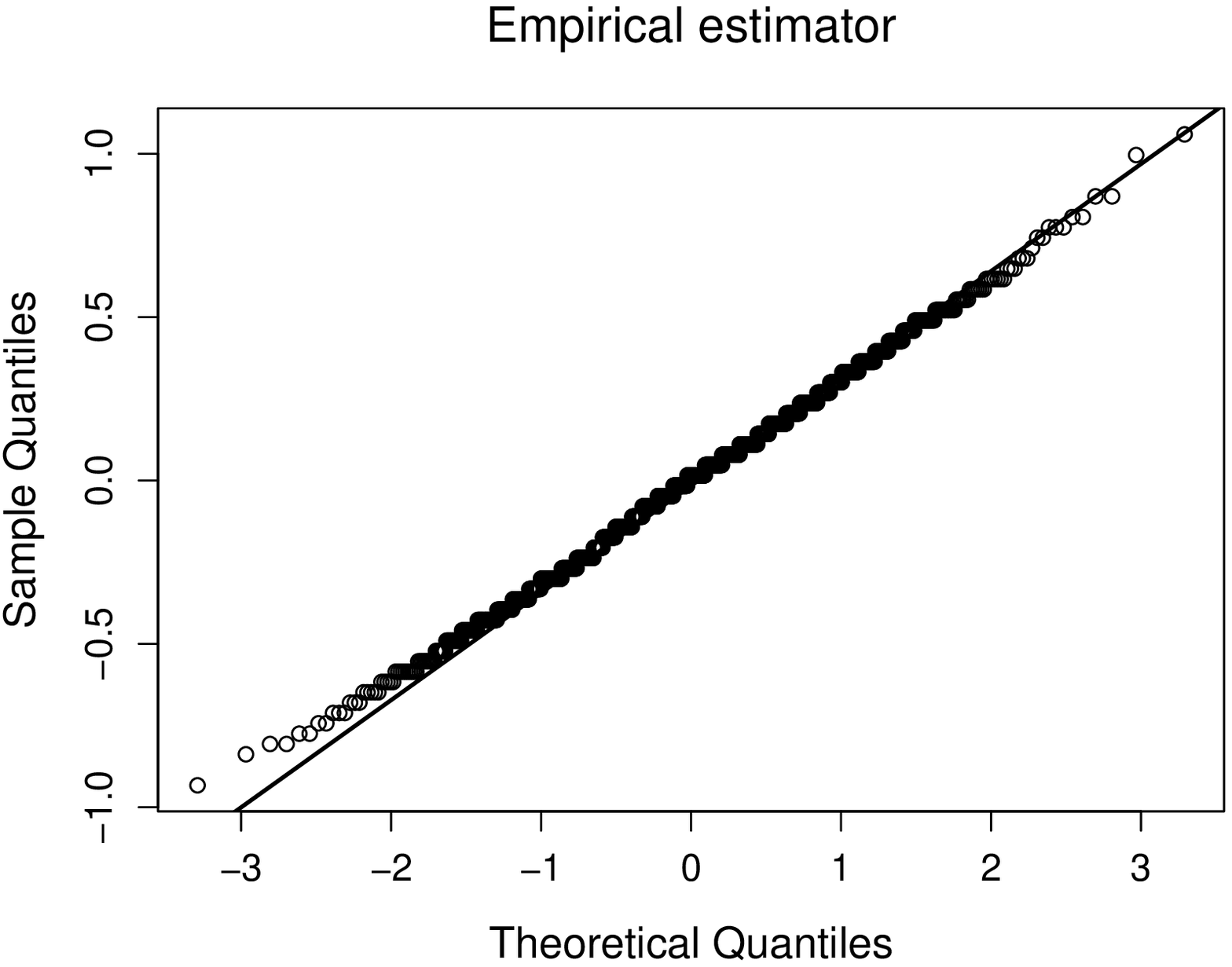}
  \end{subfigure}
  \begin{subfigure}{3.9cm}
    \centering\includegraphics[scale=0.24]{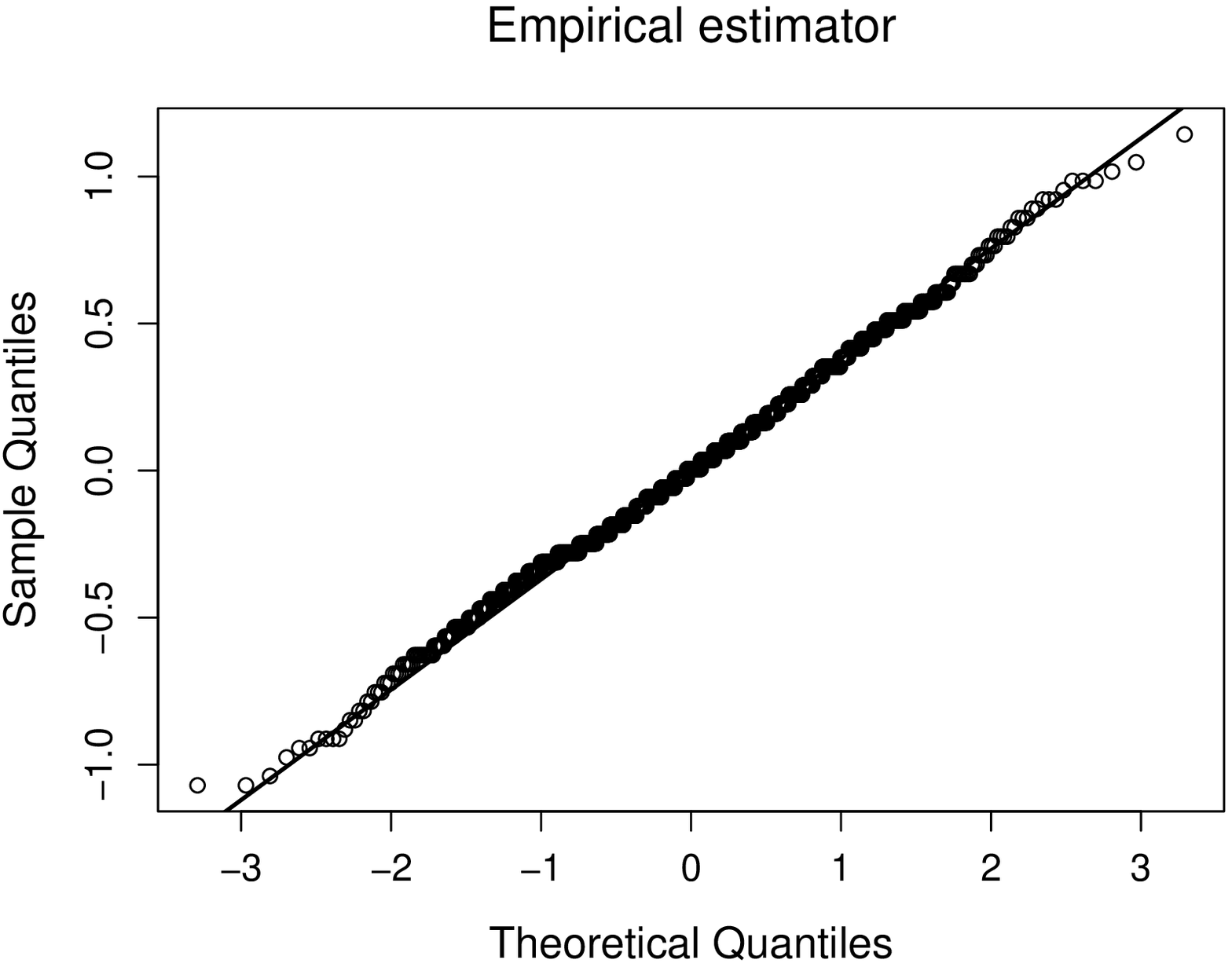}
  \end{subfigure}
  
    \begin{subfigure}{3.9cm}
    \centering\includegraphics[scale=0.24]{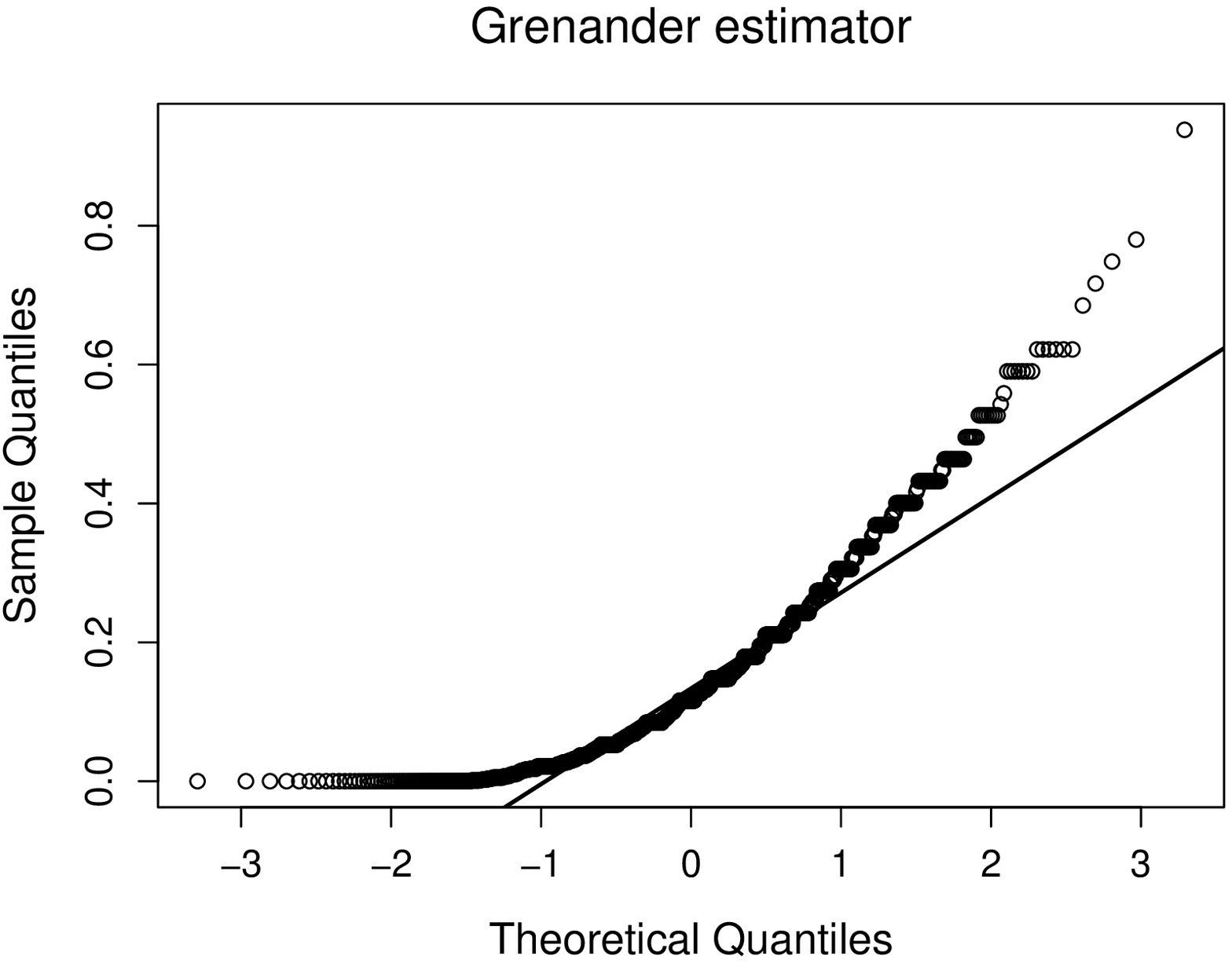}
  \end{subfigure}
  \begin{subfigure}{3.9cm}
    \centering\includegraphics[scale=0.24]{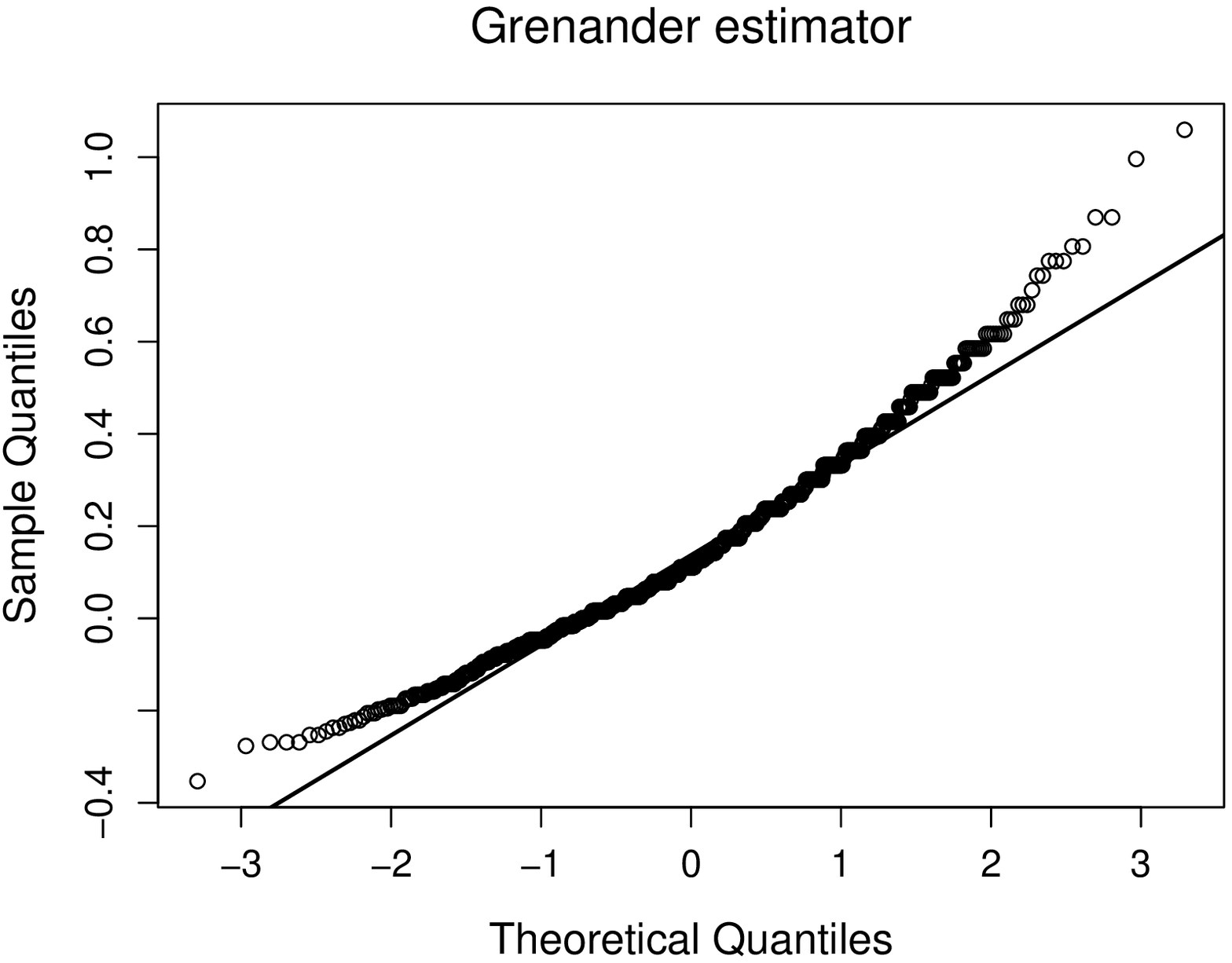}
  \end{subfigure}
  \begin{subfigure}{3.9cm}
    \centering\includegraphics[scale=0.24]{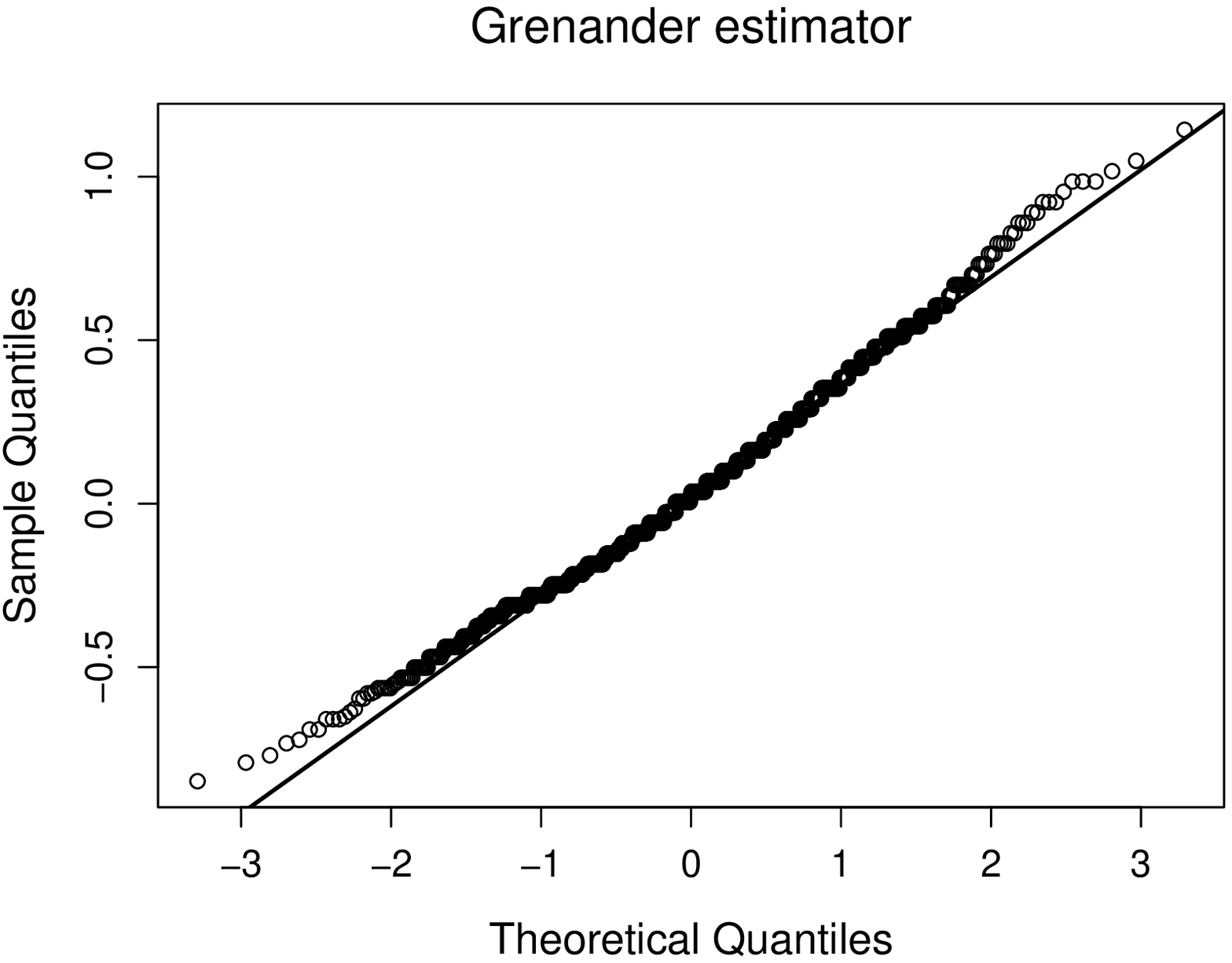}
  \end{subfigure}

    \begin{subfigure}{3.9cm}
    \centering\includegraphics[scale=0.24]{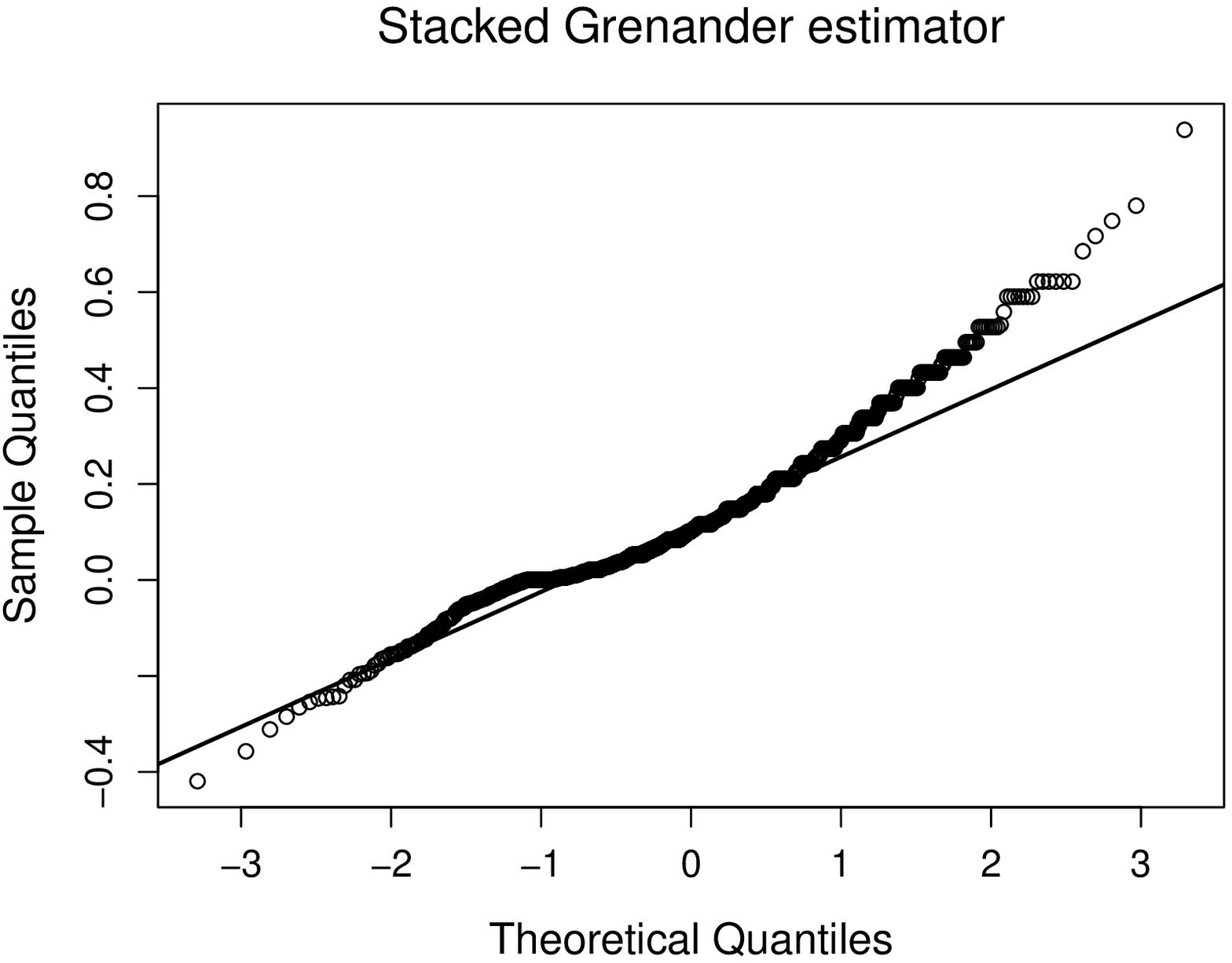}
  \end{subfigure}
  \begin{subfigure}{3.9cm}
    \centering\includegraphics[scale=0.24]{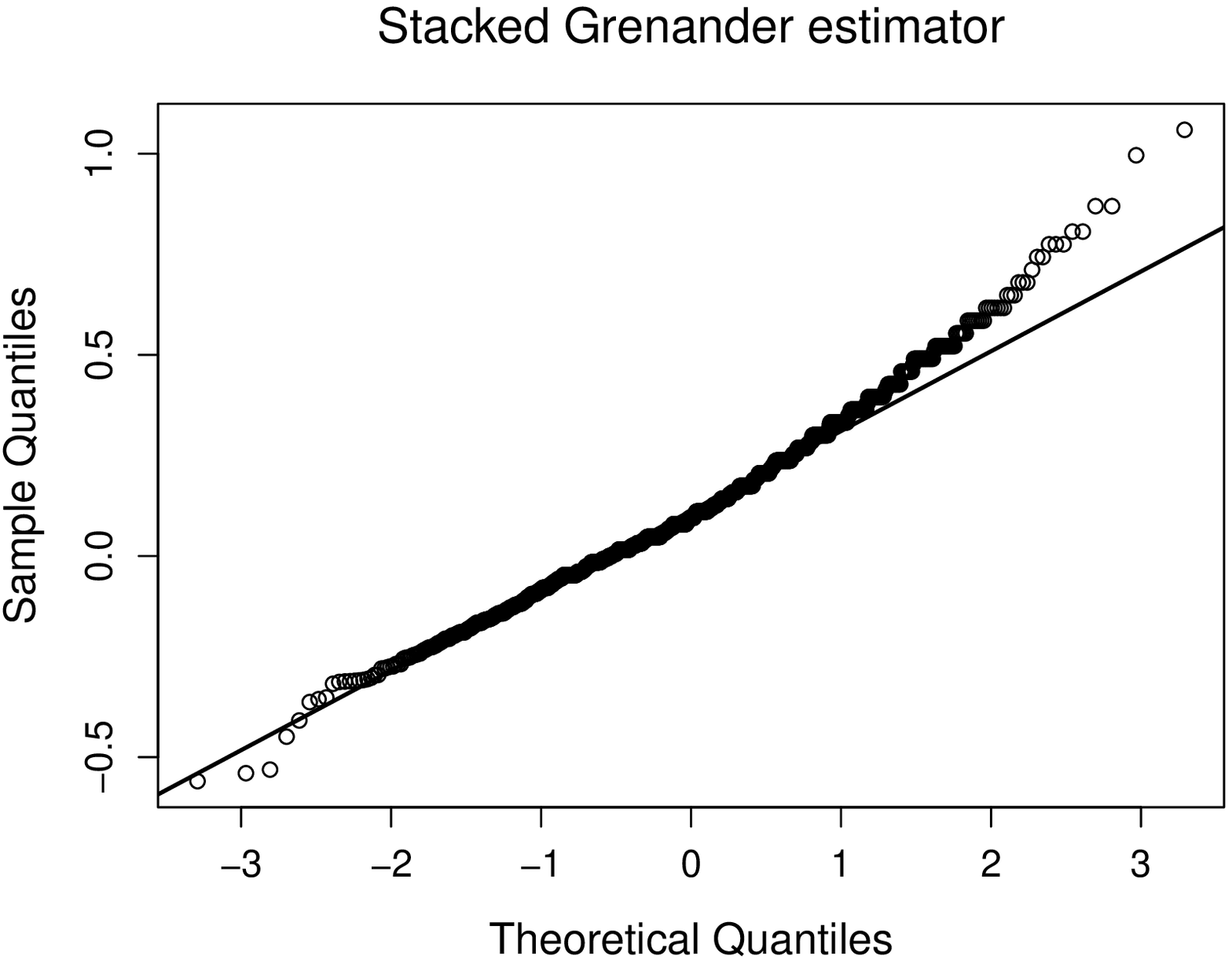}
  \end{subfigure}
  \begin{subfigure}{3.9cm}
    \centering\includegraphics[scale=0.24]{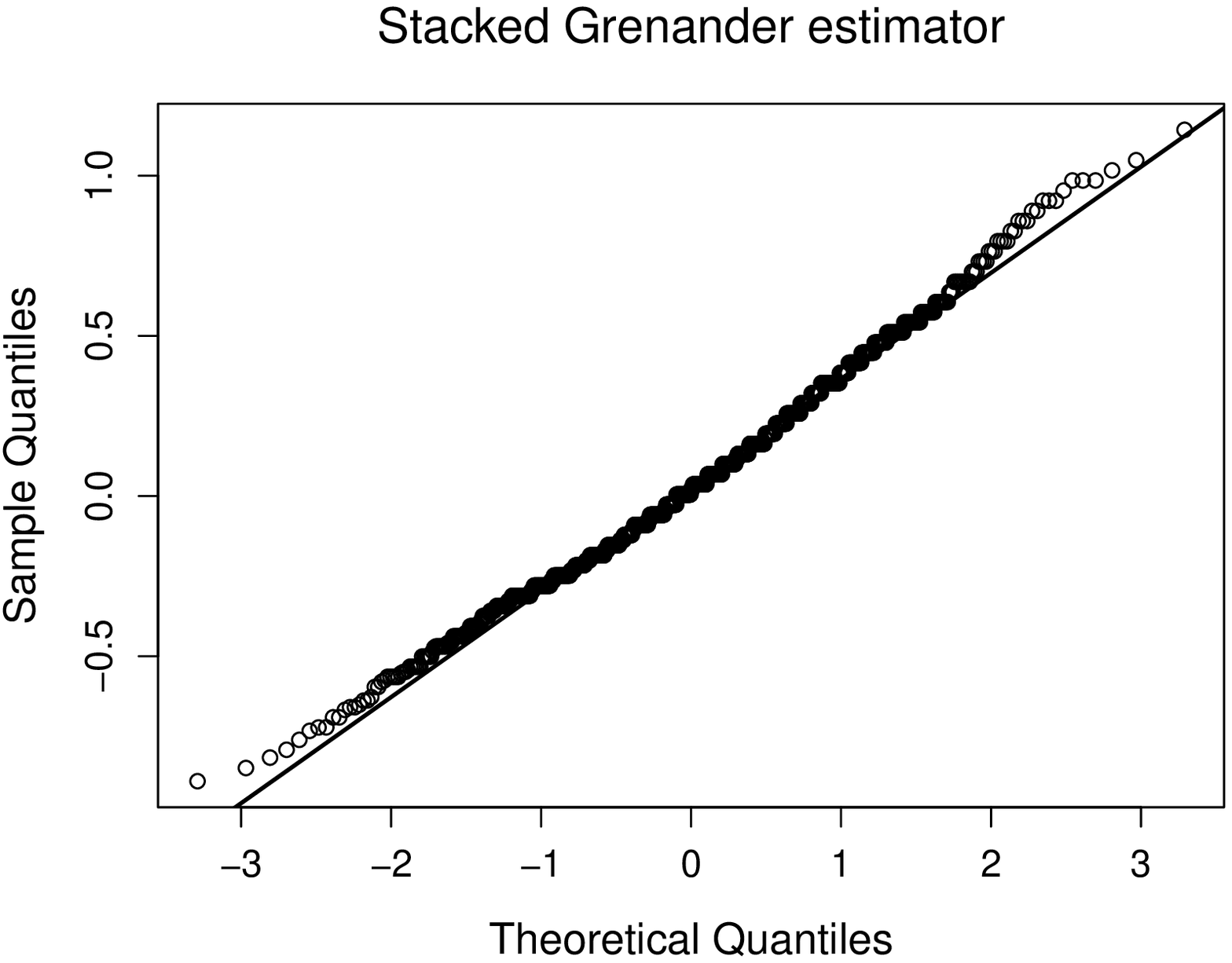}
  \end{subfigure}
  
      \begin{subfigure}{3.9cm}
    \centering\includegraphics[scale=0.24]{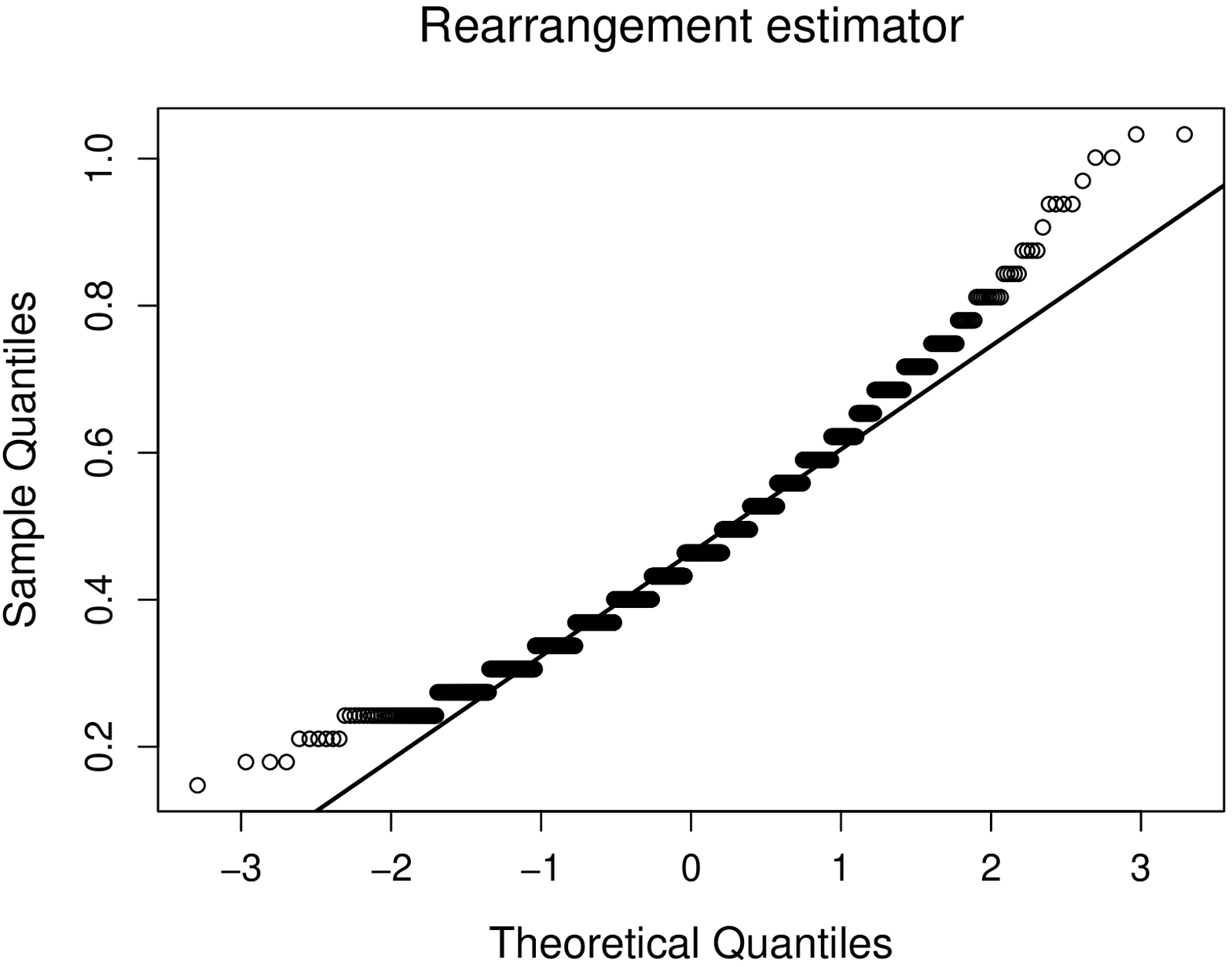}
  \end{subfigure}
  \begin{subfigure}{3.9cm}
    \centering\includegraphics[scale=0.24]{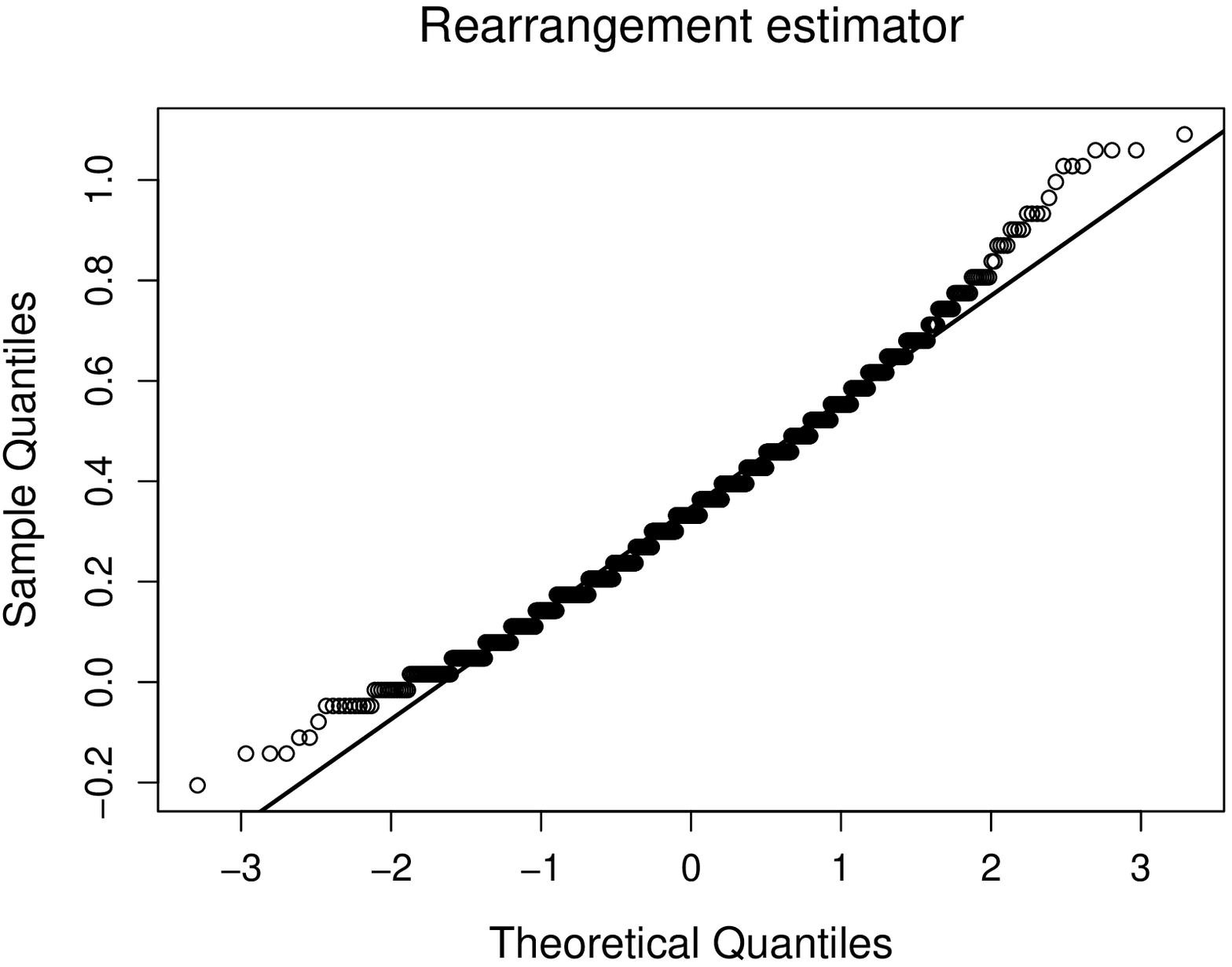}
  \end{subfigure}
  \begin{subfigure}{3.9cm}
    \centering\includegraphics[scale=0.24]{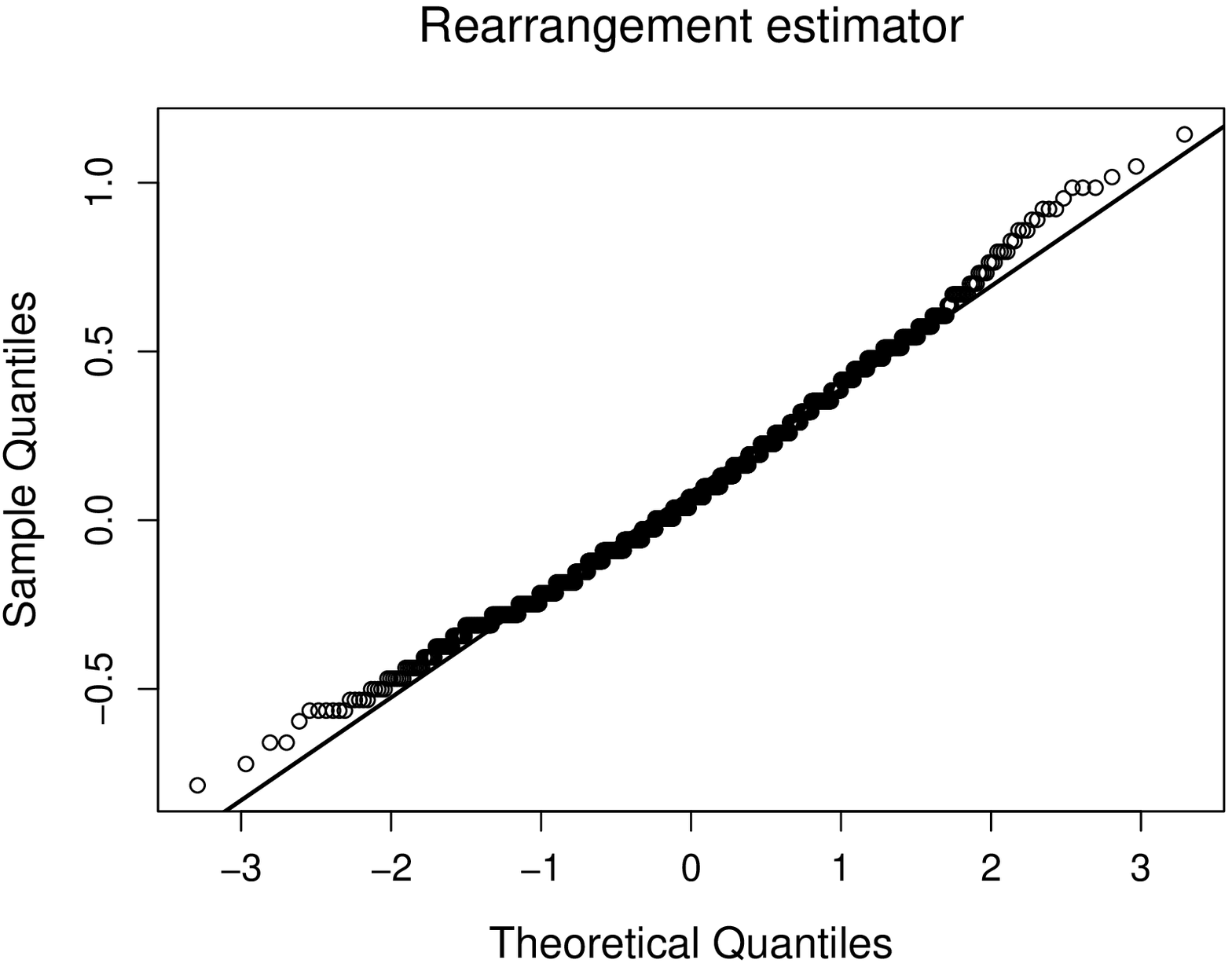}
  \end{subfigure}
  
      \begin{subfigure}{3.9cm}
    \centering\includegraphics[scale=0.24]{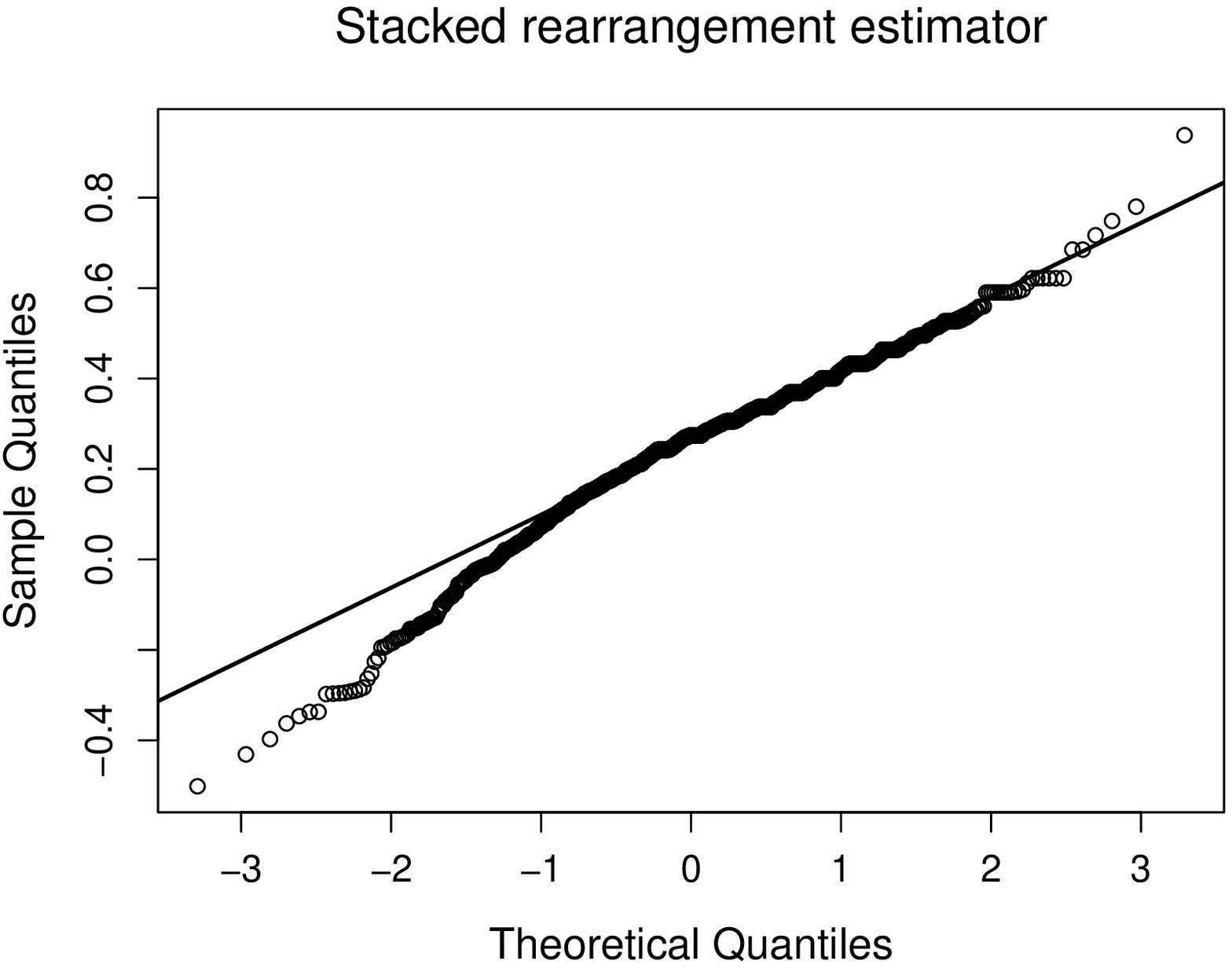}
  \end{subfigure}
  \begin{subfigure}{3.9cm}
    \centering\includegraphics[scale=0.24]{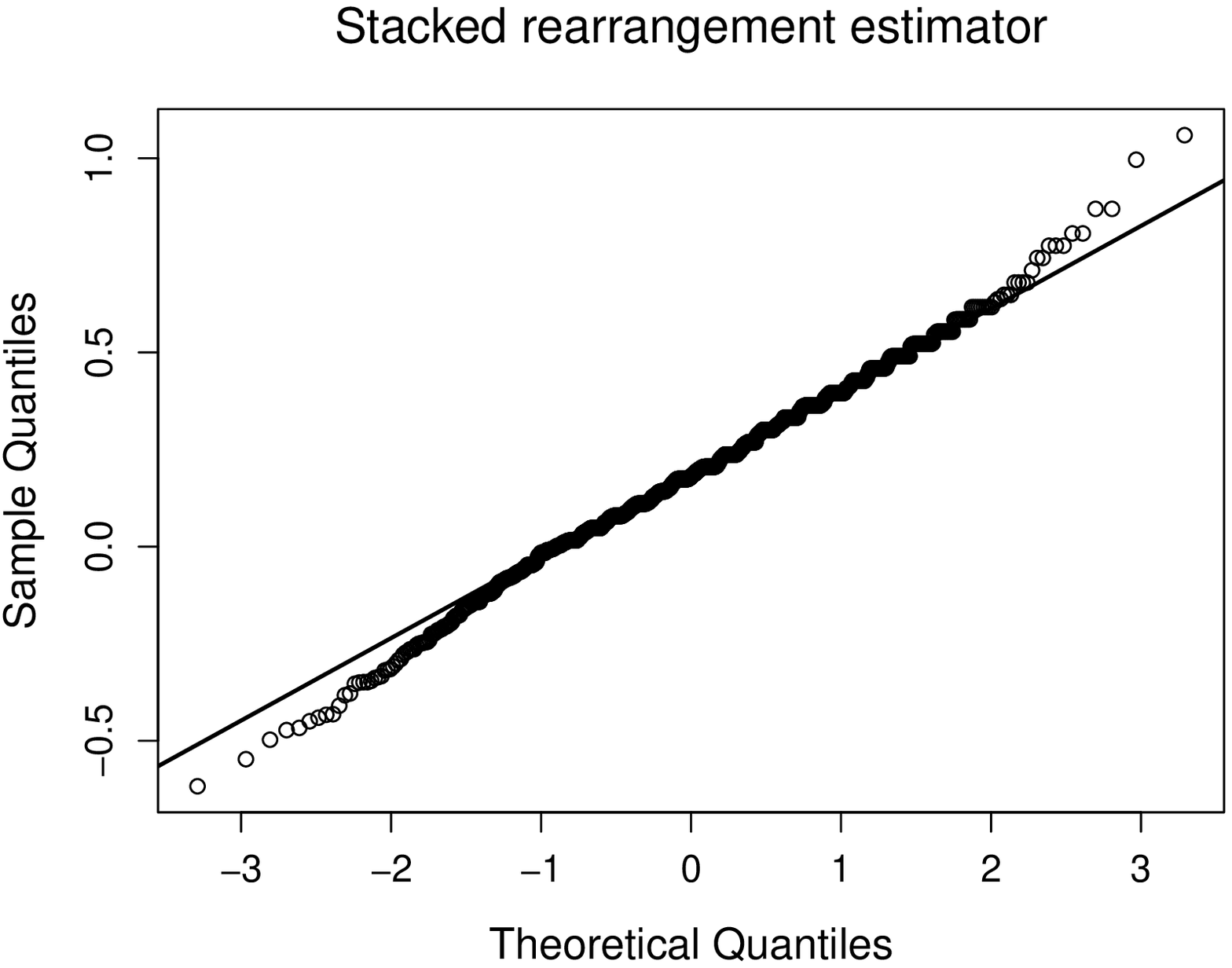}
  \end{subfigure}
  \begin{subfigure}{3.9cm}
    \centering\includegraphics[scale=0.24]{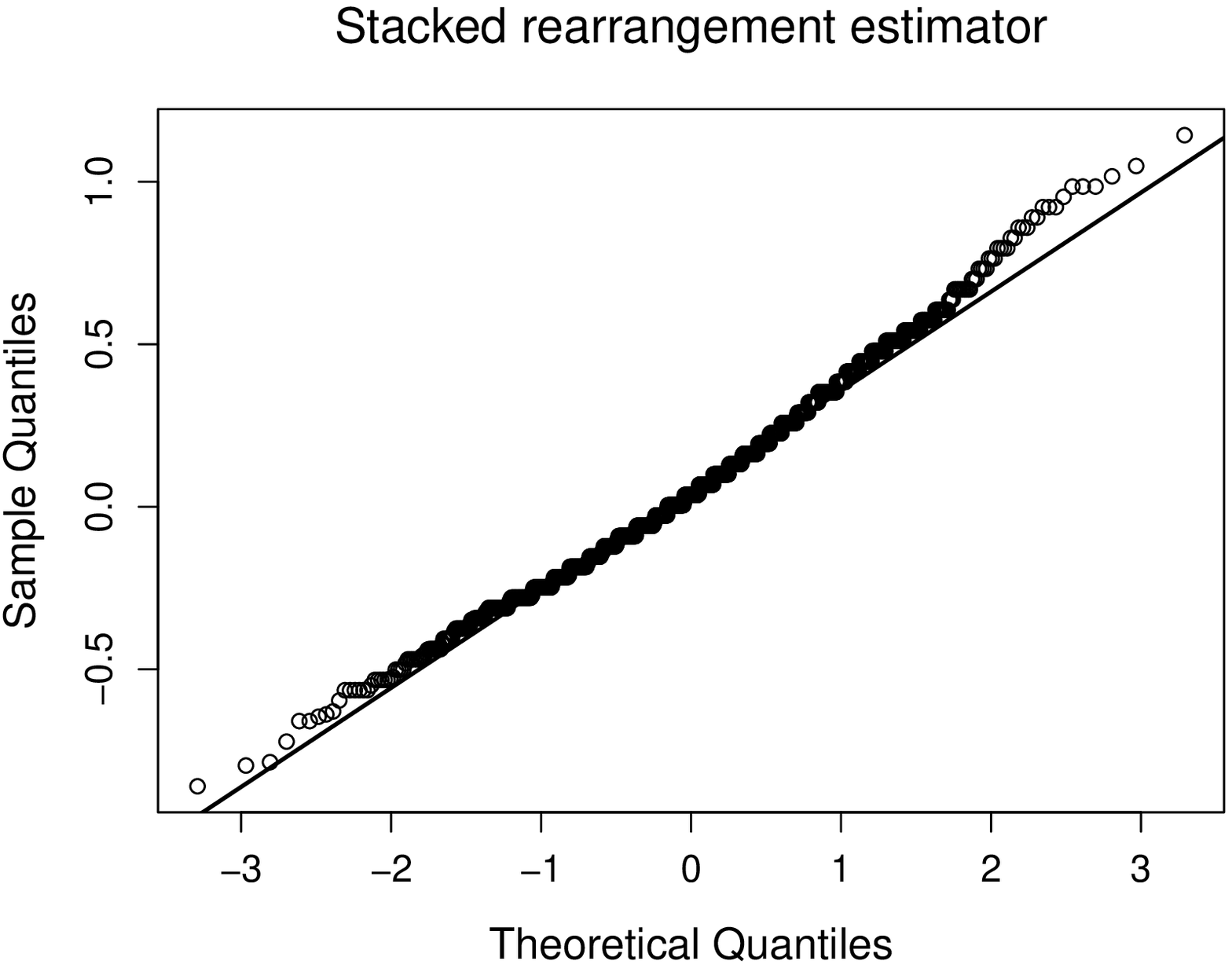}
  \end{subfigure}

   \caption{Standard normal QQ-plots of 1000 samples of $\sqrt{n}(\hat{p}_{n,1} - p_{1})$, $\sqrt{n}(\hat{g}_{n,1} - p_{1})$, $\sqrt{n}(\hat{r}_{n,1} - p_{1})$ and $\sqrt{n}(\hat{\phi}_{n,1} - p_{1})$ for both $\hat{\bm{h}}_{n} = \bm{\hat{g}}_{n}$ and $\hat{\bm{h}}_{n} = \bm{\hat{r}}_{n}$, with $n=1000$ for uniform distribution (left), decreasing distribution (middle) and strictly decreasing distribution (right).}\label{qq-plots}
\end{figure}

For the case of a general decreasing underlying p.m.f. with some constant regions the limit distribution of the stacked estimator remains an open problem. Figure \ref{qq-plots} illustrates the difference of the asymptotic distributions of the empirical estimator, monotonically constrained estimators and the stacked estimators. Let $U(s)$ denote the uniform distribution over $\{0, \dots, s \}$ and $T^{d}(s)$ be strictly decreasing triangular function with the support $\{0, \dots, s \}$ (for the definition of triangular function see e.g. \cite{durot2013}). Figure \ref{qq-plots} shows standard normal QQ-plots of 1000 samples of $\sqrt{n}(\hat{p}_{n,1} - p_{1})$, $\sqrt{n}(\hat{g}_{n,1} - p_{1})$, $\sqrt{n}(\hat{r}_{n,1} - p_{1})$ and $\sqrt{n}(\hat{\phi}_{n,1} - p_{1})$ for both $\hat{\bm{h}}_{n} = \bm{\hat{g}}_{n}$ and $\hat{\bm{h}}_{n} = \bm{\hat{r}}_{n}$, with $n=1000$ for the following distributions:
\begin{enumerate}[label=(\alph*)]
\item (left) $\bm{p} = U(11)$,
\item (middle) $\bm{p} = 0.15 U(3) + 0.1 U(7) + 0.75 U(11)$,
\item (right) $\bm{p} = T^{d}(11)$.
\end{enumerate}
 
From Figure \ref{qq-plots} we can conclude that, first, in the case of a decreasing p.m.f. the distributions of stacked estimators asymptotically are  not equivalent to the distribution of the empirical estimator, and, second, stacked estimators and constrained estimators have different asymptotic distribution if the underlying p.m.f. has constant regions.

For the process $\bm{Y}_{\bm{0},C}$ defined in Theorem \ref{asymdist} let $q_{\alpha}$ denote the $\alpha$-quantile of its $\ell_{\infty}$-norm, i.e.
\begin{equation*}
\mathbb{P}[||\bm{Y}_{\bm{0},C}||_{\infty} > q_{\alpha}] = \alpha.
\end{equation*} 
Then, if $\bm{p}$ is not decreasing or strictly decreasing, from Theorem \ref{asymdist} for stacked estimator we have
\begin{equation*}
\lim_{n}\mathbb{P}[\sqrt{n}||\hat{\bm{\phi}}_{n} -\bm{p}||_{\infty} \leq q_{\alpha}] =  1- \alpha.
\end{equation*} 

Next, note that in the case of a decreasing p.m.f. $\bm{p}$ from (\ref{erpsge}) it follows
\begin{equation*}
\mathbb{P}[\sqrt{n}||\hat{\bm{\phi}}_{n} -\bm{p}||_{\infty} \leq q_{\alpha}] \geq \mathbb{P}[\sqrt{n}||\hat{\bm{p}}_{n} -\bm{p}||_{\infty} \leq q_{\alpha}]
\end{equation*} 
for all $n$. Therefore, in the case of a decreasing $\bm{p}$ we have
\begin{equation*}
\liminf_{n}\mathbb{P}[\sqrt{n}||\hat{\bm{\phi}}_{n} -\bm{p}||_{\infty} \leq q_{\alpha}] \geq 1-  \alpha.
\end{equation*} 

In the same way as in \cite{baljan}, to estimate $q_{\alpha}$ we can use the stacked estimator $\hat{\bm{\phi}}_{n}$ in place of $\bm{p}$ in $\bm{Y}_{\bm{0},C}$, and then each quantile can be estimated using Monte-Carlo method. In Proposition B.7 in the supplementary material of \cite{baljan} it was proved that $\hat{q}_{\alpha} \stackrel{a.s.}\to q_{\alpha}$. Therefore, the following confidence band
\begin{equation*}
\Big[\max\Big((\hat{\phi}_{n,j} - \frac{\hat{q}_{\alpha}}{\sqrt{n}}), 0\Big),  \hat{\phi}_{n,j} + \frac{\hat{q}_{\alpha}}{\sqrt{n}}\Big], \, \text{for} \, j\in\mathbb{N}
\end{equation*} 
is asymptotically correct global confidence band if $\bm{p}$ is either not decreasing  or strictly decreasing, and it is asymptotically correct conservative global confidence band if $\bm{p}$ is decreasing with some constant regions.

\section{Simulation study of performance of the stacked estimators}\label{est_sim}
In this section we do simulation study to compare the performance of stacked estimators with the empirical, Grenander, rearrangement and the minimax estimators. For the p.m.f. with finite support $\{0, \dots, s \}$ and for a given sample size $n$ the minimax estimator of $\bm{p}$ with respect to $\ell_{2}$-loss is given by
\begin{equation}\label{estMM}
\hat{\bm{p}}^{mm}_{n} = \alpha^{mm}_{n} \bm{\lambda} + (1-\alpha^{mm}_{n})\hat{\bm{p}}_{n},
\end{equation}
with $\bm{\lambda} = (\frac{1}{s+1}, \dots, \frac{1}{s+1})$ and $\alpha^{mm}_{n} = \frac{\sqrt{n}}{n + \sqrt{n}}$, cf. \cite{tribula1958}. To the authors' knowledge, the minimax estimation with respect to $\ell_{2}$-loss for infinitely supported p.m.f. is an open problem.  With some abuse of notation, in this and next sections for infinitely supported distributions we refer the estimator defined in (\ref{estMM}) with $s = t_{n}$ as "minimax".

\subsection{Performance of the estimators}
We study the cases of decreasing and not decreasing true p.m.f. $\bm{p}$ separately.
\subsubsection{True p.m.f. is decreasing}
Let us consider the following uniform and decreasing p.m.f.:
\begin{eqnarray*}
\bm{M1}: \bm{p} &=& U(11), \\
\bm{M2}: \bm{p} &=& 0.15 U(3) + 0.1 U(7) + 0.75 U(11), \\
\bm{M3}: \bm{p} &=& 0.25 U(1) + 0.2 U(3) + 0.15 U(5) + 0.4 U(7),\\
\bm{M4}: \bm{p} &=& Geom(0.25), 
\end{eqnarray*}
where $Geom(\theta)$ is Geometric distribution, i.e. $p_{j} = (1-\theta)\theta^{j}$ for $j \in \mathbb{N}$ with $ 0< \theta < 1$.
\begin{figure}[t!] 
  \begin{subfigure}{2.9cm}
    \centering\includegraphics[scale=0.21]{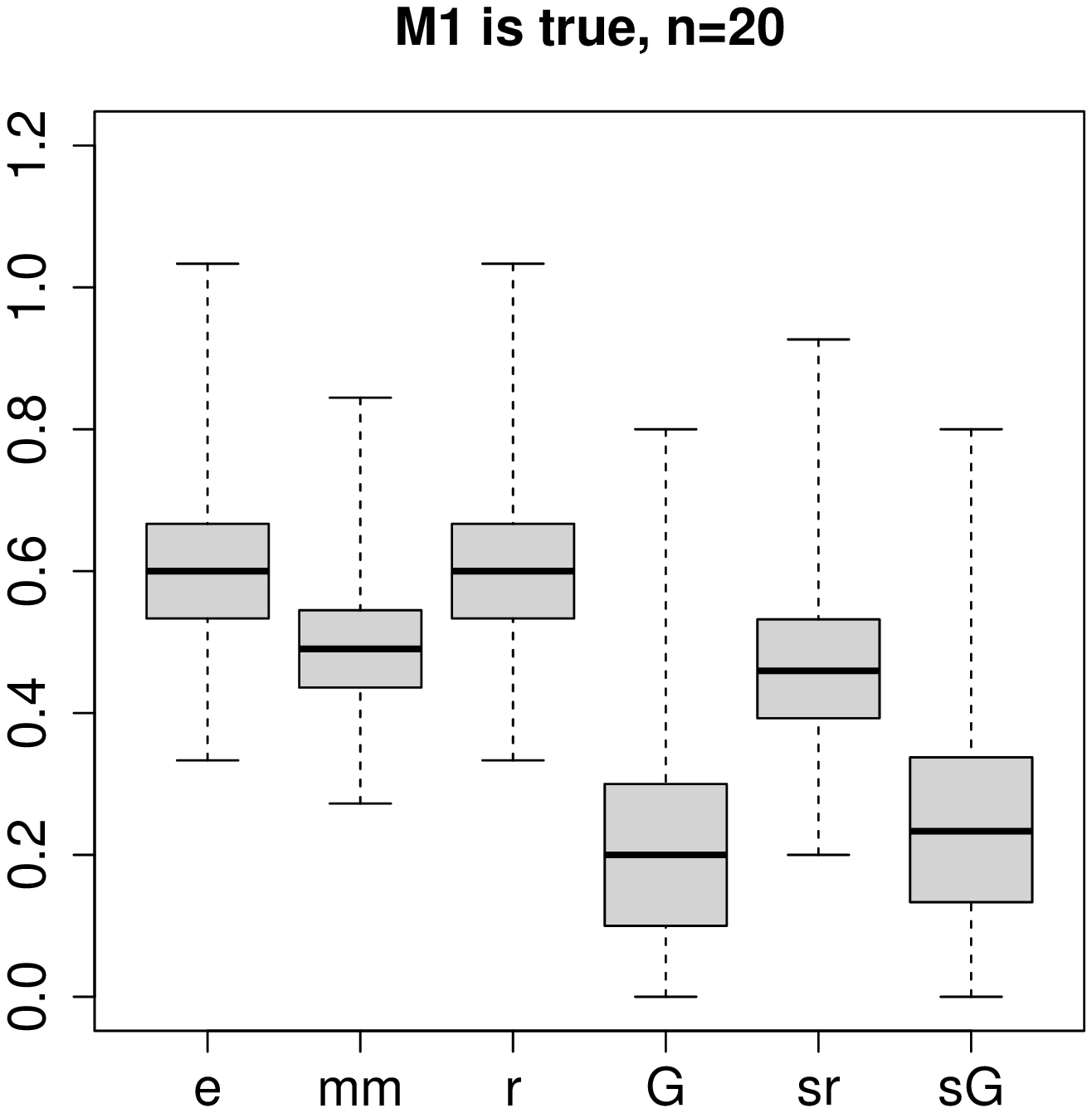}
  \end{subfigure}
  \begin{subfigure}{2.9cm}
    \centering\includegraphics[scale=0.21]{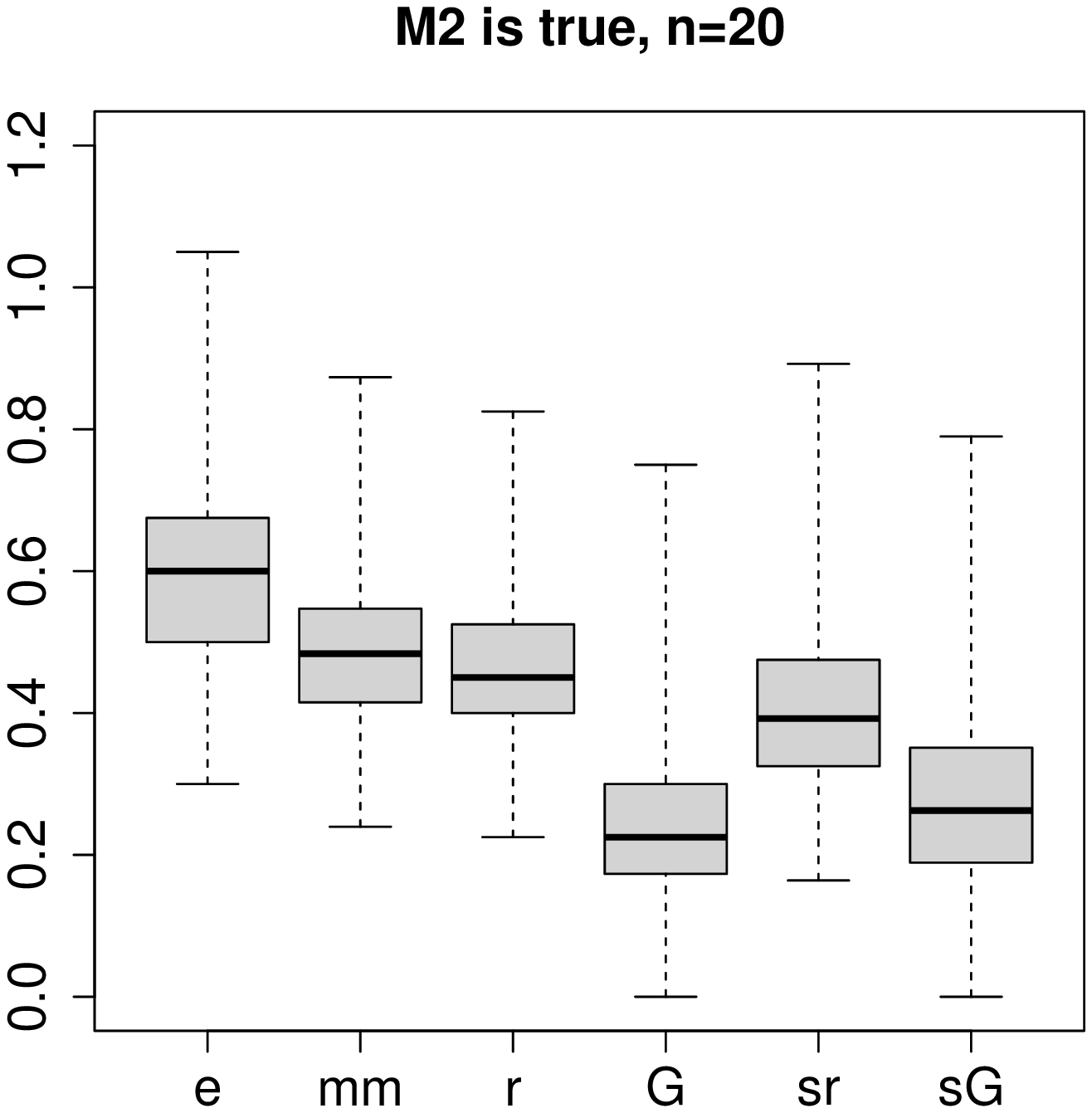}
  \end{subfigure}
  \begin{subfigure}{2.9cm}
    \centering\includegraphics[scale=0.21]{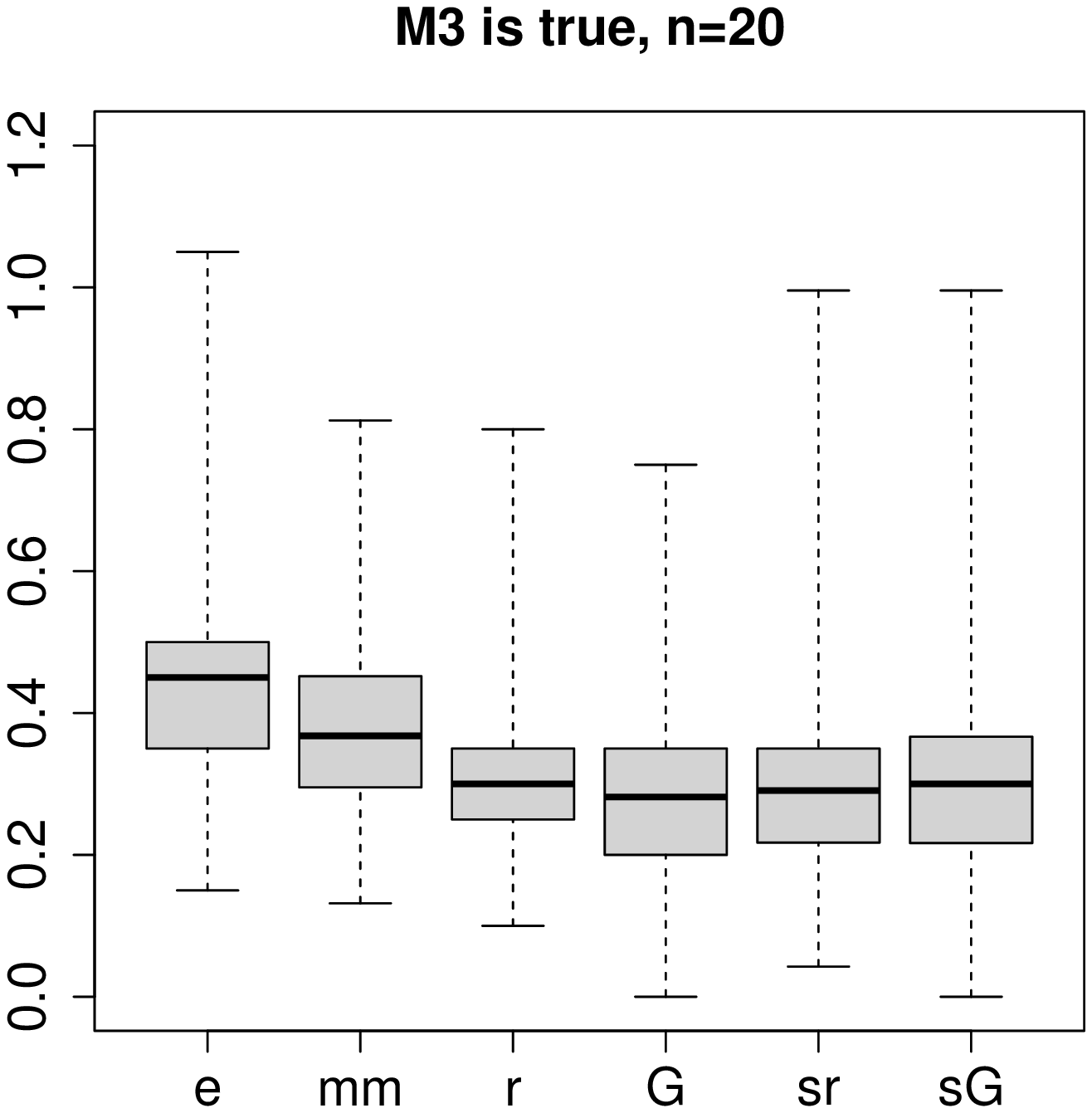}
  \end{subfigure}
   \begin{subfigure}{2.9cm}
    \centering\includegraphics[scale=0.21]{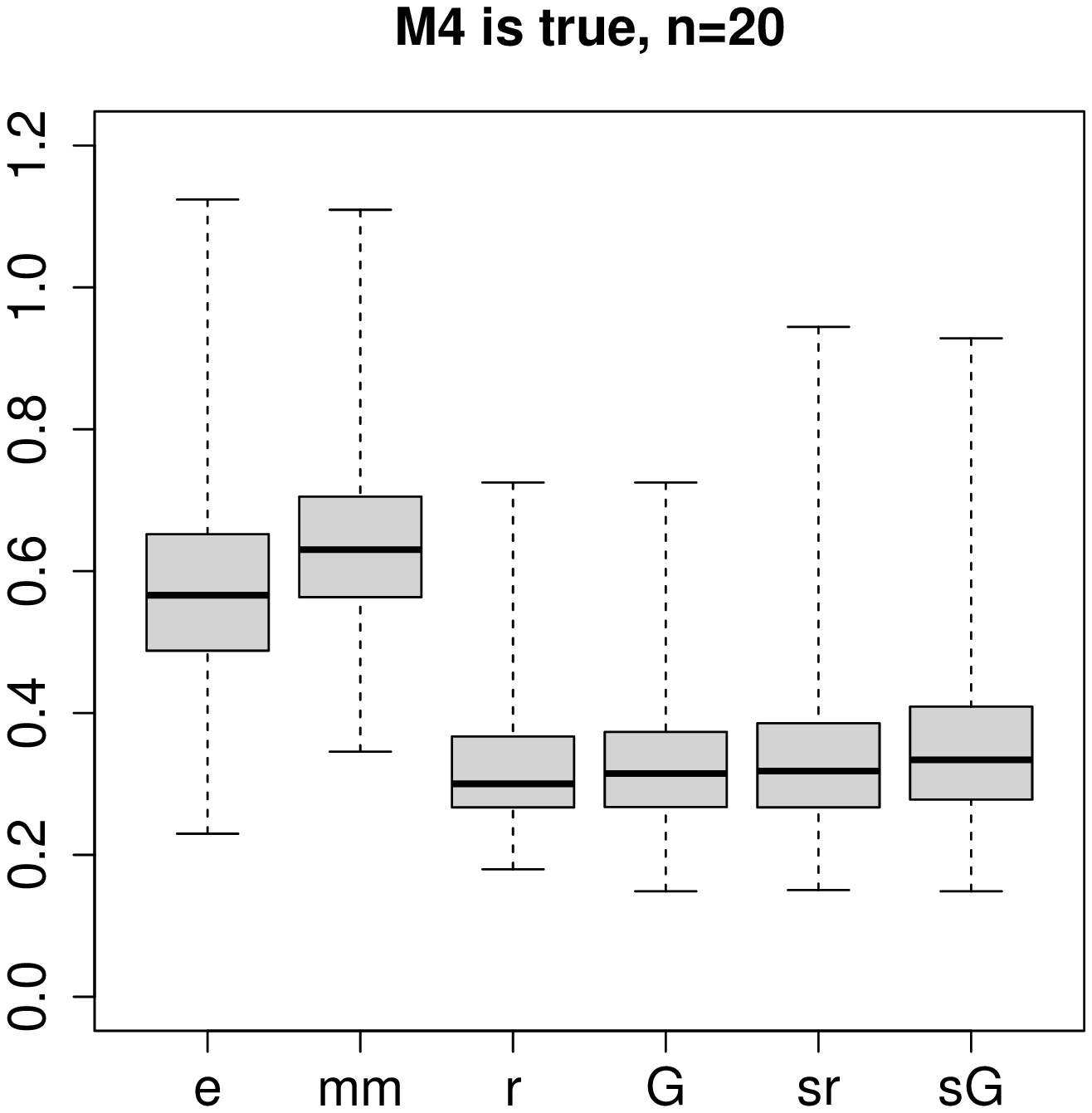}
  \end{subfigure}
  
    \begin{subfigure}{2.9cm}
    \centering\includegraphics[scale=0.21]{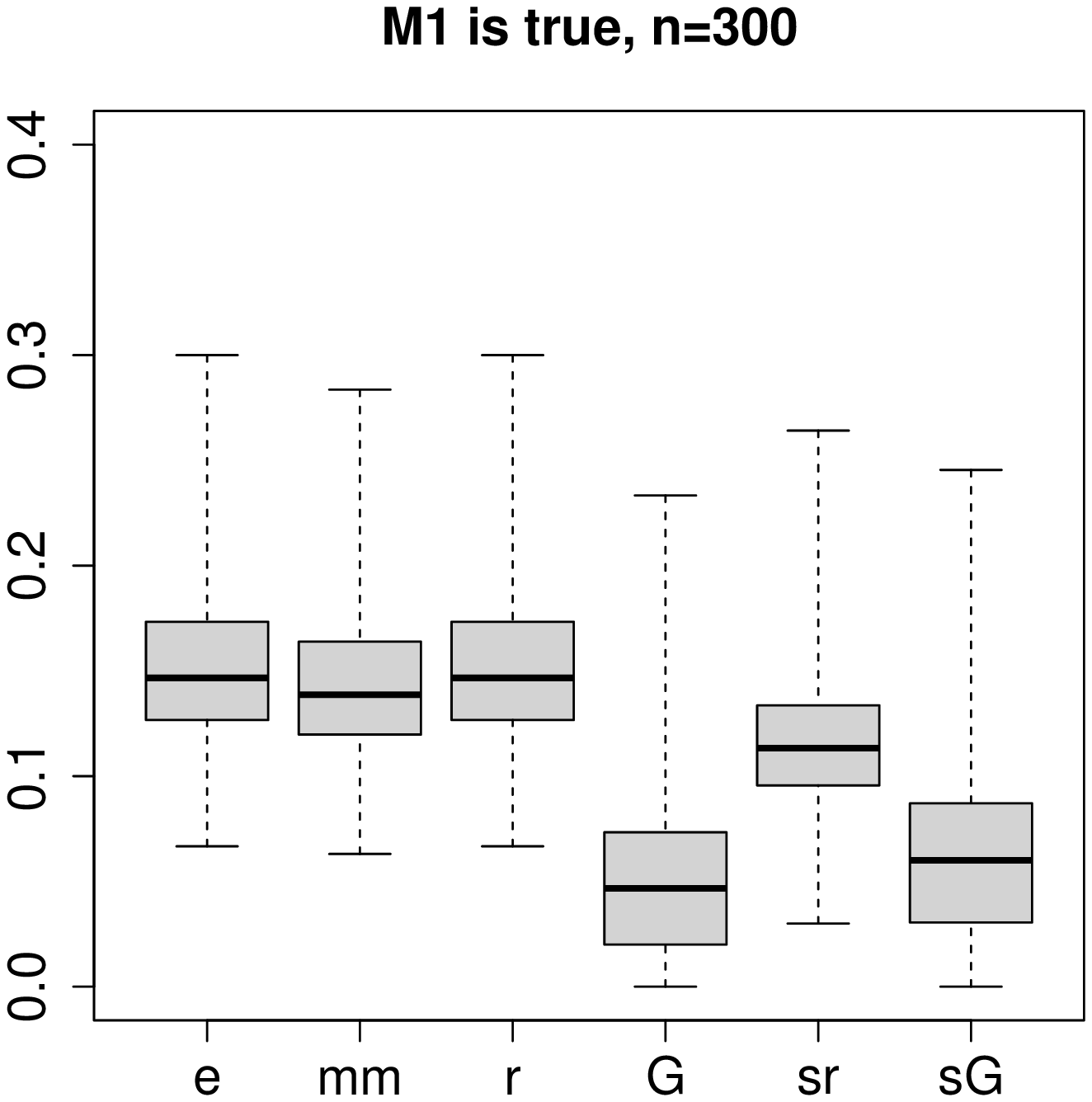}
  \end{subfigure}
  \begin{subfigure}{2.9cm}
    \centering\includegraphics[scale=0.21]{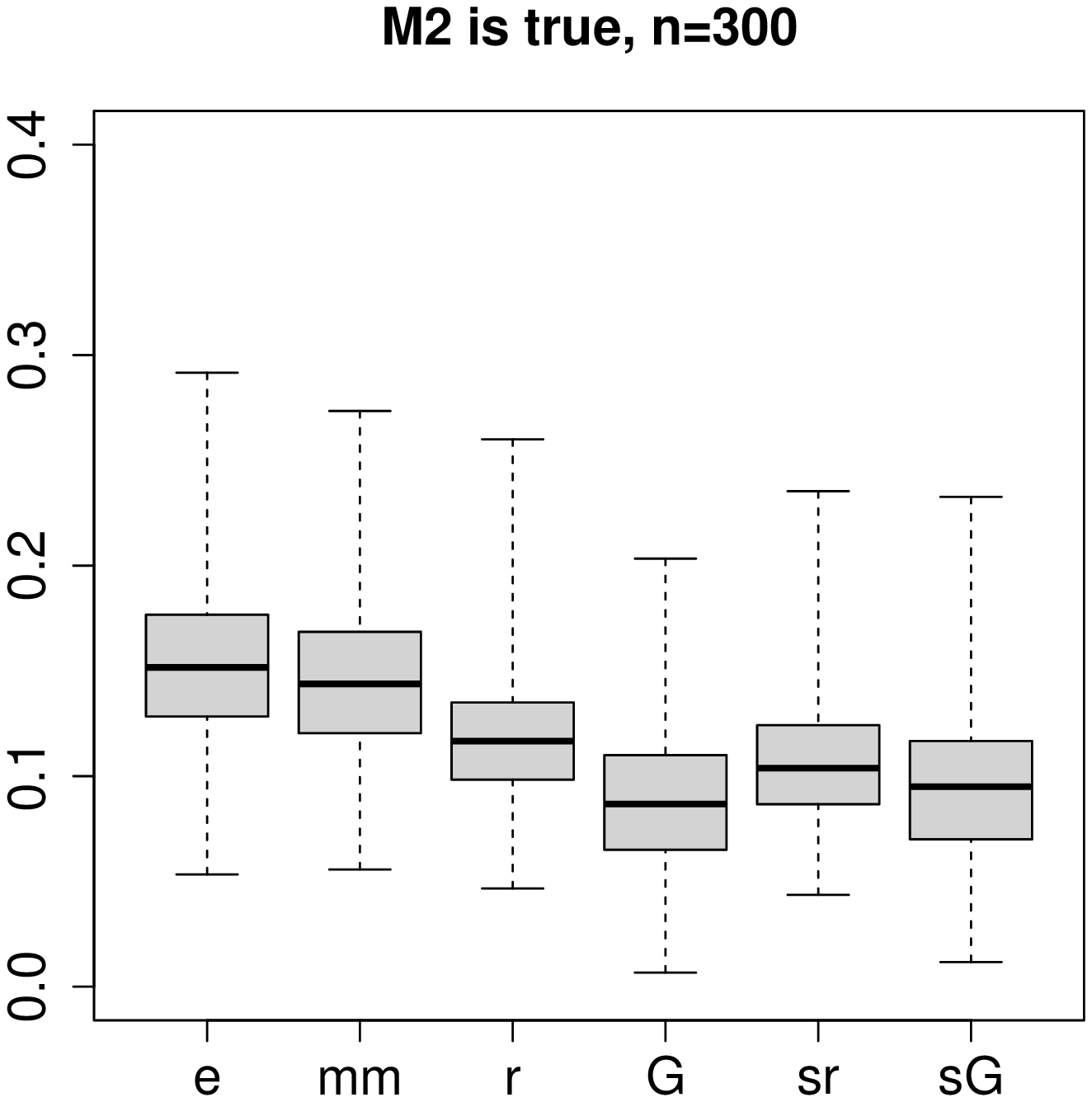}
  \end{subfigure}
  \begin{subfigure}{2.9cm}
    \centering\includegraphics[scale=0.21]{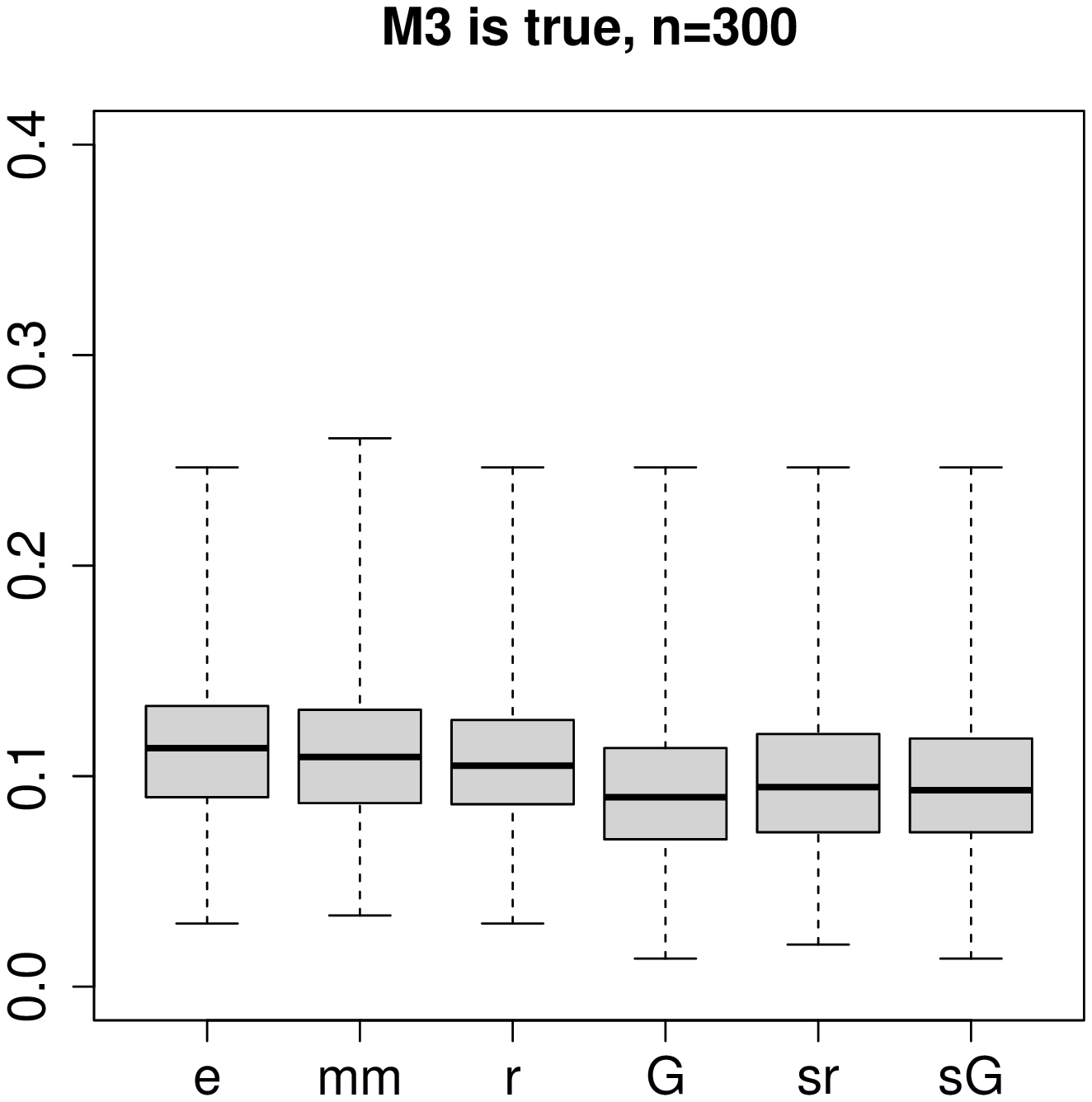}
  \end{subfigure}
   \begin{subfigure}{2.9cm}
    \centering\includegraphics[scale=0.21]{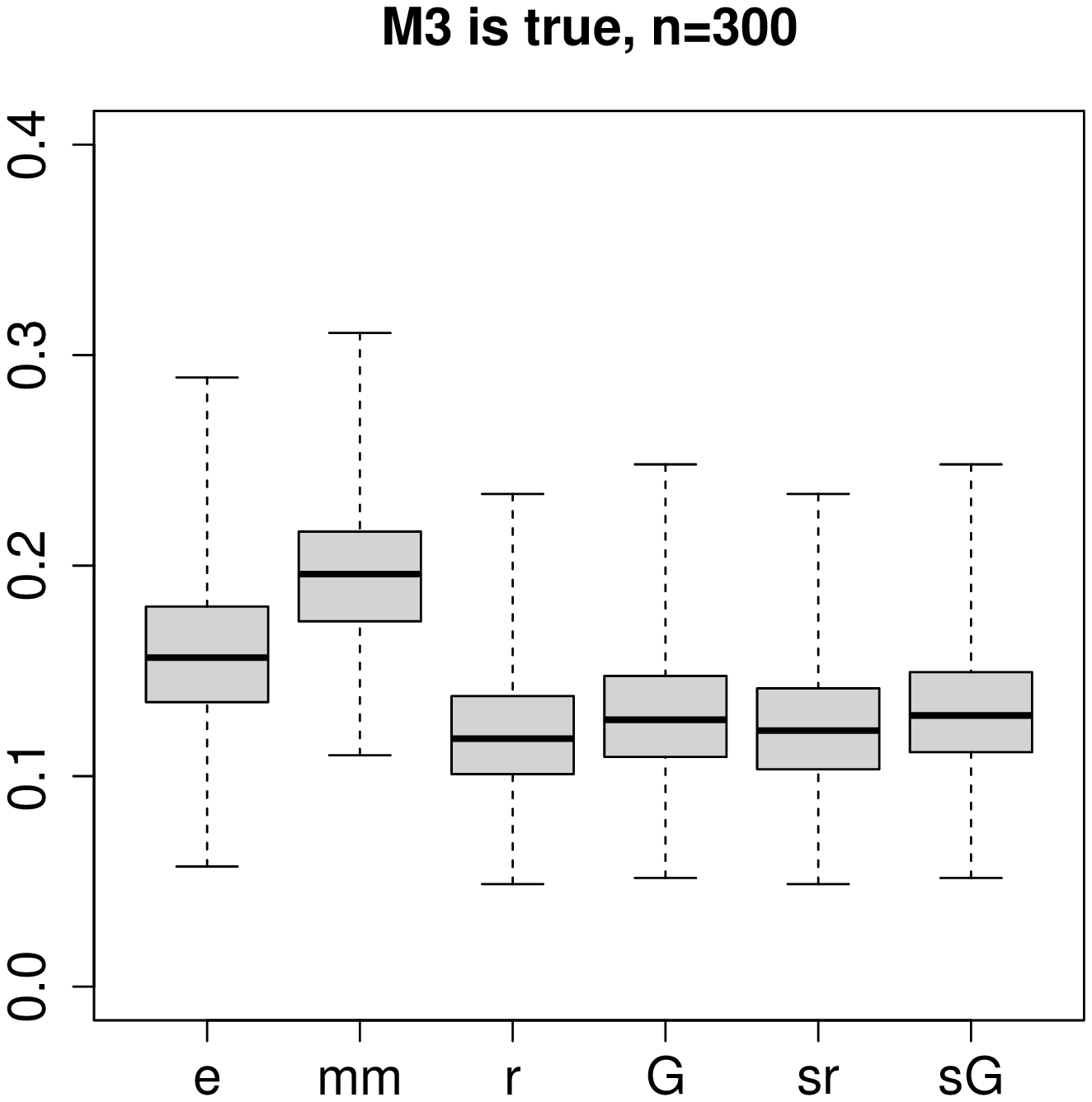}
  \end{subfigure}
   \caption{The boxplots for $\ell_1$-distances of the estimators: the empirical estimator $(e)$, minimax estimator $(mm)$, rearrangement estimator (r), Grenander estimator $(G)$, the stacked rearrangement estimator (sr) and the stacked Grenander estimator $(sG)$ for the models \textbf{M1}, \textbf{M2}, \textbf{M3} and \textbf{M4}.}\label{decr_l1}
\end{figure}

\begin{figure}[b!] 
  \begin{subfigure}{2.9cm}
    \centering\includegraphics[scale=0.21]{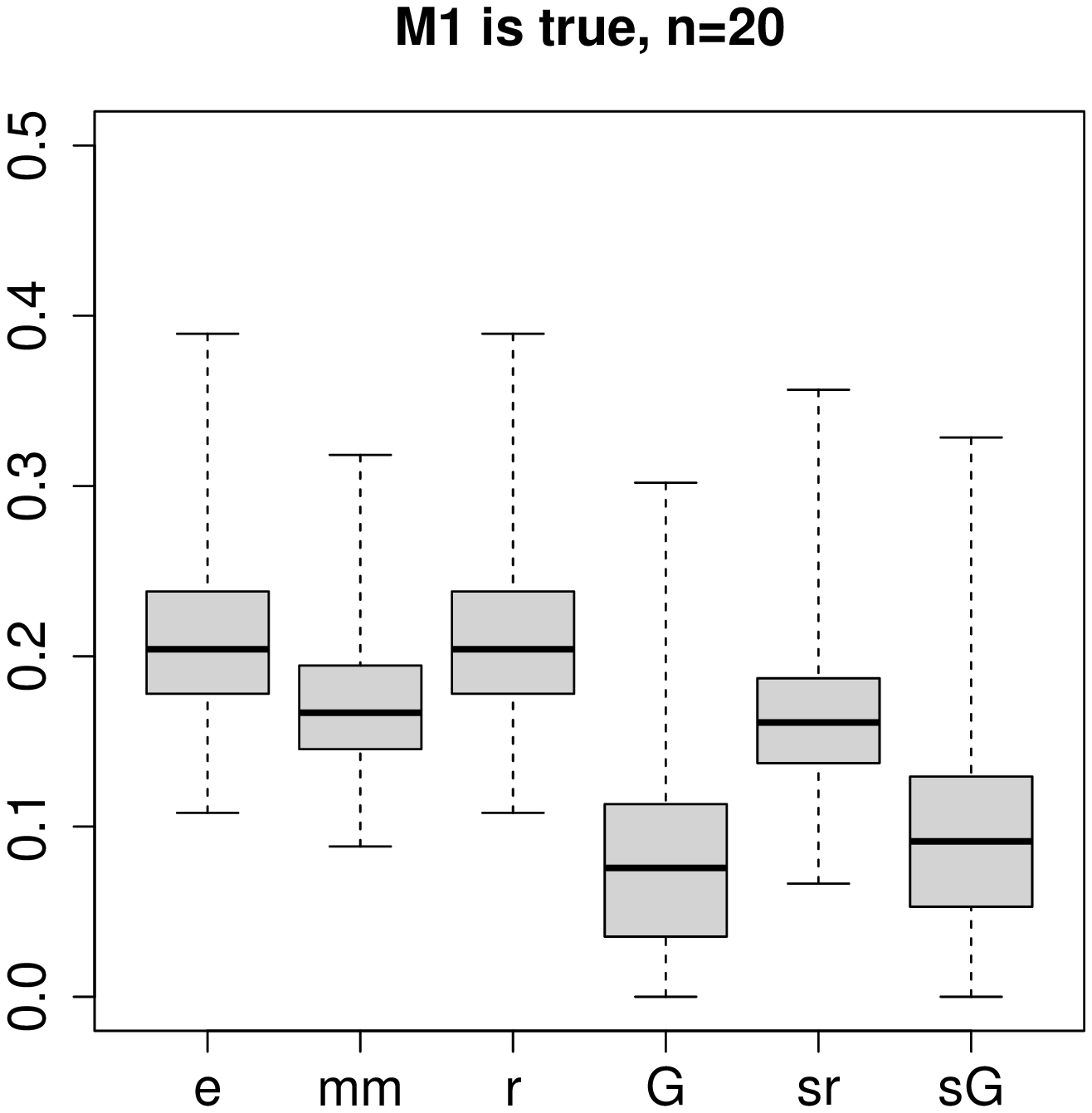}
  \end{subfigure}
  \begin{subfigure}{2.9cm}
    \centering\includegraphics[scale=0.21]{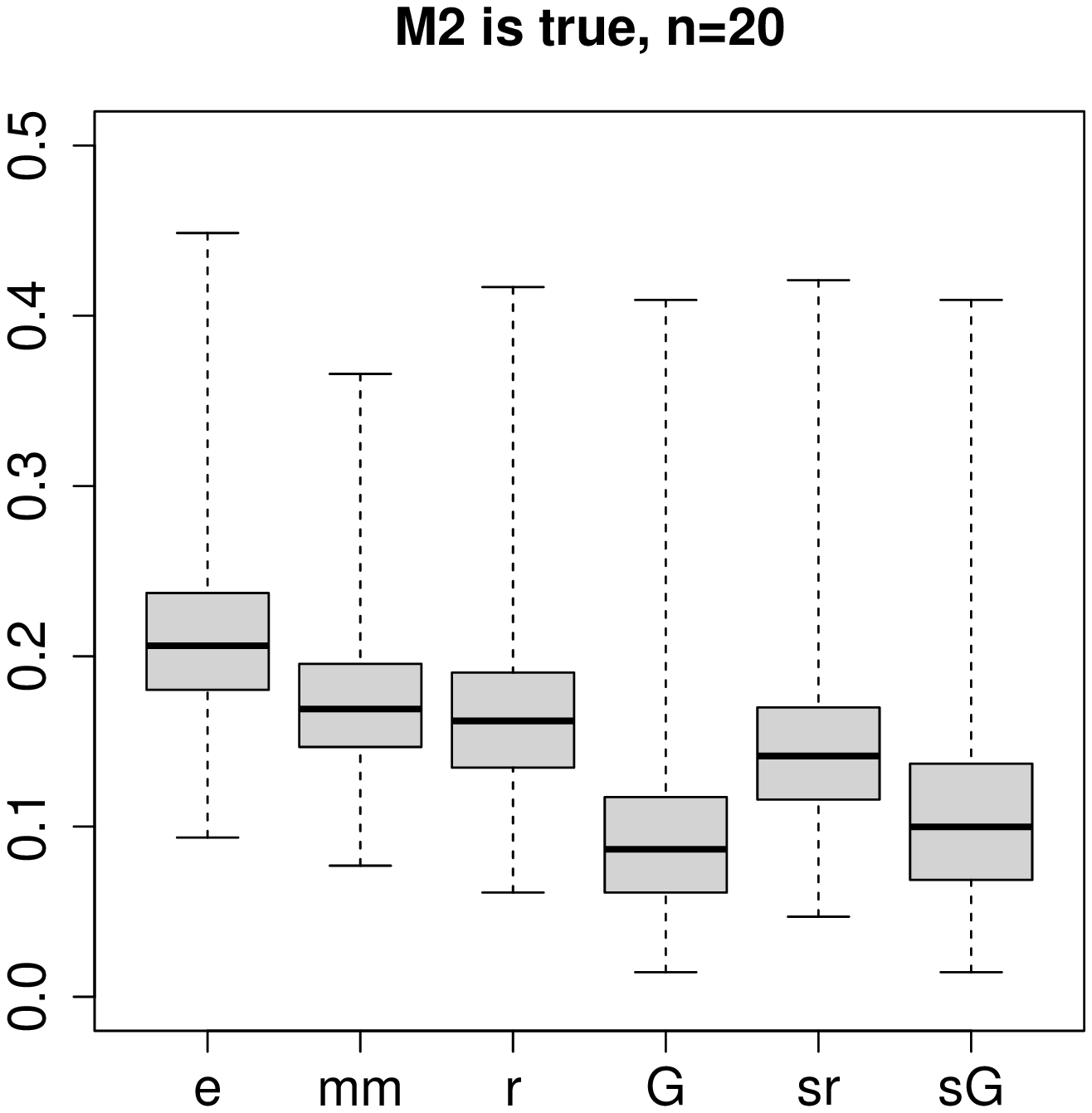}
  \end{subfigure}
  \begin{subfigure}{2.9cm}
    \centering\includegraphics[scale=0.21]{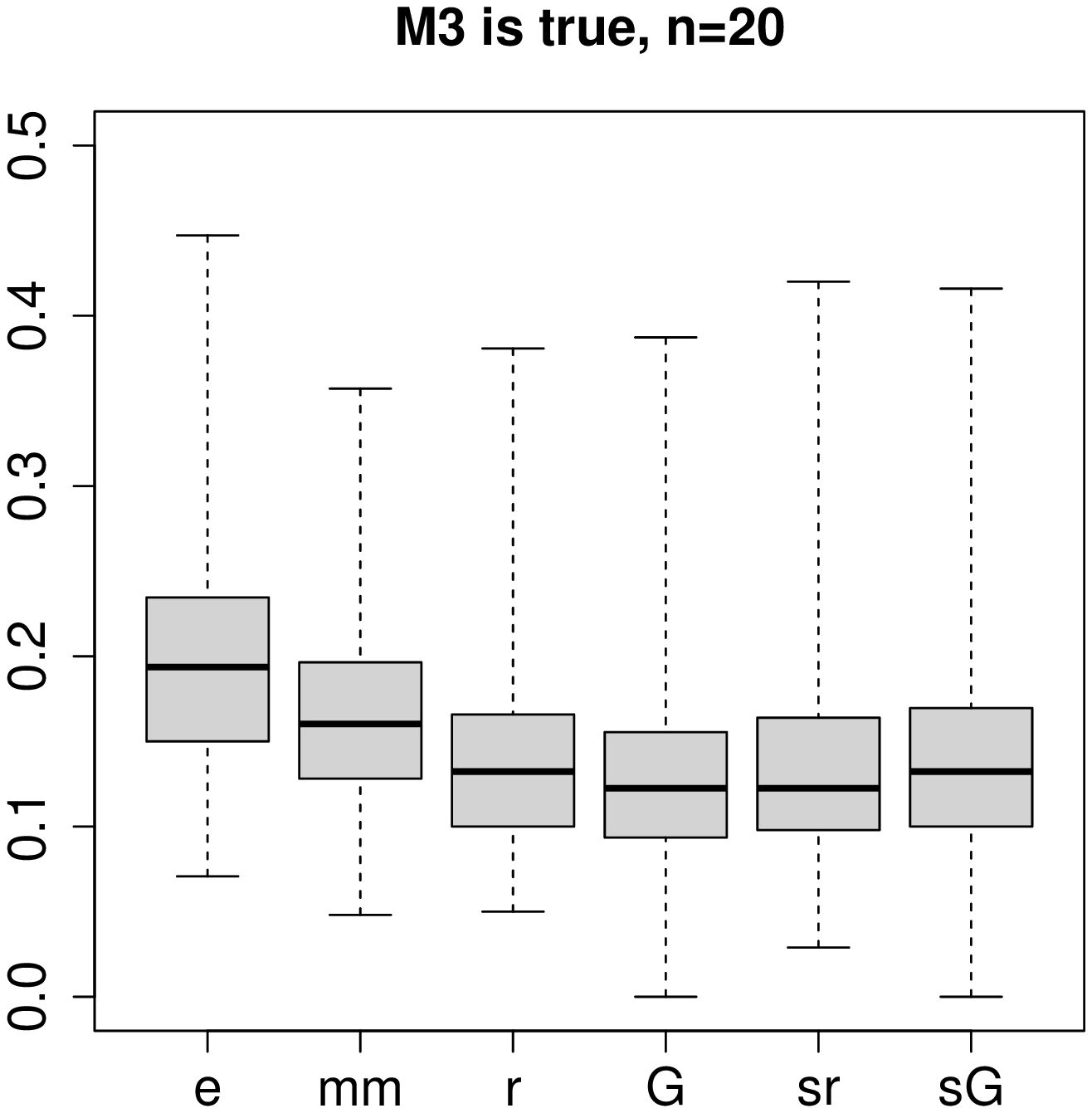}
  \end{subfigure}
   \begin{subfigure}{2.9cm}
    \centering\includegraphics[scale=0.21]{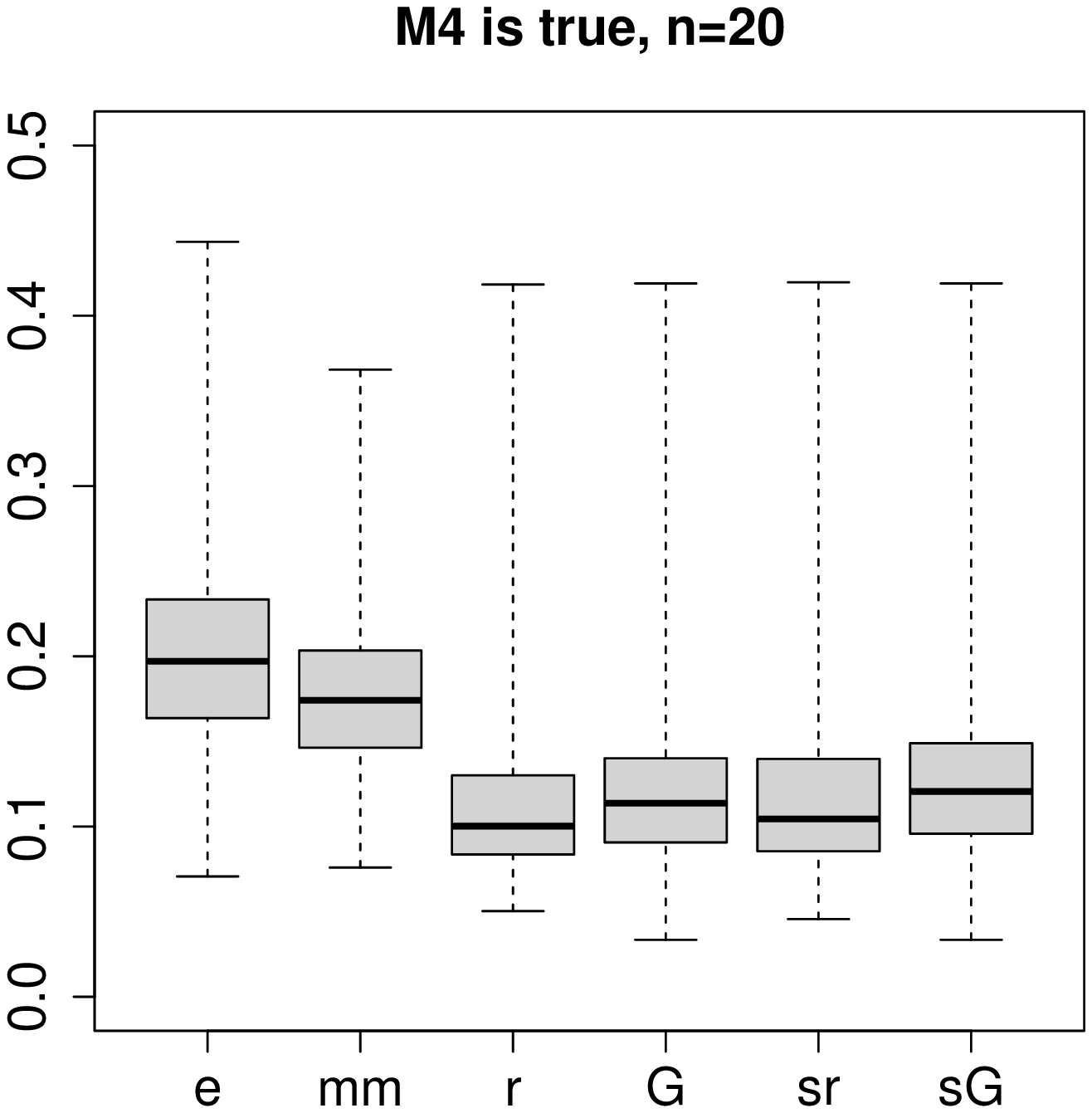}
  \end{subfigure}
  
    \begin{subfigure}{2.9cm}
    \centering\includegraphics[scale=0.21]{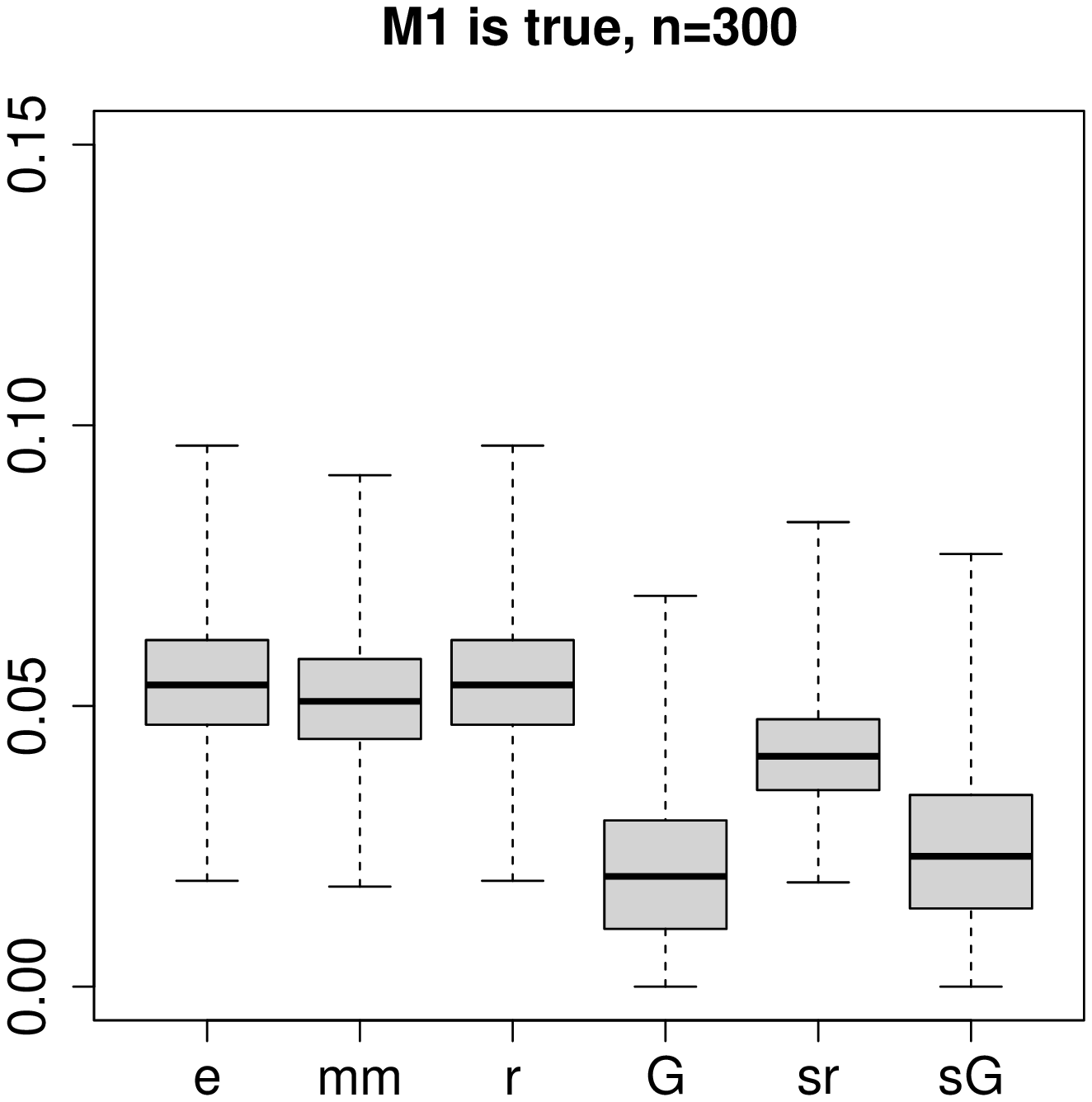}
  \end{subfigure}
  \begin{subfigure}{2.9cm}
    \centering\includegraphics[scale=0.21]{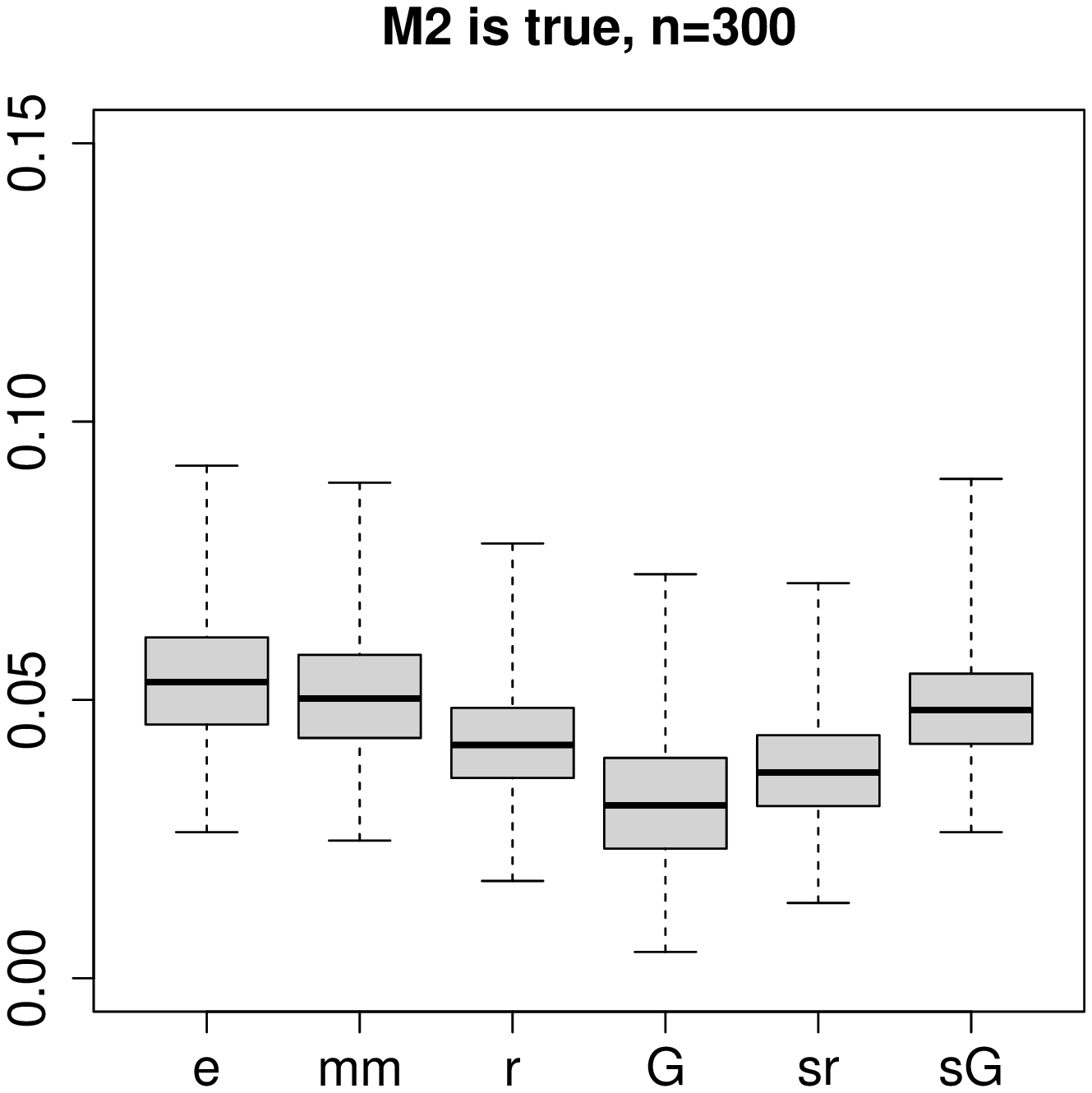}
  \end{subfigure}
  \begin{subfigure}{2.9cm}
    \centering\includegraphics[scale=0.21]{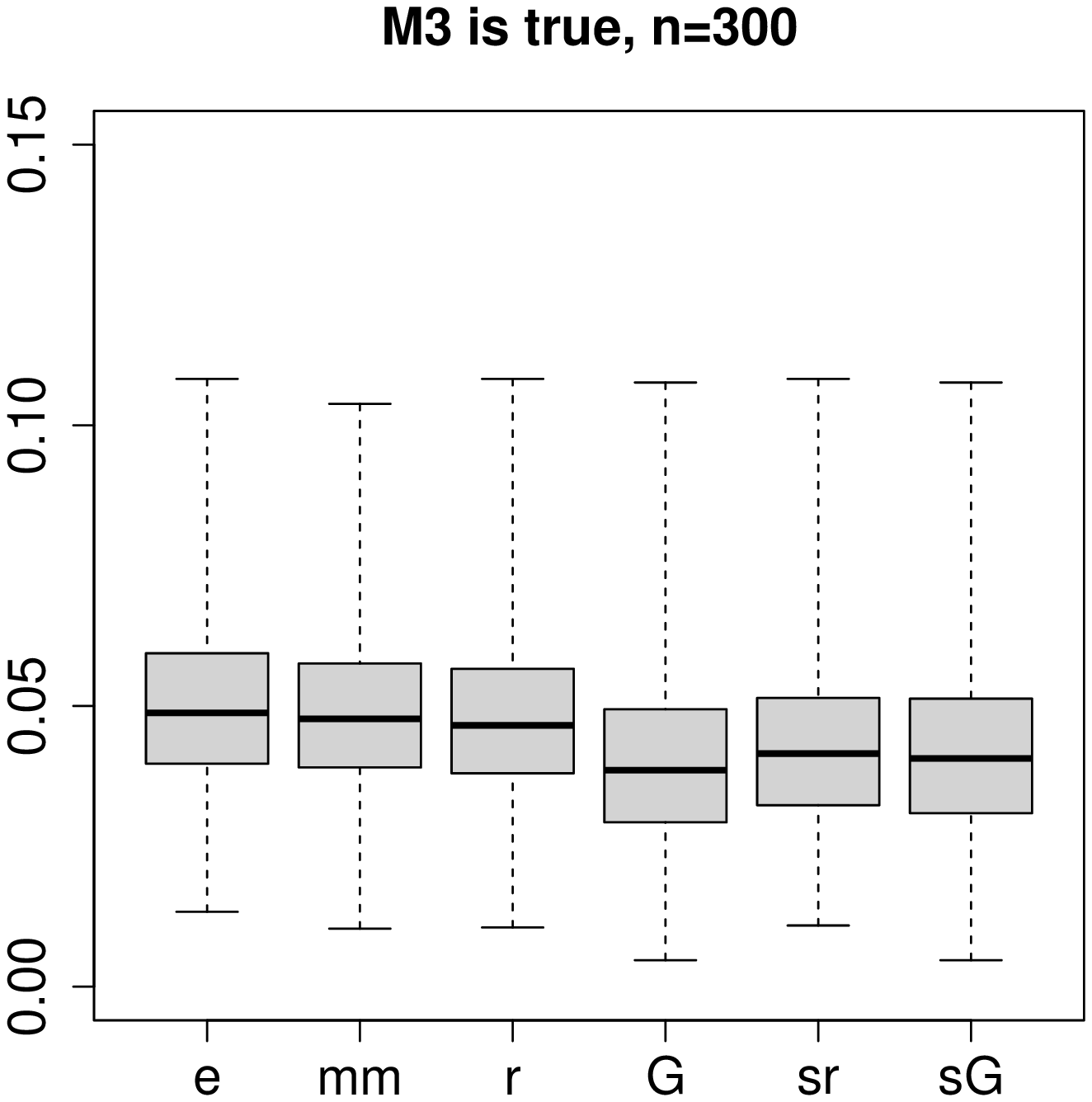}
  \end{subfigure}
   \begin{subfigure}{2.9cm}
    \centering\includegraphics[scale=0.21]{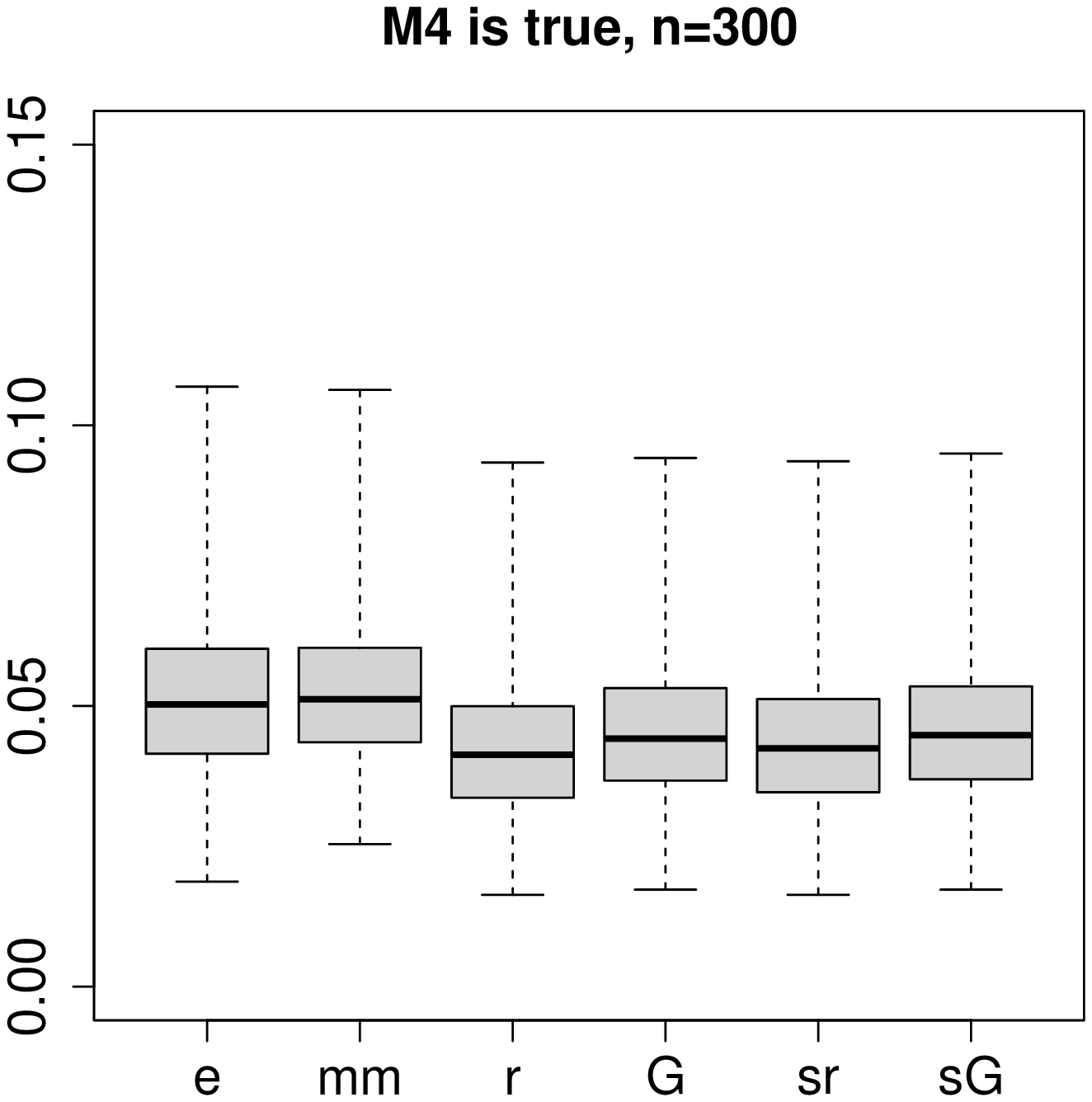}
  \end{subfigure}
   \caption{The boxplots for $\ell_2$-distances of the estimators: the empirical estimator $(e)$, minimax estimator $(mm)$, rearrangement estimator (r), Grenander estimator $(G)$, the stacked rearrangement estimator (sr) and the stacked Grenander estimator $(sG)$ for the models \textbf{M1}, \textbf{M2}, \textbf{M3} and \textbf{M4}.}\label{decr_l2}
\end{figure} 

The models $\bm{M2}$, $\bm{M3}$ and $\bm{M4}$ were used in  \cite{jankowski2009estimation} to assess the performance of Grenander estimator and compare its performance with empirical and rearrangement estimators. First, we compare the performance of the estimators in $\ell_{1}$ (Figure \ref{decr_l1}) and $\ell_{2}$ (Figure \ref{decr_l2}) distances  for small $n=20$ and moderate $n=300$ sample sizes with 1000 Monte Carlo simulations.

From the boxplots at Figure \ref{decr_l1} and Figure \ref{decr_l2} we can conclude that for both small and moderate sized data sets stacked Grenander estimator outperforms in $\ell_{1}$ and $\ell_{2}$ norms both the empirical estimator and minimax estimator ("minimax" for the case of Geometric distribution). Further, stacked Grenander estimator outperforms stacked rearrangement estimator when the underlying distribution has constant regions and it performs almost the same in the case of strictly decreasing p.m.f. The superiority of Grenander estimator over the rearrangement estimator was proved in \cite{jankowski2009estimation}.

\begin{figure}[t!] 
  \begin{subfigure}{2.9cm}
    \centering\includegraphics[scale=0.21]{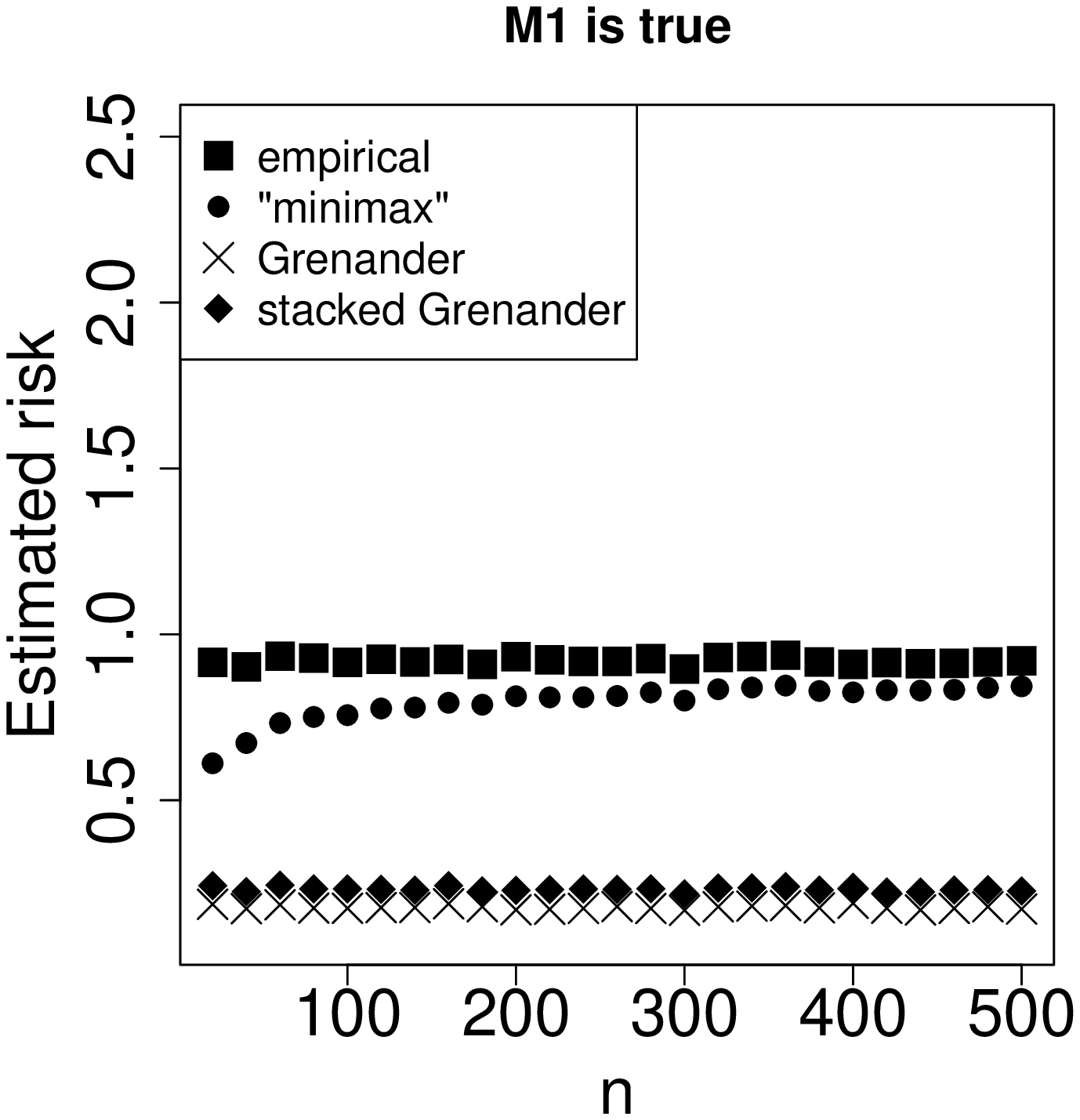}
  \end{subfigure}
  \begin{subfigure}{2.9cm}
    \centering\includegraphics[scale=0.21]{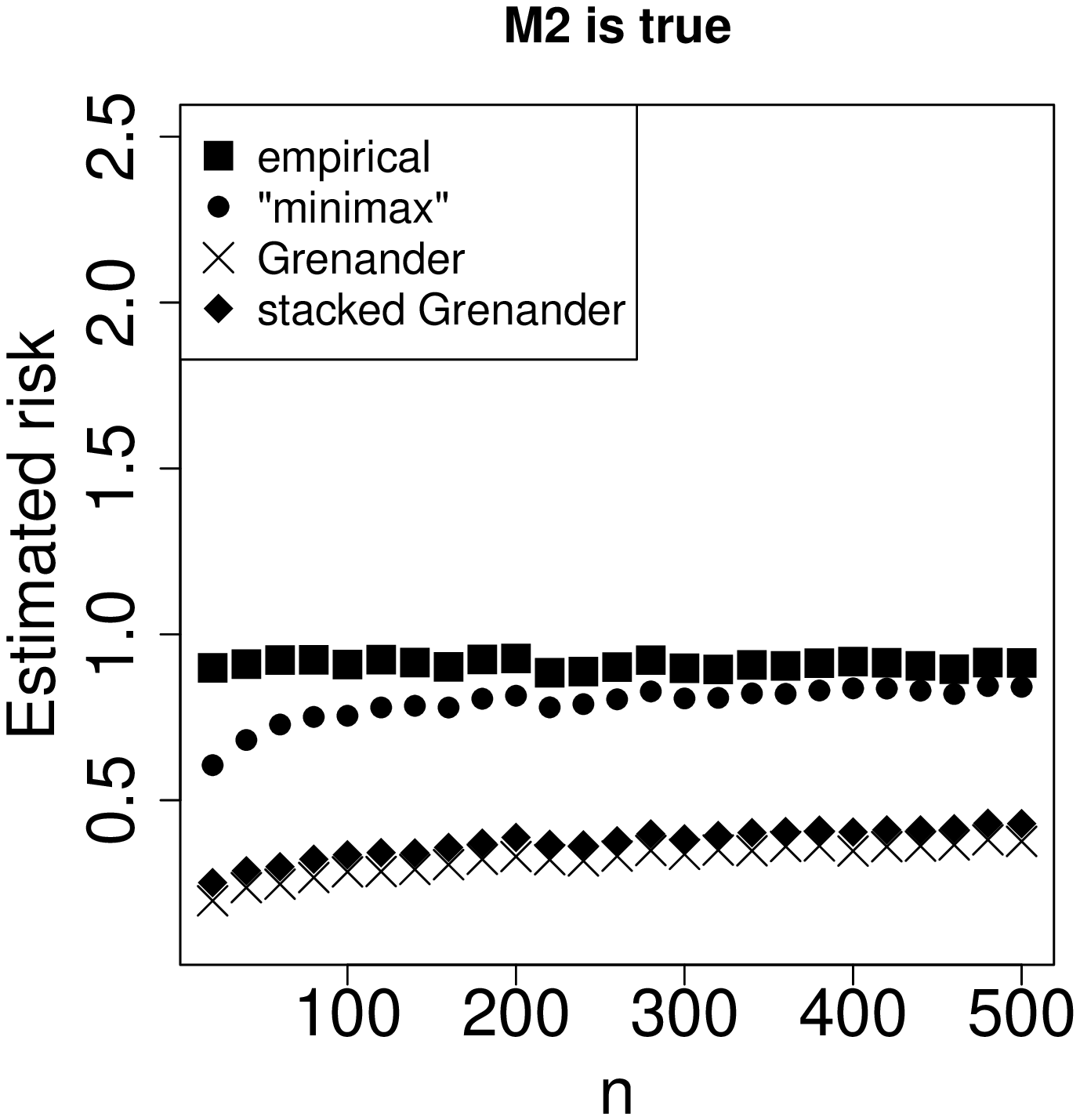}
  \end{subfigure}
  \begin{subfigure}{2.9cm}
    \centering\includegraphics[scale=0.21]{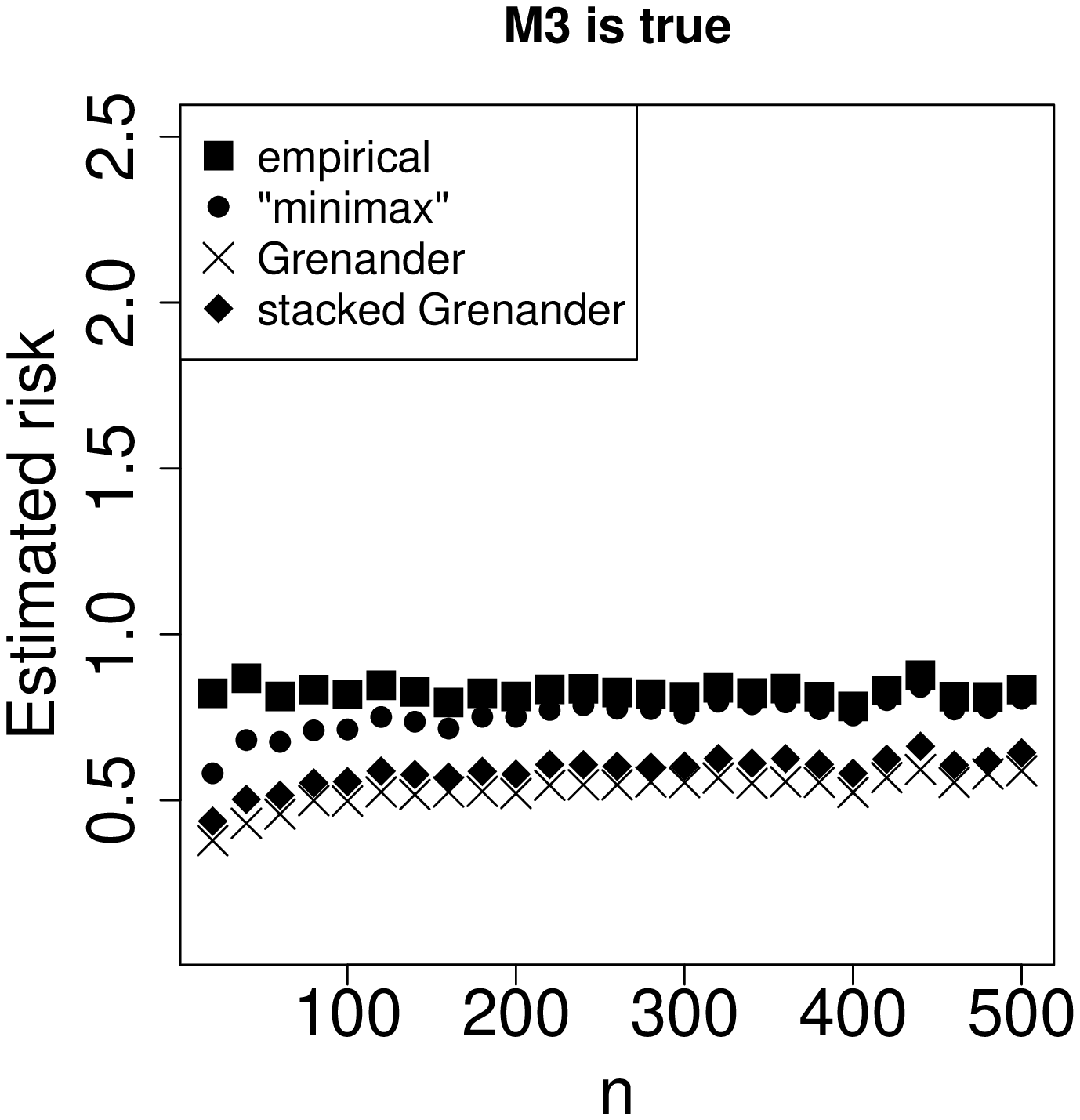}
  \end{subfigure}
   \begin{subfigure}{2.9cm}
    \centering\includegraphics[scale=0.21]{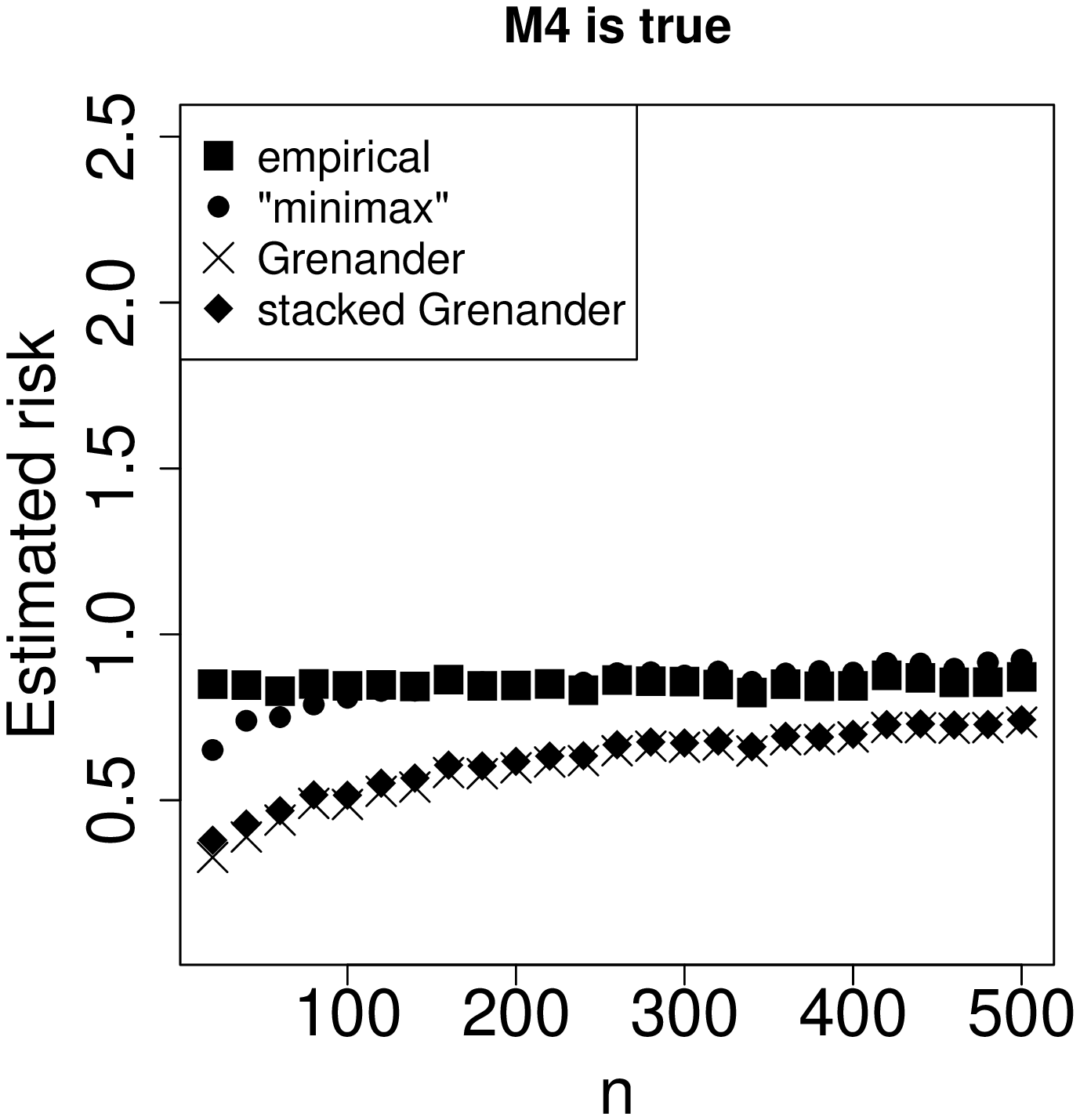}
  \end{subfigure}
   \caption{The estimates of the scaled risk for the models \textbf{M1}, \textbf{M2}, \textbf{M3} and \textbf{M4}.}\label{decrRisk}
\end{figure}

Next, in order to summarise the results and demonstrate the superiority of stacked Grenander estimator we plot the estimates of scaled risk $n\mathbb{E}[||\hat{\bm{\xi}}_{n} - \bm{p}||_{2}^{2}]$ (with $\hat{\bm{\xi}}_{n}$ one of the following estimators: empirical, minimax Grenander or stacked Grenander estimator) versus the sample size $n$, based on $1000$ Monte Carlo simulations, cf. Figure \ref{decrRisk}. We can conclude that in the case of a decreasing underlying distribution stacked Grenander estimator performs almost as good as Grenander estimator and it performs significantly better than the empirical and the minimax estimators.

\subsubsection{True p.m.f. is not decreasing}
Now let us consider the case when the underlying distributions are not decreasing:
\begin{eqnarray*}
\bm{M5}: \bm{p} &=& T^{i}(11), \\
\bm{M6}: \bm{p} &=& NBin(7, 0.4), \\
\bm{M7}: \bm{p} &=& \frac{3}{8}Pois(2) + \frac{5}{8}Pois(15),
\end{eqnarray*}
where $T^{i}(s)$ stands for strictly increasing triangular function; $NBin(r,\theta)$ is the negative binomial distribution with $r$ the number of failures until the experiment is stopped and $\theta$ the success probability; $Pois(\lambda)$ is Poisson distribution with rate $\lambda$. Therefore, we consider very non-monotonic distributions. Indeed, model $\bm{M5}$ is a strictly increasing p.m.f., $\bm{M6}$ is a unimodal distribution, and $\bm{M7}$ is bimodal.

\begin{figure}[!htbp] 
  \begin{subfigure}{3.9cm}
    \centering\includegraphics[scale=0.211]{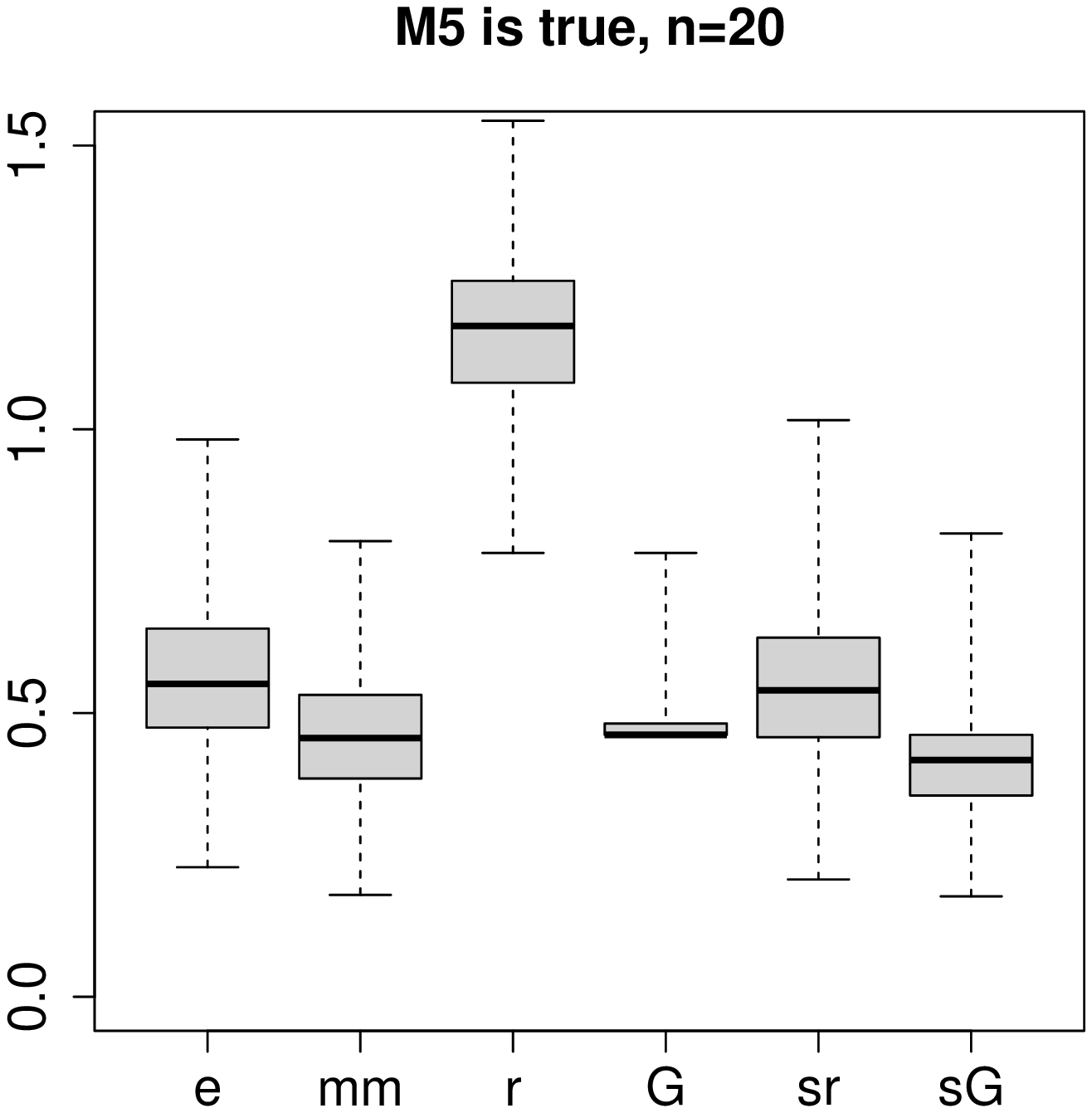}
  \end{subfigure}
  \begin{subfigure}{3.9cm}
    \centering\includegraphics[scale=0.211]{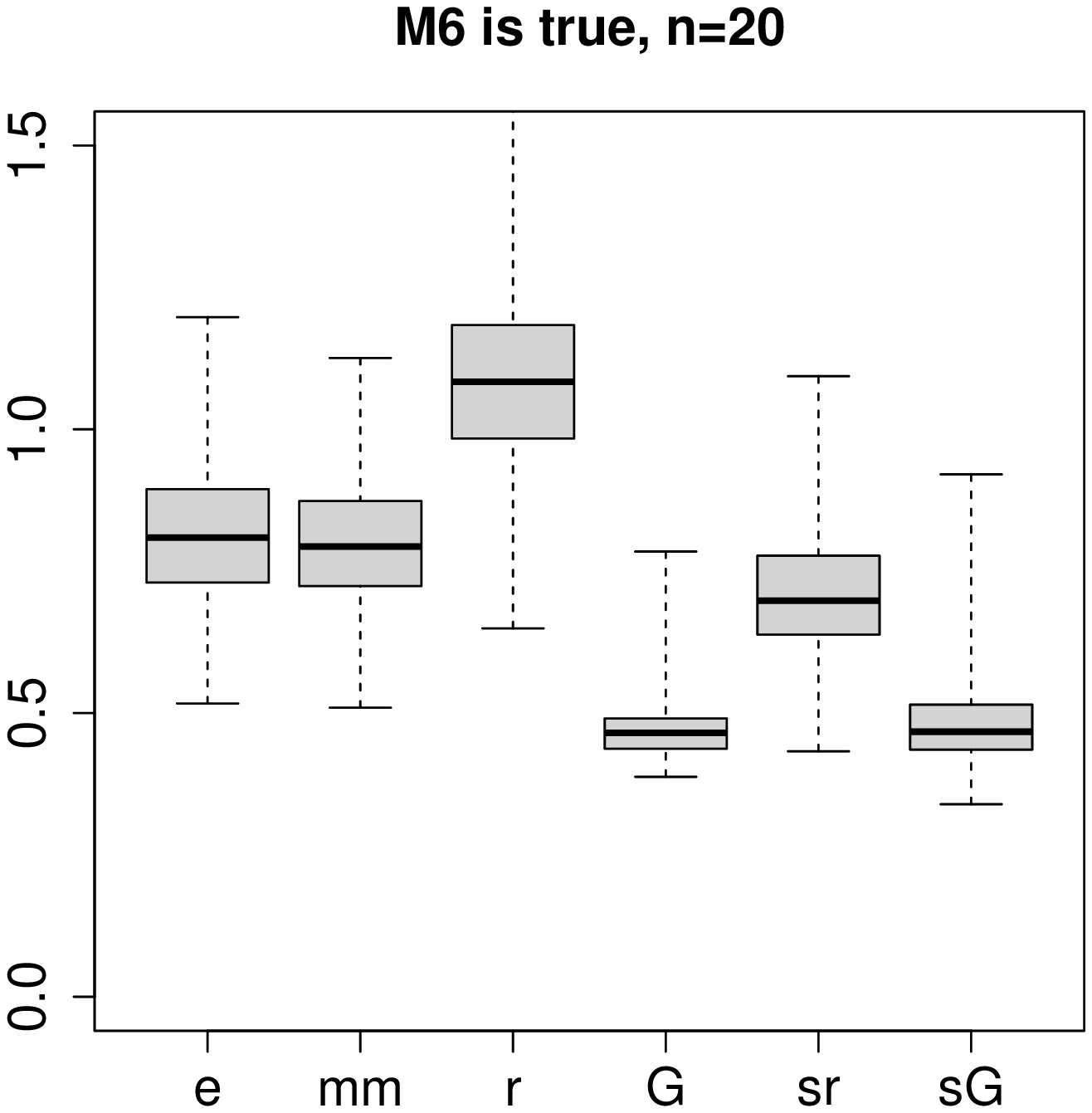}
  \end{subfigure}
  \begin{subfigure}{3.9cm}
    \centering\includegraphics[scale=0.211]{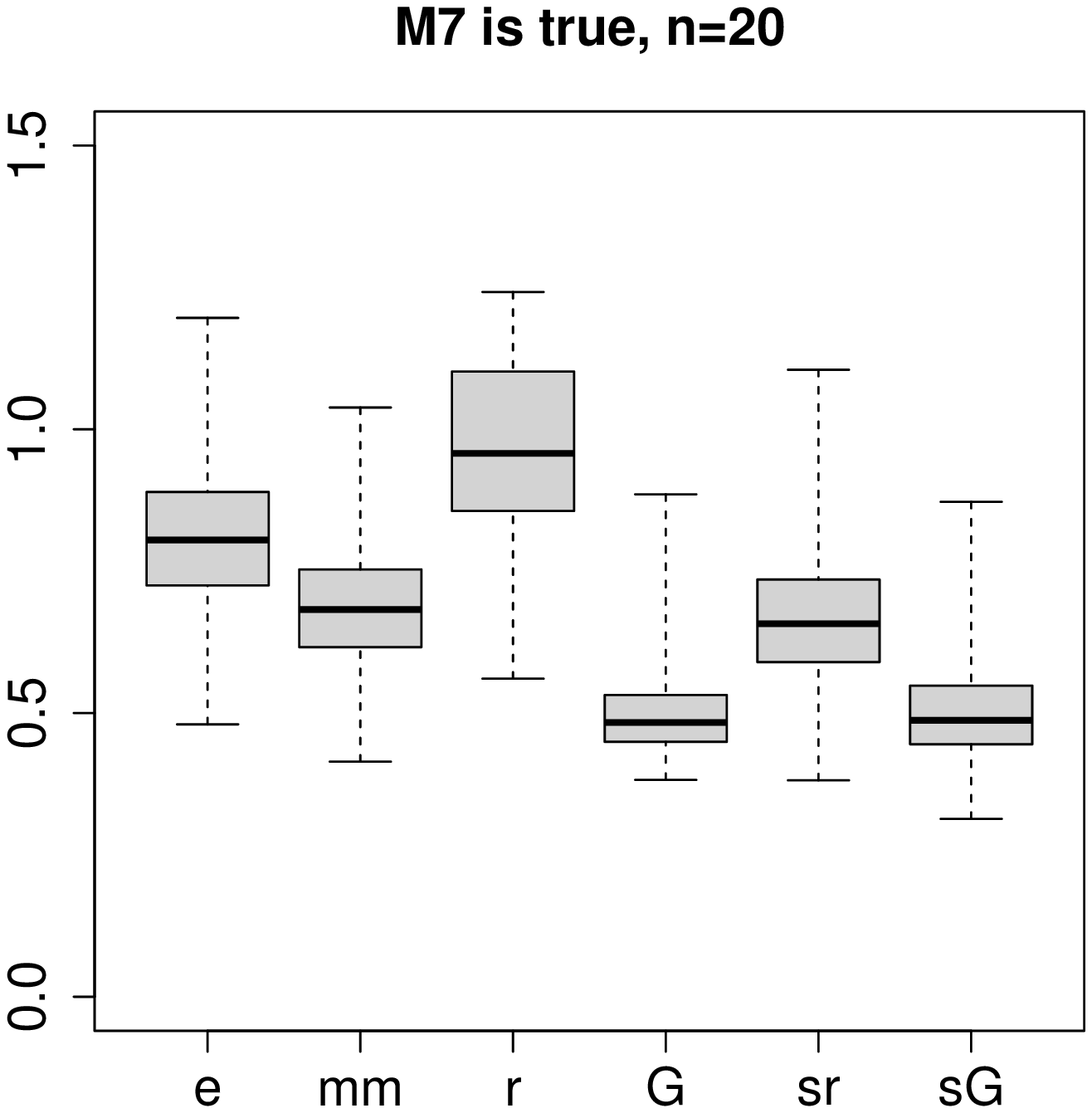}
  \end{subfigure}
  
    \begin{subfigure}{3.9cm}
    \centering\includegraphics[scale=0.211]{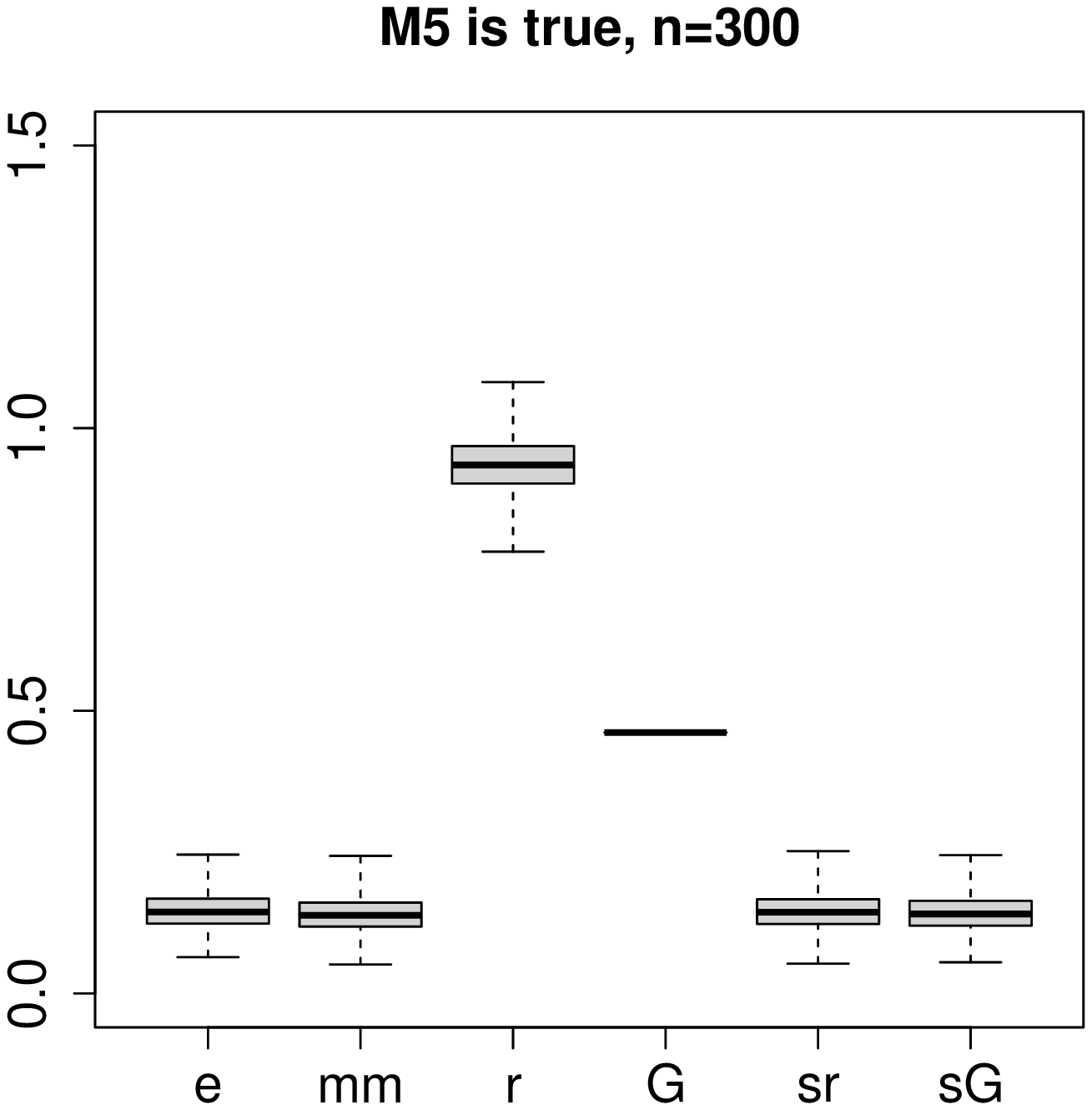}
  \end{subfigure}
  \begin{subfigure}{3.9cm}
    \centering\includegraphics[scale=0.211]{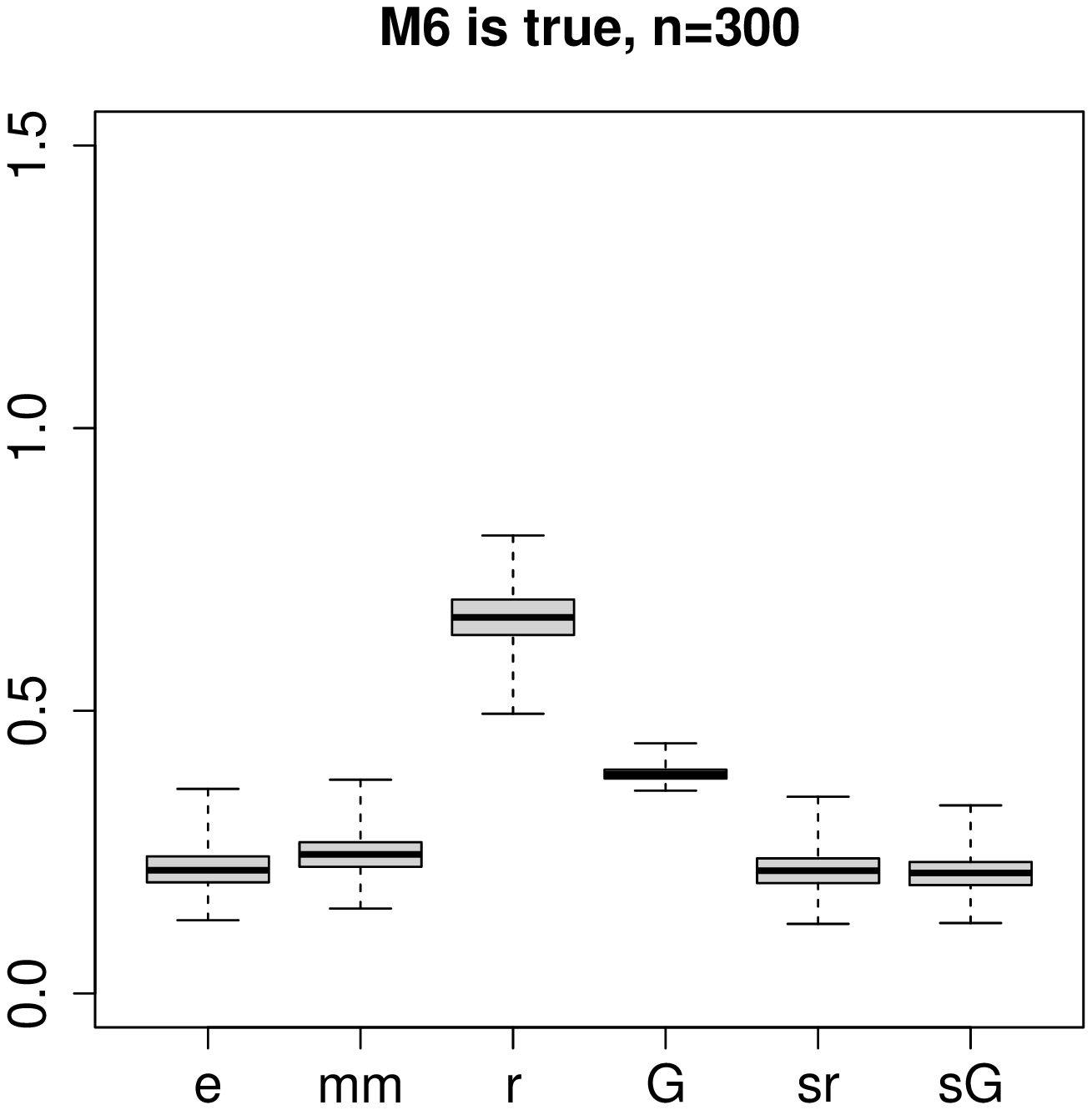}
  \end{subfigure}
  \begin{subfigure}{3.9cm}
    \centering\includegraphics[scale=0.211]{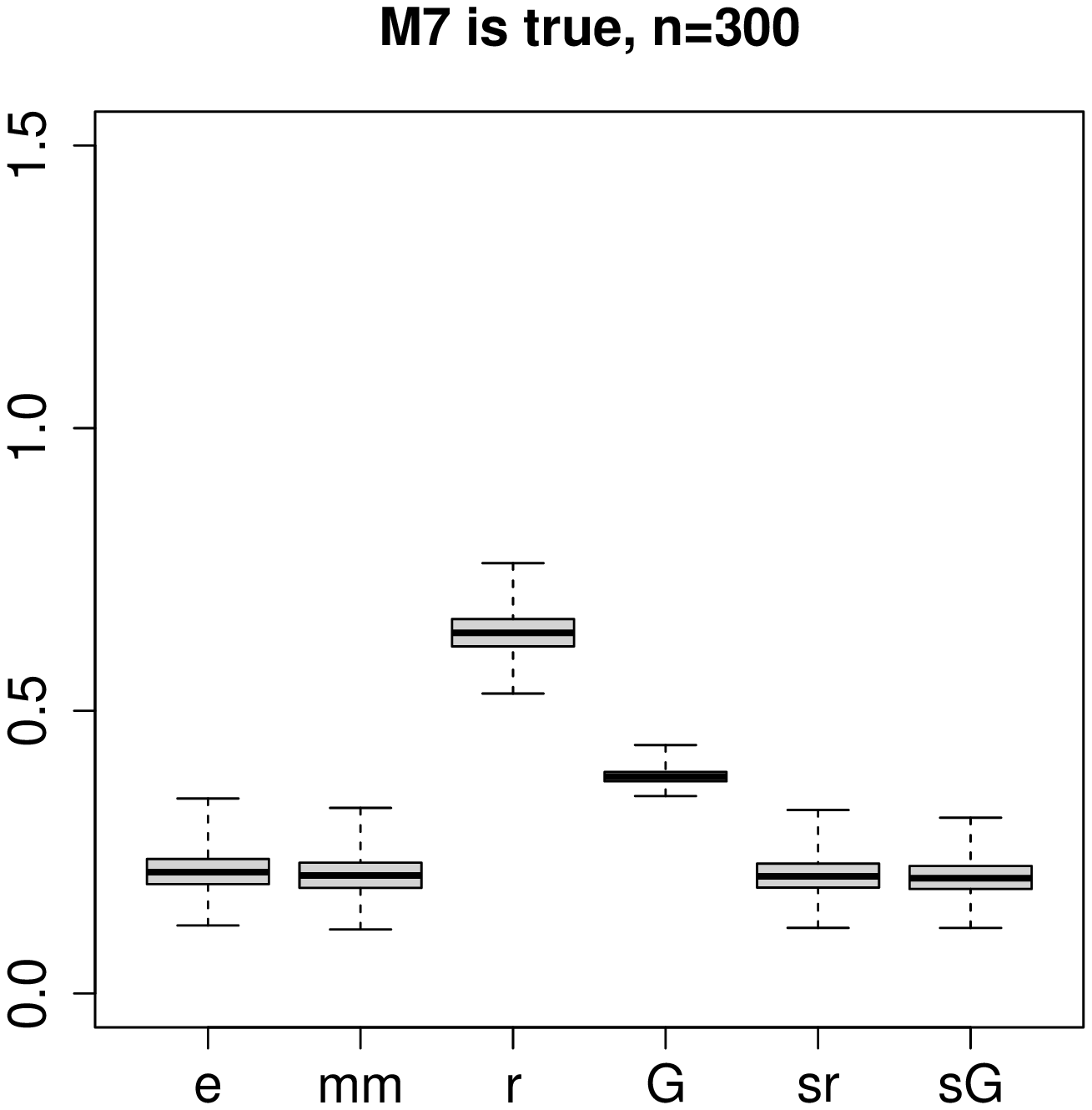}
  \end{subfigure}
   \caption{The boxplots for $\ell_1$-distances of the estimators: the empirical estimator $(e)$, minimax estimator $(mm)$, rearrangement estimator (r), Grenander estimator $(G)$, the stacked rearrangement estimator (sr) and the stacked Grenander estimator $(sG)$ for the models \textbf{M5}, \textbf{M6} and \textbf{M7}.}\label{nondecr_l1}
\end{figure}

\begin{figure}[!htbp] 
  \begin{subfigure}{3.9cm}
    \centering\includegraphics[scale=0.211]{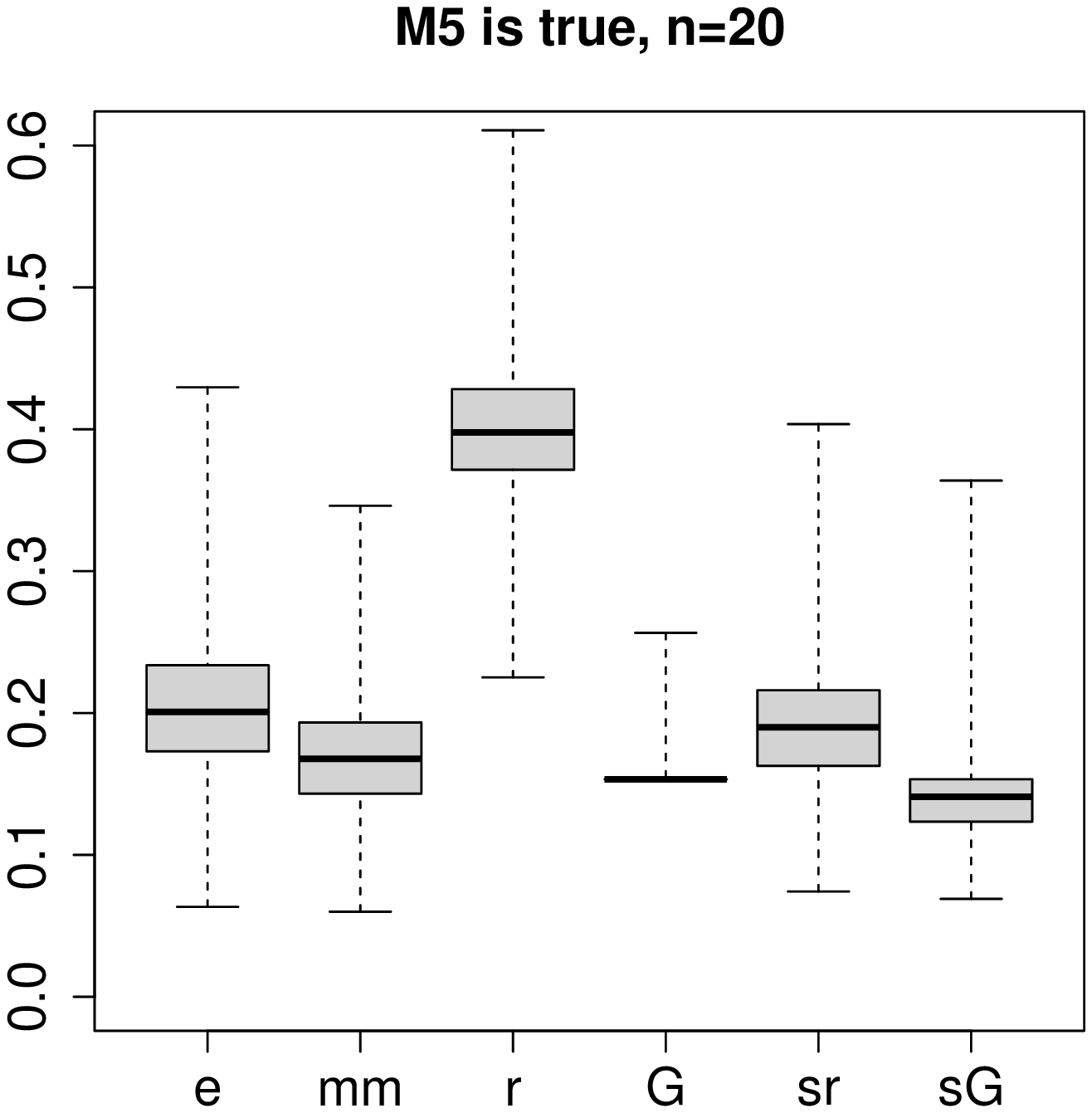}
  \end{subfigure}
  \begin{subfigure}{3.9cm}
    \centering\includegraphics[scale=0.211]{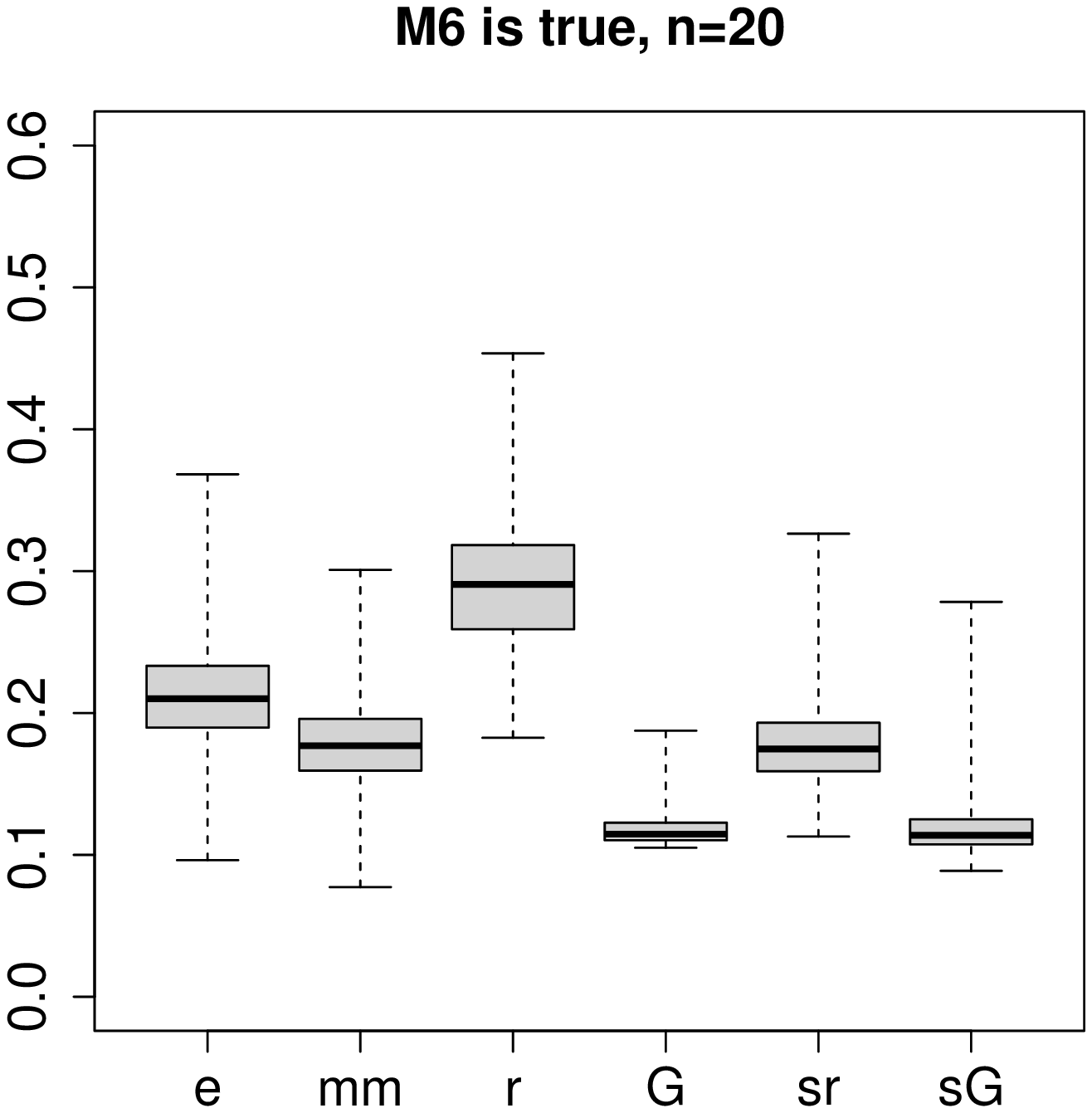}
  \end{subfigure}
  \begin{subfigure}{3.9cm}
    \centering\includegraphics[scale=0.211]{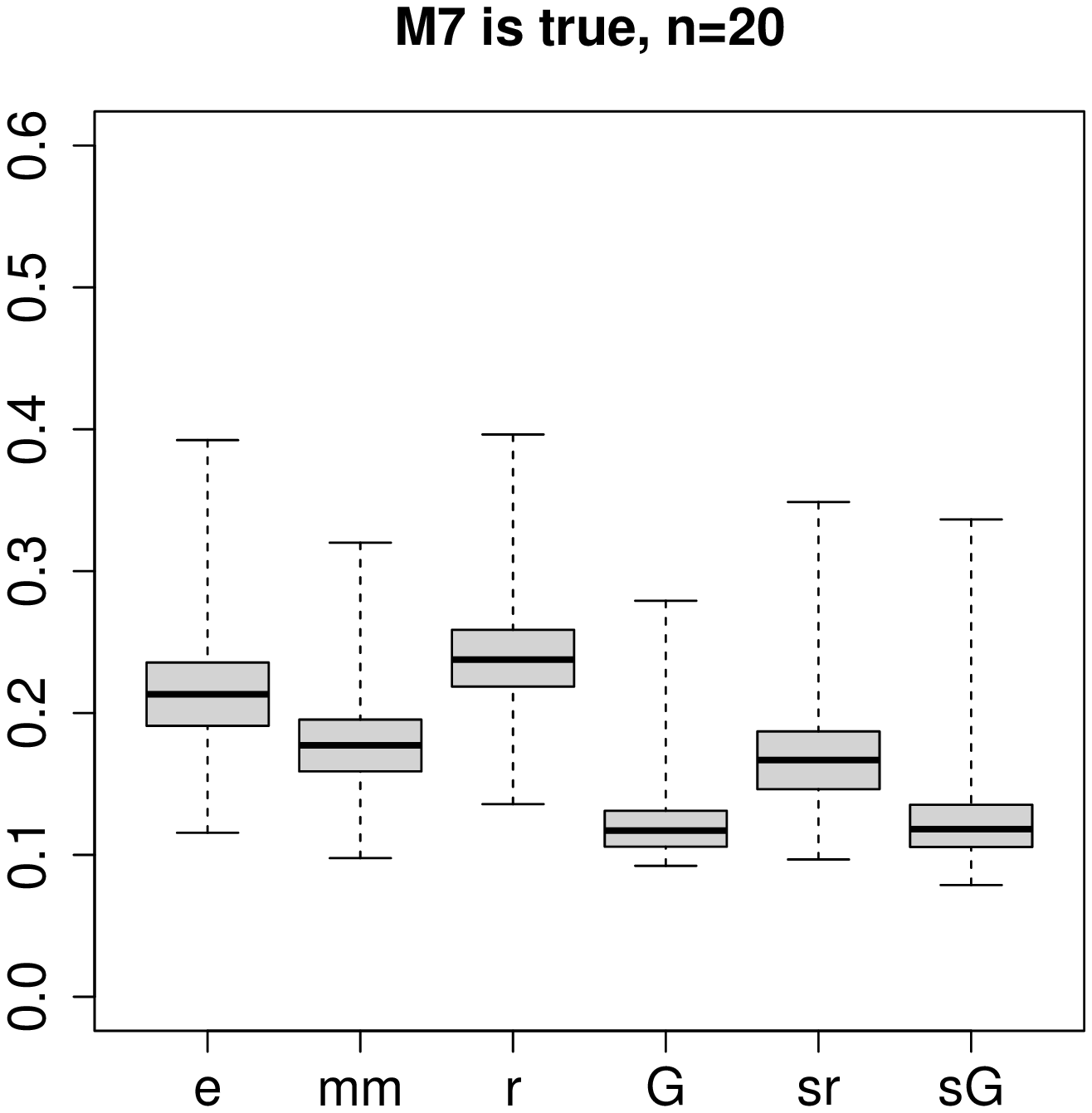}
  \end{subfigure}
  
    \begin{subfigure}{3.9cm}
    \centering\includegraphics[scale=0.211]{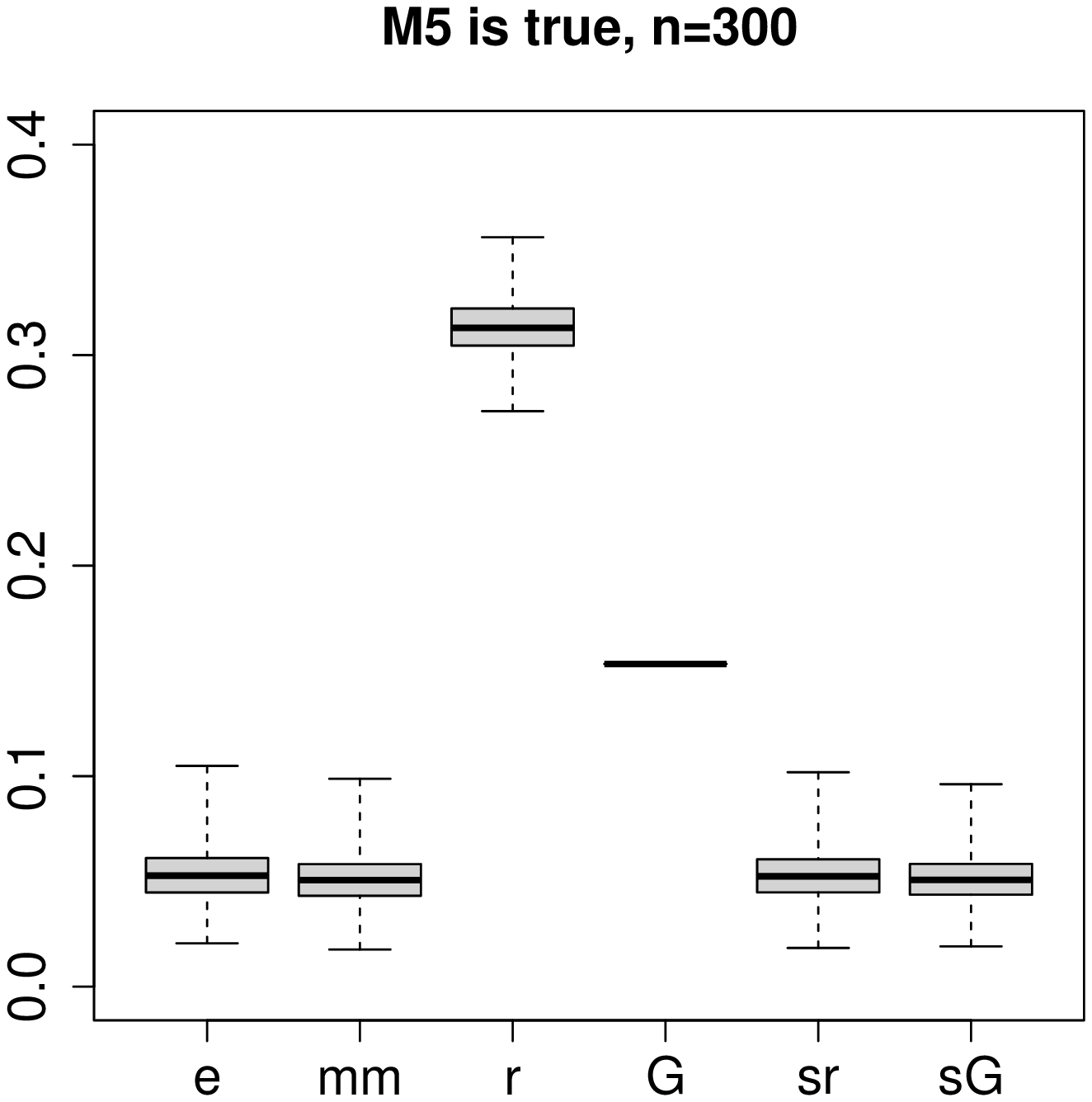}
  \end{subfigure}
  \begin{subfigure}{3.9cm}
    \centering\includegraphics[scale=0.211]{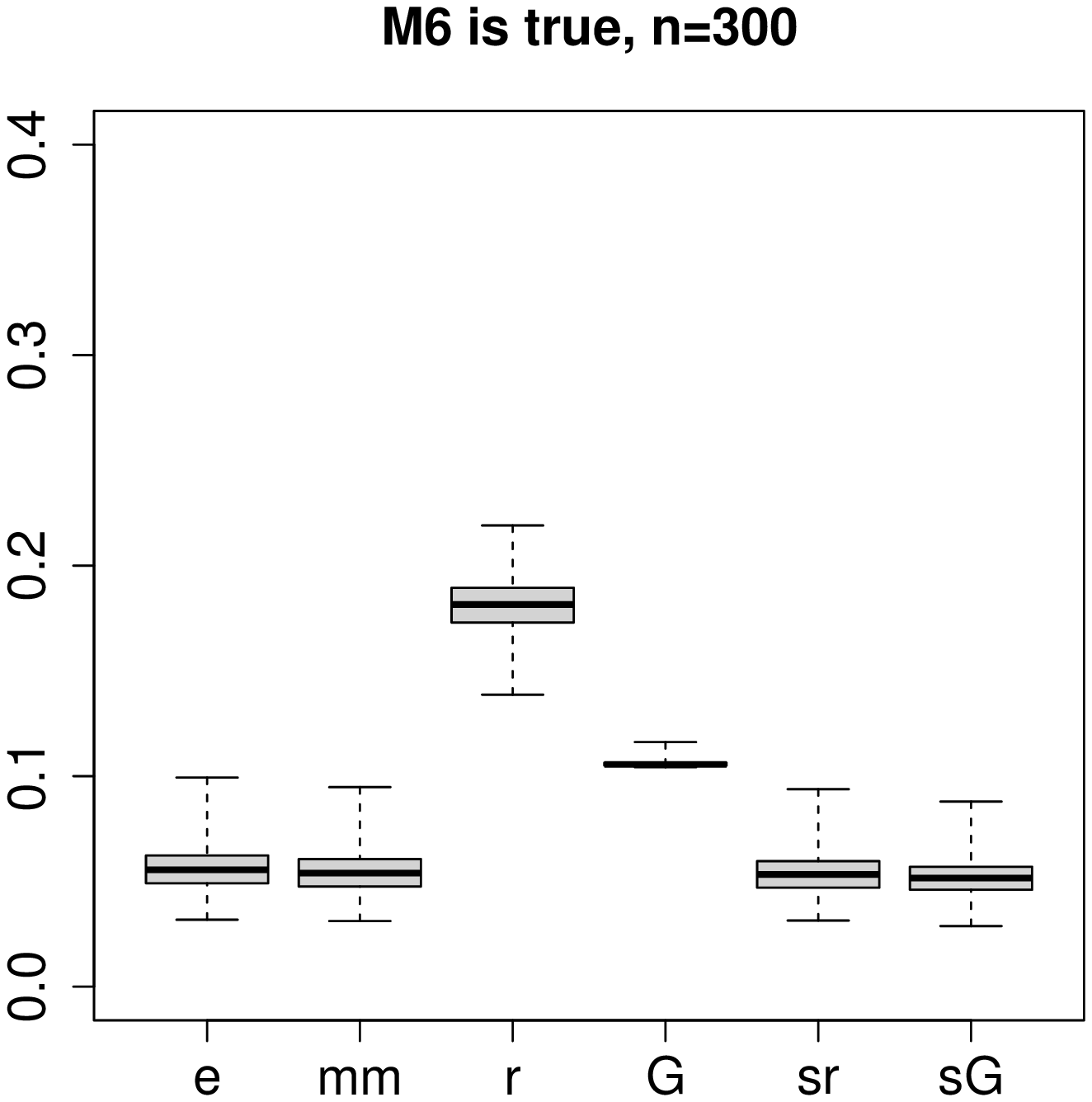}
  \end{subfigure}
  \begin{subfigure}{3.9cm}
    \centering\includegraphics[scale=0.211]{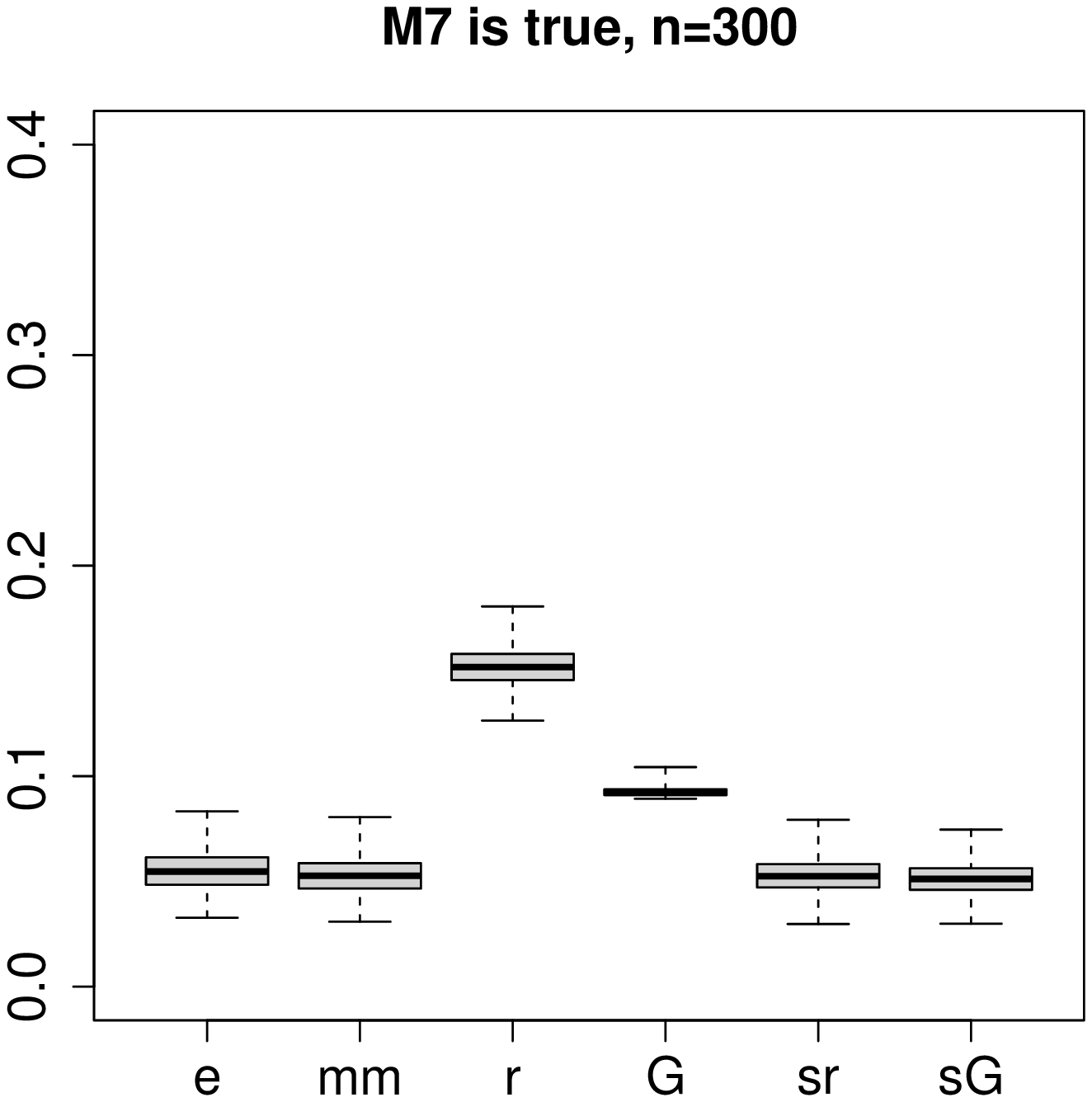}
  \end{subfigure}
   \caption{The boxplots for $\ell_2$-distances of the estimators:the empirical estimator $(e)$, minimax estimator $(mm)$, rearrangement estimator (r), Grenander estimator $(G)$, the stacked rearrangement estimator (sr) and the stacked Grenander estimator $(sG)$ for the models \textbf{M5}, \textbf{M6} and \textbf{M7}.}\label{nondecr_l2}
\end{figure}

From Figure \ref{nondecr_l1} and Figure \ref{nondecr_l2} we can conclude that stacked Grenander estimator outperforms in $\ell_{1}$ and $\ell_{2}$ norms the empirical, rearrangement and minimax estimators ("minimax" for the cases of Negative Binomial and Poisson mixture).

Next, it is interesting to note that even if the underlying distribution is not monotone, Grenander estimator can still outperform the empirical estimator in both $\ell_{1}$ and $\ell_{2}$ norms for small sample size. This happens because the isotonisation decreases the variance of the estimator though bias becomes larger. 
 

\begin{figure}[!htbp] 
  \begin{subfigure}{3.9cm}
    \centering\includegraphics[scale=0.25]{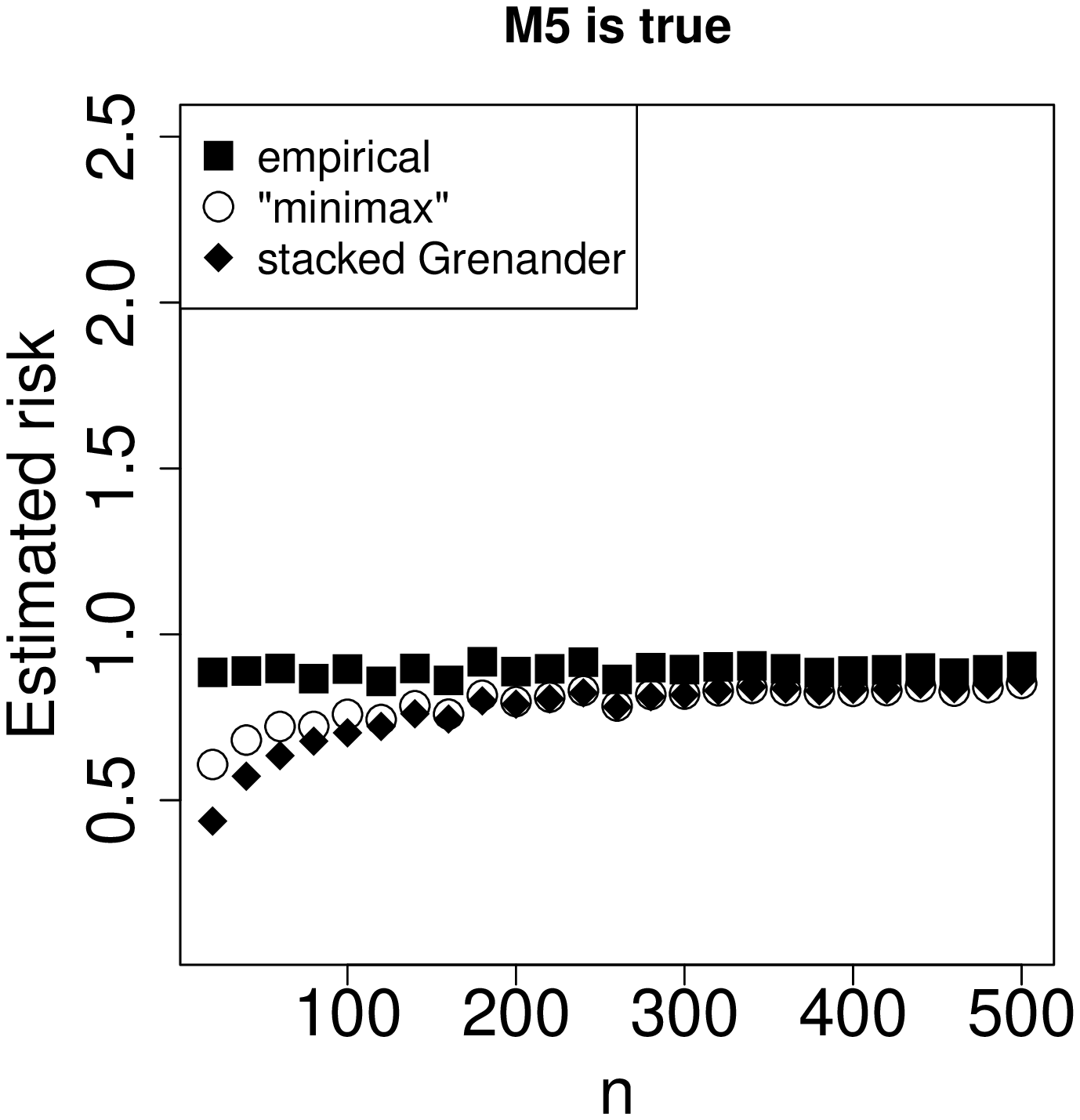}
  \end{subfigure}
  \begin{subfigure}{3.9cm}
    \centering\includegraphics[scale=0.25]{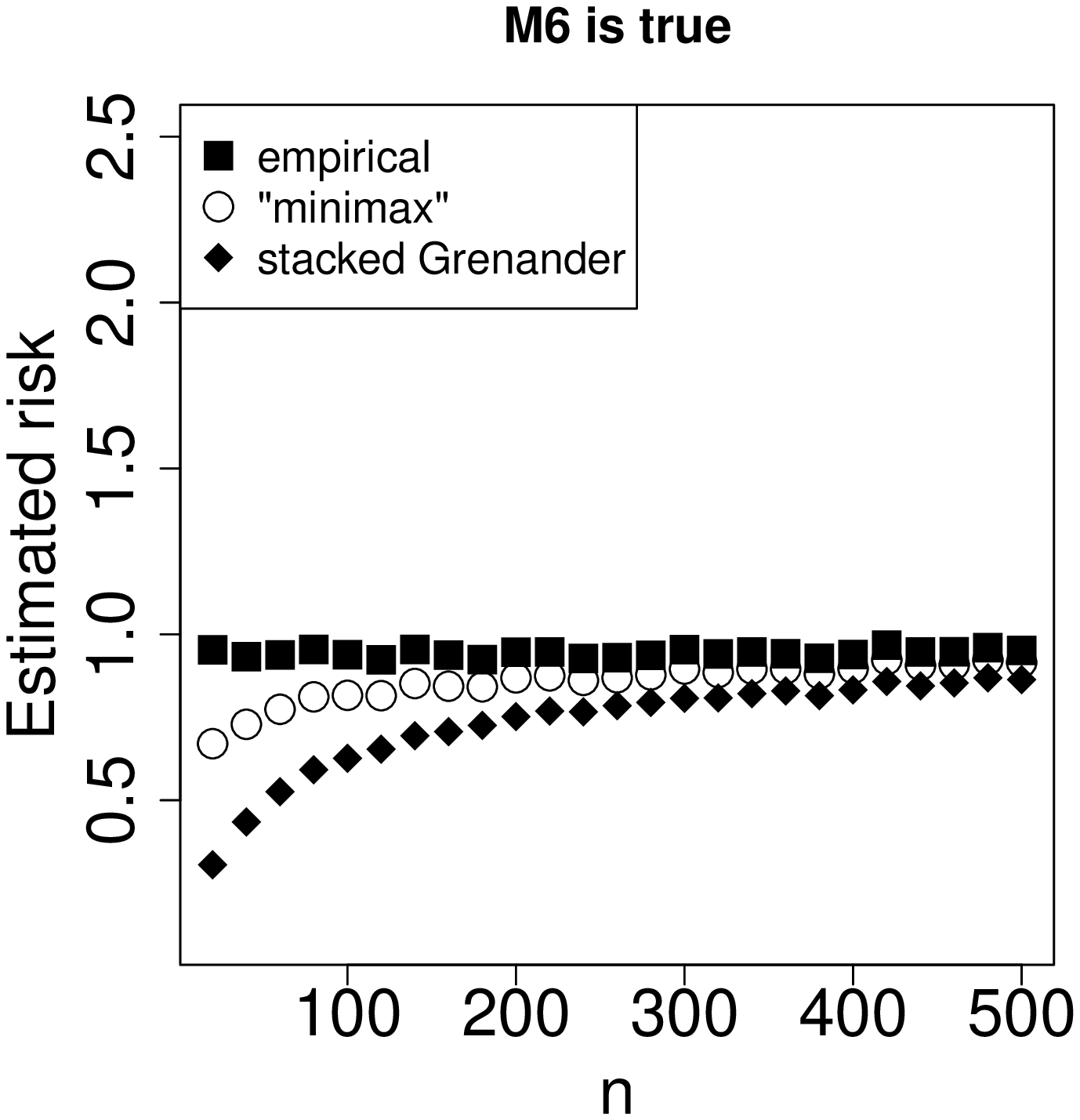}
  \end{subfigure}
  \begin{subfigure}{3.9cm}
    \centering\includegraphics[scale=0.25]{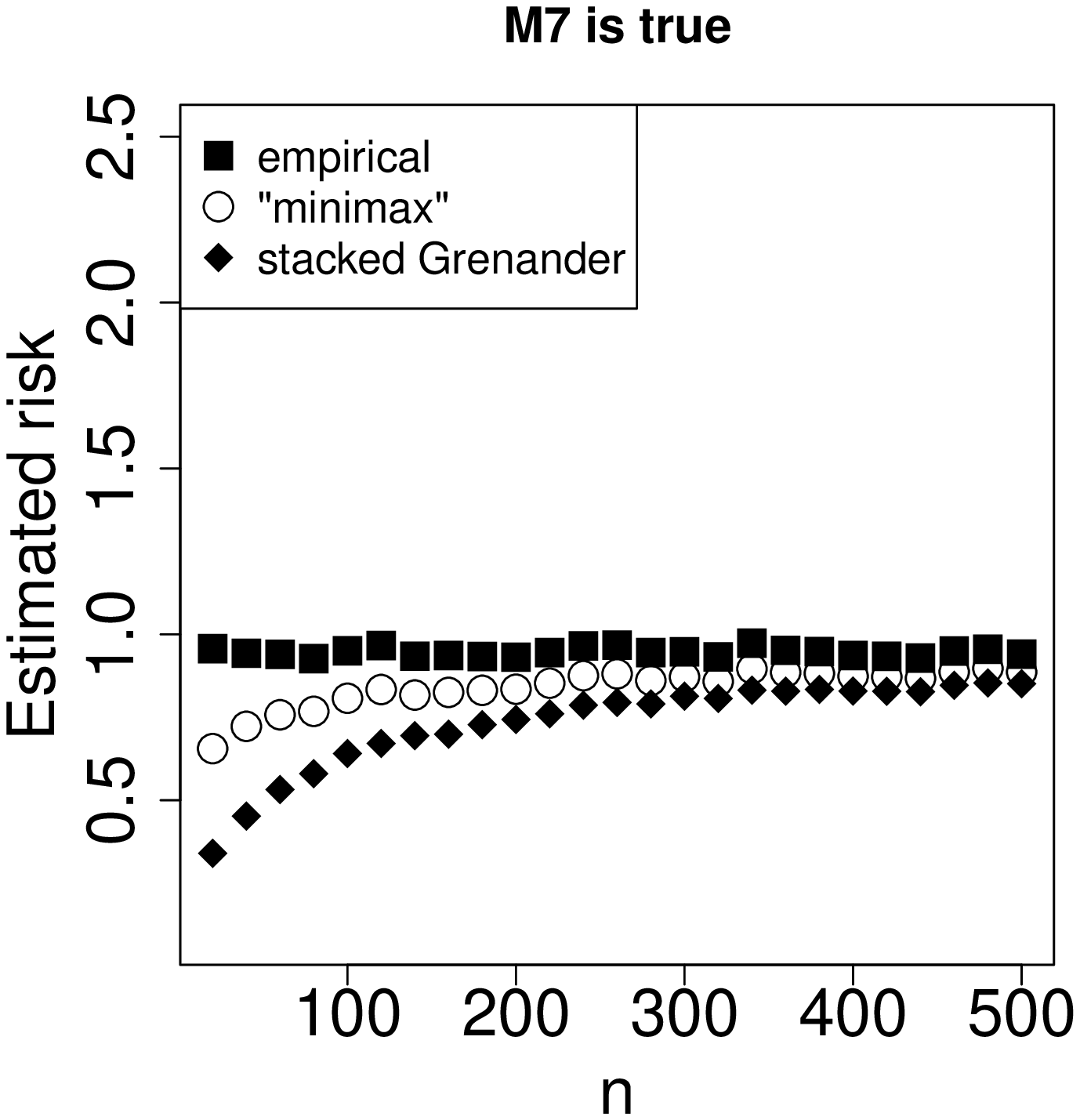}
  \end{subfigure}
  
   \caption{The estimates of the scaled risk for the models \textbf{M5}, \textbf{M6} and \textbf{M7}.}\label{ndecrRisk}
\end{figure}

Let us summarise the results at Figure \ref{ndecrRisk} by plotting the estimates of the scaled risk $n\mathbb{E}[||\hat{\bm{\xi}}_{n} - \bm{p}||_{2}^{2}]$ (with $\hat{\bm{\xi}}_{n}$ one of the following estimators: empirical, minimax or stacked Grenander estimator). Note that in the case of non-decreasing true p.m.f. we do not plot the risk for Grenander estimator, because, obviously, in the miss-specified case the scaled risks of the constrained estimators are worse than the risk of consistent estimators. Based on the simulations we can conclude that stacked Grenander estimator performs better than empirical and minimax estimators even when the underlying distribution is not decreasing. 

The result might look surprising at the first sight. Nevertheless, the explanation of the effect of $\ell_{2}$-risk reduction by stacking empirical estimator with some fixed probability vector was explained in \cite{fienberg1973}. Further, let us consider the case of model \textbf{M5}, i.e. very non-decreasing case when the underlying distribution is strictly increasing. Then, since the empirical estimator is strongly consistent there exist a random $n_{1}$ such that for all $n>n_{1}$ the vector is $\hat{\bm{p}}_{n}$ is strictly increasing almost surely. Next, note that from Lemma \ref{maxsetalg} it follows that for all $n>n_{1}$ we have $\hat{g}_{j} = 1/12$, for all $j=0,\dots 11$, almost surely. Therefore, for $n>n_{1}$ stacked Grenander estimator becomes the stacking of the empirical estimator with a uniform distribution $U(11)$ almost surely, which is similar to what, for example, minimax estimator in  (\ref{estMM}) does. One can also see from Figure \ref{ndecrRisk} that in the case of model \textbf{M5} stacked Grenander estimator performs very similarly to the minimax estimator in a sense of $\ell_{2}$-risk. 

\subsection{Coverage probabilities for the confidence bands}
The Table \ref{emp-cov} presents the proportion of times that
\begin{equation*}
\max\Big((\hat{\phi}_{n,j} - \frac{\hat{q}_{\alpha}}{\sqrt{n}}), 0\Big) \leq p_{j} \leq  \hat{\phi}_{n,j} + \frac{\hat{q}_{\alpha}}{\sqrt{n}}, \, \text{for all } j\in\mathbb{N}
\end{equation*}
among 1000 runs for the models \textbf{M1}--\textbf{M7}. The quantiles $\hat{q}_{\alpha}$ are estimated based on 100000 Monte-Carlo simulations.

First, one can see that the proposed global confidence band performs well. Second, note that for the decreasing p.m.f (models \textbf{M1}-\textbf{M4}) the coverage probabilities mostly larger than $0.95$, while for non-decreasing p.m.f (models \textbf{M5}-\textbf{M7}) the coverage probabilities are closer to $0.95$ when $n$ becomes large, because in the former case the confidence band is asymptotically conservative, while in the later case it is asymptotically correct. 

\begin{table}[!htbp]
\caption{\label{emp-cov} Empirical coverage probabilities for the confidence bands for $\alpha = 0.05$ of the empirical estimator (e), stacked rearrangement estimator (sr) and stacked Grenander estimator (sG).}
\centering
\begin{tabular}{ |l|l|l|l|l|l|l|l|l| }
\hline
\textbf{Estimator} & n & \textbf{M1} & \textbf{M2} & \textbf{M3} & \textbf{M4} & \textbf{M5} & \textbf{M6} & \textbf{M7} \\ \hline
\multirow{3}{*}{e} 
 & 100 & 0.961 & 0.961 & 0.957 & 0.956 & 0.963 & 0.973 & 0.971 \\
 & 1000 & 0.945 & 0.945 & 0.952 & 0.949 & 0.953 & 0.964 & 0.956 \\
 & 5000 & 0.955 & 0.943 & 0.95 & 0.955 & 0.945 & 0.953 & 0.951 \\ \hline
\multirow{3}{*}{sr} 
 & 100 & 0.994 & 0.994 & 0.981 & 0.982 & 0.969 & 0.996 & 0.996 \\
 & 1000 & 0.994 & 0.985 & 0.972 & 0.952 & 0.95 & 0.973 & 0.959 \\
 & 5000 & 0.996 & 0.981 & 0.97 & 0.959 & 0.945 & 0.954 & 0.949 \\ \hline
\multirow{3}{*}{sG} 
 & 100 & 0.996 & 0.994 & 0.979 & 0.981 & 0.989 & 0.999 & 0.997 \\
 & 1000 & 0.998 & 0.984 & 0.971 & 0.951 & 0.953 & 0.976 & 0.963 \\
 & 5000 & 0.997 & 0.984 & 0.97 & 0.959 & 0.945 & 0.954 & 0.953 \\ \hline
\end{tabular}
\end{table}

\subsection{Computational times}
First, note, that in general the complexity of the solution for the mixture parameter $\hat{\beta}_{n}$ depends on the largest order statistic $t_{n}$. In Table \ref{comp-times} we provide the "worst case" computational times, i.e. we compute $\hat{\beta}_{n}$ for the estimator based on the following strictly increasing frequency data vector $\bm{x}' = (x'_{0}, \dots, x'_{s})$, with $x'_{j} = j+1$ for the different values of  $s$, averaged over $10$ runs for every $s$. 

\begin{table}[!h]
\caption{\label{comp-times} The "worst case" averaged over 10 runs computational times of the mixture parameter $\hat{\beta}_{n}$ for stacked rearrangement (SR) and stacked Grenander (SG) estimators for different sizes $s$ of the frequency data vector $\bm{x}$.}
\centering
\begin{tabular}{ |l|l|l|l|l|l| }
\hline
\textbf{Estimator} & \textbf{s=500} & \textbf{s=1000} & \textbf{s=3000} & \textbf{s=5000} \\ \hline
\multirow{1}{*}{SR} 
 & 0.4 s & 2.6 s & 3.1 m & 14.1 m  \\ \hline
\multirow{1}{*}{SG} 
 & 0.3 s & 1.6 s & 3.0 m & 14.0 m  \\ \hline
\end{tabular}
\end{table}

Second, recall that in order to compute the confidence band for a given estimated distribution $\hat{\bm{\theta}}_{n}$ for estimation of the coverage probability in Table \ref{emp-cov} we performed 100000 Monte-Carlo simulations of the multivariate normal distribution $\mathcal{N}(\bm{0}, \bm{\Sigma}(\hat{\bm{\theta}}_{n}))$, with $\Sigma_{i,j}(\hat{\bm{\theta}}_{n}) = \hat{\theta}_{n, j}\delta_{i, j} - \hat{\theta}_{n, i}\hat{\theta}_{n, j}$ $(i, j = 0, \dots, t_{n})$ to estimate the quantile $\hat{q}_{\alpha}$. The Table \ref{comp-times-q} shows the averaged over $10$ runs computational times of the estimation of $\hat{q}_{\alpha}$ of $\mathcal{N}(\bm{0}, \bm{\Sigma}(\bm{\theta}))$ for a fixed non-random p.m.f. vector $ \bm{\theta} = T^{d}(s)$ (recall that $T^{d}(s)$ is a strictly decreasing triangular function), for different values of $s$ based on 100000 Monte-Carlo simulations.
\begin{table}[!h]
\caption{\label{comp-times-q} The averaged over 10 runs computational times of the quantile $\hat{q}_{\alpha}$ for different values of the support size $s$.}
\centering
\begin{tabular}{ |l|l|l|l|l| }
\hline
\textbf{s=500} & \textbf{s=1000} & \textbf{s=3000} & \textbf{s=5000} \\ \hline
 14.9 s & 49.6 s & 7.8 m & 22.0 m  \\ \hline
\end{tabular}
\end{table}

All the computations were performed on MacBook Air (Apple M1 chip), 16 GB RAM. We can conclude that both stacked rearrangement and stacked Grenander estimators are computationally feasible. 

\section{Conclusion and discussion}\label{concl}
In this paper we introduced and studied estimation of a discrete infinitely supported distribution by stacking the empirical estimator with Grenander estimator and the empirical estimator with rearrangement estimator. 

The main results of the paper: the stacked Grenander estimator is computationally feasible, it outperforms the empirical estimator, and it is almost as good as Grenander estimator for the case of decreasing true p.m.f. Also, stacked Grenander estimator outperforms the stacked rearrangement estimator, except for the case of a strictly decreasing p.m.f. The same effect was shown in \cite{jankowski2009estimation} for rearrangement and Grenander estimators in the case when underlying p.m.f. is decreasing. We proved that even when the true distribution is not decreasing, the estimator remains strongly consistent with $\sqrt{n}$-rate of convergence. Therefore, the stacked Grenander estimator provides a trade-off between goodness of fit and monotonicity. 


The first natural generalisation of stacked Grenander estimator could be stacking with isotonic regression for a general isotonic constraint (cf. Appendix for the definition). Throughout the paper, in almost all the proofs we used properties of a general isotonic regression, cf. Lemma \ref{propisot}. However, the proof of Lemma \ref{loois} is based on the maximum upper sets algorithm, which is given in Lemma \ref{maxsetalg} in Appendix, and this algorithm is valid only for one dimensional monotonic case. Therefore, the generalisation of stacked Grenander estimator to the general isotonic case for finite support is straightforward, though the case of an infinite support remains an open problem. 

Second, it is also important to consider other shape constraints, such as unimodal, convex and log-concave cases. Stacking these estimators is, in effect, similar to the generalisation of nearly-isotonic regression to the nearly-convex regression in \cite{tibshirani2011}.

Third, in this work we studied the case of discrete distribution with infinite support. The empirical estimator is closely related to estimation of probability density functions via histograms. Therefore, another direction is stacking the histogram estimators with isotonised histogram. 

Forth, as mentioned in the introduction, the constrained stacked estimators have not been investigated for the case of continuous density. The interesting property of Grenander estimator in a continuous case is that the distributional pointwise rate of convergence depends on the local behaviour of the underlying distribution: if the true distribution is flat, the Grenander estimator has $n^{1/2}$-rate of convergence cf. \cite{carolan1999}, and  $n^{1/3}$-rate otherwise, cf. \cite{rao1969}.  Therefore, in the case of a continuous support it would be interesting to study stacking, for example, Grenander estimator and kernel density estimator.

Another interesting direction of research concerns the stacking with a cross-validation based on other loss functions. For the overview and theoretical properties of different loss functions for evaluation of discrete distributions we refer to the paper \cite{haghtalab2019}.

Finally, as we mentioned in the introduction, the problem of stacking shaped constrained regression estimators has not been studied much.  Therefore, since stacked Grenander estimator performs quite well, it would be interesting to explore, for example, the prediction performance of stacked isotonic regression. 

\section{Appendix}
We start with the definition of a general isotonic regression. Let $\mathcal{J} = \{j_{1}, \dots, j_{s}\}$, with $s \leq \infty$, be some index set. Next, let us define the following binary relations on $\mathcal{J}$:

A binary relation $\preceq$ on $\mathcal{J}$ is a simple order if 
\begin{enumerate}[label=(\roman*)]
\item it is reflexive, i.e. $j \preceq j$ for $j \in \mathcal{J}$;
\item it is transitive, i.e. $j_{1}, j_{2}, j_{3} \in \mathcal{J}$, $j_{1} \preceq j_{2}$ and $j_{2} \preceq j_{3}$ imply $j_{1} \preceq j_{3}$;
\item it is antisymmetric, i.e. $j_{1}, j_{2} \in \mathcal{J}$, $j_{1} \preceq j_{2}$ and $j_{2} \preceq j_{1}$ imply $j_{1} = j_{2}$;
\item every two elements of $\mathcal{J}$ are comparable, i.e. $j_{1}, j_{2} \in \mathcal{X}$ implies that  either $j_{1} \preceq j_{2}$ or $j_{2} \preceq j_{1}$.
\end{enumerate}
A binary relation $\preceq$ on $\mathcal{J}$ is a partial order if it is reflexive, transitive and antisymmetric, but there may be noncomparable elements. A pre-order is reflexive and transitive but not necessary antisymmetric and the set $\mathcal{J}$ can have noncomparable elements. Note, that in some literature the pre-order is called as a quasi-order.

Next, a vector $\bm{v}$ with the elements indexed by $\mathcal{J}$ is isotonic if $j_{1} \preceq j_{2}$ implies $v_{j_{1}} \leq v_{j_{2}}$. We denote the set of all isotonic square summable vectors  by $\bm{\mathcal{F}}^{is}$, which is also called isotonic cone.

Furthermore, a vector $\bm{v}^{*}\in \mathbb{R}^{s}$, with $s \leq \infty$, is the isotonic regression of an arbitrary vector $\bm{v} \in \mathbb{R}^{s}$ (or $\bm{v} \in \ell_{2}$, if $s=\infty$) over the pre-ordered index set $\mathcal{J}$ if 
\begin{eqnarray*}
\bm{v}^{*} = \underset{\bm{f} \in \bm{\mathcal{F}}^{is}}{\argmin} \sum_{j \in \mathcal{J}}(f_{j} - v_{j})^{2}.
\end{eqnarray*}

In Lemma \ref{propisot} we provide properties of a general isotonic regression which are referred to in the paper.
\blem\label{propisot}[Properties of a general isotonic regression]
Let $\bm{v}^{*}_{n}\in \ell_{2}$ be the isotonic regressions of some set of vectors $\bm{v}_{n} \in \ell_{2}$, for $n = 1, 2 \dots$.
Then, the following holds.
\begin{enumerate}[label=(\roman*)]
\item $\bm{v}^{*}_{n}$ exists and it is unique.
\item $\sum_{j} v_{n,j} = \sum_{j} v^{*}_{n, j}$, for all $n = 1, 2, \dots$. 
\item $\bm{v}^{*}_{n}$, viewed as a mapping from $\ell_{2}$ into $\ell_{2}$, is continuous. 
\item $\bm{v}^{*}_{n}$ satisfies the same bounds as the basic estimator, i.e. $a \leq v^{*}_{n, j} \leq b$, for all $n = 1, 2, \dots $ and $j = 1, 2, \dots$. 
\item $\Pi(a\bm{v}_{n} |\mathcal{F}^{is})= a\Pi(\bm{v}_{n}|\mathcal{F}^{is})$ for all $a \in \mathbb{R}_{+}$.\\
\end{enumerate}  
\elem
\textbf{Proof}. Statements $(i)$, $(ii)$ and $(iii)$ follow from Theorem 8.2.1, Corollary B of Theorem 8.2.7 and Theorem 8.2.5, respectively, in \cite{robertsonorder}, statements $(iv)$, $(v)$ and $(vi)$ follow from Corollary B of Theorem 7.9, Theorems 7.5, respectively, in \cite{barlowstatistical}.
\eop

In the next lemma we describe the maximum upper sets algorithm for the solution to the isotonic regression in the monotone case.
\begin{lem}\label{maxsetalg}[Maximum upper sets algorithm]
For a given $\bm{x}\in\mathbb{R}^{t+1}_{+}$ the solution $\bm{x}^{*}$ of a simple order isotonic regression, i.e. 
\begin{equation*}
\bm{x}^{*} =  \underset{f_{0} \geq f_{1} \geq  \dots \geq f_{t_{n}}}{\argmin}\sum_{j=0}^{t_{n}}[x_{j} - f_{j}]^{2}
\end{equation*}
is given by the following algorithm. First, let us define $m(-1) = -1$. Second, we choose $m(0) > m(-1)$ to be the largest integer which maximizes the following mean 
\begin{equation*}
\frac{\sum\limits_{k=m(-1)+1}^{m(0)}x_{k}}{m(0)-m(-1)}.
\end{equation*}
Next, let us choose $m(1) > m(0)$ to be the largest integer which maximizes
\begin{equation*}
\frac{\sum\limits_{k=m(0)+1}^{m(1)}x_{k}}{m(1)-m(0)}.
\end{equation*}
We continue this process and get
\begin{equation*}
-1=m(-1) < m(1) < \dots < m(l) = t_{n}.
\end{equation*}
The solution $\bm{x}^{*}$ (i.e. the isotonic regression of $\bm{x}$) is given by
\begin{equation*}
x^{*}_{j} =  \frac{\sum\limits_{k=m(r-1)+1}^{m(r)}x_{k}}{m(r)-m(r-1)}
\end{equation*}
for $j\in [m(r-1)+1, m(r)]$ and $r \in [0,  l]$.
\end{lem} 
\textbf{Proof}. The proof is given on p. 77 in \cite{barlowstatistical} and p. 26 in \cite{robertsonorder}, and, also, for simpler explanation of the algorithm we refer to \cite{wrigh1982}.\eop

\textbf{Proof of Lemma \ref{pwconv}.}
In order to prove the statement of the lemma, we show that the pointwise convergence almost surely of $\hat{\bm{p}}^{\backslash [j]}_{n}$, $\hat{\bm{r}}^{\backslash [j]}_{n}$ and $\hat{\bm{g}}^{\backslash [j]}_{n}$ for a fixed $j$ holds. First, note that for $j$ such that $p_{j} = 0$ the statement holds, since in this case we have
\begin{equation*}
\hat{\pi}_{n,j}=\hat{\rho}_{n,j}=\hat{\gamma}_{n,j} = 0
\end{equation*}
for all $n$ almost surely.

Second, let us fix some $0 \leq j \leq t_{n}$, such that $p_{j} \neq 0$. Next, clearly, 
\begin{equation}\label{pwclooem}
\hat{\bm{p}}^{\backslash [j]}_{n} \stackrel{a.s.}{\to} \bm{p}
\end{equation}
in $\ell_{k}$-norm for $1 \leq k \leq \infty$. Next, from (\ref{pwclooem}) for the sequence $\hat{\bm{g}}^{\backslash [j]}_{n}$ we have 
\begin{equation*}
\hat{\bm{g}}^{\backslash [j]}_{n} = \Pi \big(\hat{\bm{p}}^{\backslash [j]}_{n}|\mathcal{F}^{decr}\big) \stackrel{a.s.}{\to} \bm{g}
\end{equation*}
in $\ell_{2}$-norm, since the isotonic regression is a continuous map (cf. statement (iii) in Lemma \ref{propisot}). Therefore, we have proved the statement of the lemma for the sequences  $\hat{\bm{\pi}}_{n}$ and $\hat{\bm{\gamma}}_{n}$. 

Next, we prove the statement for $\hat{\bm{\rho}}_{n}$. Let us fix some $s > j$ such that $p_{k} < p_{j}$ for all $k > s$. Next, let
\begin{equation*}
\hat{\mathfrak{p}}_{(k)} = \text{the $k$th largest of $\{\hat{p}^{\backslash [j]}_{n,0}, \dots, \hat{p}^{\backslash [j]}_{n,t_{n}} \}$}. 
\end{equation*}
Further, from (\ref{pwclooem}) it follows that there exists $n_{1}$ such that for all $n > n_{1}$
\begin{equation*}
[\hat{\bm{r}}^{\backslash [j]}_{n}]^{(0,j)} = \{\hat{\mathfrak{p}}_{(1)}, \dots, \hat{\mathfrak{p}}_{(j)} \} \subset \{\hat{p}^{\backslash [j]}_{n,0}, \dots, \hat{p}^{\backslash [j]}_{n,s}  \},
\end{equation*}
almost surely, where $[\cdot]^{(0,j)}$ denotes the first $(j+1)$ elements of the vector. Finally, since the rearrangement operator is continuous map in a finite dimensional case (Lemma 6.1 in \cite{jankowski2009estimation}), the result of the lemma follows from continuous mapping theorem.
\eop

\textbf{Proof of Theorem \ref{thmLSCV}.} Recall that the least-squares cross-validation criterion is given by
\begin{equation*}\label{}
\begin{aligned}
CV(\beta) ={} & \sum_{j=0}^{t_{n}}\hat{\phi}_{n, j}^{2}  - 2\sum_{j=0}^{t_{n}} \hat{p}_{n,j}\hat{\phi}^{\backslash [j]}_{n, j} =\\
  &\sum_{j=0}^{t_{n}}(\beta \, \hat{h}_{n, j} + (1-\beta)\hat{p}_{n, j})^{2}- 2\sum_{j=0}^{t_{n}} \hat{p}_{n,j}(\beta \, \hat{h}^{\backslash [j]}_{n, j} + (1-\beta)\hat{p}^{\backslash [j]}_{n, j}).
\end{aligned}
\end{equation*}

Then, after simplification we get 
\begin{equation*}\label{}
\begin{aligned}
CV(\beta) ={} &a_{n}\beta^{2} - 2b_{n}\beta + c_{n},
\end{aligned}
\end{equation*}
where  the term $c_{n}$ does not depend on $\beta$, and 
\begin{equation*}
a_{n} = \sum_{j=0}^{t_{n}}(\hat{h}_{n, j} - \hat{p}_{n, j})^{2},
\end{equation*}
and 
\begin{equation*}
b_{n} = \sum_{j=0}^{t_{n}}\hat{p}_{n, j}(\hat{h}_{n, j}^{\backslash[j]} - \hat{p}_{n, j}^{\backslash[j]}) - \sum_{j=0}^{t_{n}}\hat{p}_{n, j}(\hat{h}_{n, j} - \hat{p}_{n, j}).
\end{equation*}
Assume, that $a_{n} \neq 0$. Then, $CV(\beta)$ is minimised by
\begin{equation*}\label{}
  \beta_{n}=\begin{cases}
    \frac{b_{n}}{a_{n}}, & \text{if } \, 0 \leq b_{n} \leq a_{n}, \\
    1, & \text{if } \,  a_{n} \leq b_{n}, \\
    0, & \text{if } \, b_{n} \leq 0. 
  \end{cases}
\end{equation*}
Next, note that if $\hat{\bm{p}}_{n} = \hat{\bm{h}}_{n}$, then $\hat{\bm{\phi}}_{n} = \hat{\bm{p}}_{n} = \hat{\bm{h}}_{n}$ for any $0 \leq \beta_{n} \leq 1$, and, therefore, for consistency of notation we define $\hat{\beta}_{n} = 0$ when $a_{n} = 0$.
\eop

\textbf{Proof of Lemma \ref{loois}}.
First, we prove the statement for $\hat{\bm{\pi}}_{n}$. Assume that for some $j$ we have $\hat{p}^{\backslash [j]}_{n,j}\neq0$ and recall that 
\begin{equation*}
\hat{\pi}_{n, j} = \hat{p}^{\backslash [j]}_{n, j} =\frac{x_{j}-1}{n-1}.
\end{equation*}
Next, note that
\begin{equation*}
\frac{x_{j}-1}{n-1} - \frac{x_{j}}{n} = \frac{-n+x_{j}}{n(n-1)} < 0.
\end{equation*}

Let us study the case of $\hat{\bm{\gamma}}_{n}$. 
To prove the statement of the lemma we will use maximum upper sets algorithm, which is given in Lemma \ref{maxsetalg} in the Appendix. Let $\bm{x} = (x_{0}, \dots, x_{t_{n}})$ be frequency data from $\bm{p}$. Next, let  $\bm{x}^{*} = (x_{0}^{*}, \dots, x_{t_{n}}^{*})$ be the isotonic regression of $\bm{x}$ and assume that $\bm{x}^{*}$ has $(l+1)$ constant regions. Let 
\begin{equation*}
m(0) < \dots < m(l) = t_{n}
\end{equation*}
be the indices of the last elements in the constant regions of $\bm{x}^{*}$ and $m(-1) = -1$. Therefore, we have
\begin{equation*}
x^{*}_{j} =  \frac{\sum\limits_{k=m(r-1)+1}^{m(r)}x_{k}}{m(r)-m(r-1)}
\end{equation*}
for $j\in [m(r-1)+1, m(r)]$ and $r \in [0, l]$.

Let us consider the first constant region of $\bm{x}^{*}$ and for some integer $q\in[0, m(0)]$ define vector $\bm{y}\in\mathbb{R}^{t+1}_{+}$
\begin{equation*}\label{}
  y_{j}=\begin{cases}
    x_{j}-1, &\text{ if } \, j=q, \\
    x_{j}, &\text{ otherwise},
  \end{cases}
\end{equation*}
and let $\bm{y}^{*}$ be isotonic regression of $\bm{y}$.

Recall, $m(0)$ is the largest non-negative integer which maximizes the following mean 
\begin{equation*}
S_{1} = \frac{\sum\limits_{k=0}^{m(0)}x_{k}}{m(0) + 1}.
\end{equation*}

Further, let $m'(0)$ be the largest non-negative integer which maximizes the following mean for the vector $\bm{y}$ 
\begin{equation*}
S_{2} = \frac{\sum\limits_{k=0}^{m'(0)}y_{k}}{m'(0) + 1}.
\end{equation*}
    
Let us prove that $S_{2} \leq S_{2}$. First, assume that $m'(0) = m(0)$, then, clearly, $S_{2} \leq S_{1}$ since $y_{j} \leq x_{j}$. Second, let us assume that $m'(0) \neq m(0)$. Then, from the definitions of $m(0)$ and $m'(0)$ it follows
\begin{equation*}
S_{2} = \frac{\sum\limits_{k=0}^{m'(0)}y_{k}}{m'(0) + 1} \leq \frac{\sum\limits_{k=0}^{m'(0)}x_{k}}{m'(0) + 1} \leq \frac{\sum\limits_{k=0}^{m(0)}x_{k}}{m(0) + 1} = S_{1}.
\end{equation*}

Next, assume that $q$ is not in the first constant region. Then in this case from maximum upper sets algorithm it follows that the constant regions in the isotonic regressions $\bm{x}^{*}$ and $\bm{y}^{*}$ are the same up to the region which contains element with index $m$. Then, we can use the same approach as for the first region. Therefore, we have proved that $y^{*}_{q} \leq x^{*}_{q}$. 

Next, from statement $(v)$ of Lemma \ref{propisot} for $\hat{\bm{g}}_{n}$ and $\hat{\bm{\gamma}}_{n}$ we have
\begin{equation*}
\hat{g}_{n, j} = \frac{x^{*}_{j}}{n},
\end{equation*}
and
\begin{equation*}
\hat{\gamma}_{n, j} = \frac{y^{*}_{j}}{n-1},
\end{equation*}
therefore, we proved that
\begin{equation*}
\hat{\gamma}_{n, j} \leq \frac{n}{n-1} \hat{g}_{n, j}.
\end{equation*}

Finally, we prove the inequality for $\hat{\bm{\rho}}_{n}$. Analogously to the case of $\hat{\bm{\gamma}}_{n}$, let us consider the vectors $\bm{x}$ and $\bm{y}$, discussed above. Note that $y_{j} \leq x_{j}$ for all $j$, therefore, the same componentwise inequality holds for the sorted vectors $rear(\bm{x})$ and $rear(\bm{y})$. Next, using the definition of $\hat{\bm{r}}_{n}$ and $\hat{\bm{\rho}}_{n}$ we prove that
\begin{equation*}
\hat{\rho}_{n,j} \leq \frac{n}{n-1} \hat{r}_{n,j}.
\end{equation*}
\eop

\textbf{Proof of Theorem \ref{asymdist}}.
Assume that the p.m.f. $\bm{p}$ is not decreasing. Note that
\begin{equation*}
\begin{aligned}
&||\sqrt{n}(\hat{\bm{\phi}}_{n} - \bm{p}) -\sqrt{n}(\hat{\bm{p}}_{n} - \bm{p}))||_{2} ={}  \sqrt{n}||\hat{\bm{\phi}}_{n} - \hat{\bm{p}}_{n}||_{2} \leq\\
&\hat{\beta}_{n}\sqrt{n}||\hat{\bm{h}}_{n} - \hat{\bm{p}}_{n}||_{2} + (1-\hat{\beta}_{n})\sqrt{n}||\hat{\bm{p}}_{n} - \hat{\bm{p}}_{n}||_{2} =
 \hat{\beta}_{n}\sqrt{n}||\hat{\bm{h}}_{n} - \hat{\bm{p}}_{n}||_{2}.
\end{aligned}
\end{equation*} 
Then, since
\begin{equation*}
\begin{aligned}
&||\hat{\bm{r}}_{n} - \hat{\bm{p}}_{n}||_{2}\stackrel{a.s.}{\to} ||\bm{r} - \bm{p}||_{2} < \infty,\\
&||\hat{\bm{g}}_{n} - \hat{\bm{p}}_{n}||_{2}\stackrel{a.s.}{\to} ||\bm{g} - \bm{p}||_{2} < \infty,
\end{aligned}
\end{equation*} 
and using (\ref{convlbd}) we have
\begin{equation*}
\begin{aligned}
&\hat{\beta}_{n}\sqrt{n}||\hat{\bm{r}}_{n} - \hat{\bm{p}}_{n}||_{2} \stackrel{a.s.}\to 0,\\
&\hat{\beta}_{n}\sqrt{n}||\hat{\bm{g}}_{n} - \hat{\bm{p}}_{n}||_{2} \stackrel{a.s.}\to 0,
\end{aligned}
\end{equation*} 
which leads to
\begin{equation*}
\sqrt{n}||\hat{\bm{\phi}}_{n} - \hat{\bm{p}}_{n}||_{2}  \stackrel{a.s.}\to 0.
\end{equation*} 
The statement of the theorem now follows from Theorem 3.1 in \cite{billingsley}.

Assume that $\bm{p}$ is a strictly decreasing p.m.f. over $\{ 0, \dots, s \}$, with $s < \infty$. Next, let $\varepsilon = \inf\{|p_{j}-p_{j+1}|: j=0,\dots, s-1\}$ and note that 
\begin{equation*}
\{\sup_{j}\{|\hat{p}_{j} - p_{j} | < \varepsilon/2\} \subseteq \{\hat{h}_{n,j} = \hat{p}_{n,j}\}
\end{equation*}
for both $\hat{\bm{h}}_{n} = \hat{\bm{g}}_{n}$ and $\hat{\bm{h}}_{n} = \hat{\bm{r}}_{n}$. Therefore, this implies that for any $j = \{0, \dots, s$ we have 
\begin{equation*}
\mathbb{P}[\hat{\phi}_{n,j} = \hat{p}_{n,j}] \geq \mathbb{P}[\sup_{j}\{|\hat{p}_{n,j} - p_{j} | < \varepsilon/2] \to 1,
\end{equation*}
since the empirical estimator is strongly consistent. The statement of the theorem follows from Theorem 3.1 in \cite{billingsley}.
\eop

\section{Acknowledgments}
This work was partially supported by the Wallenberg AI, Autonomous Systems and Software Program (WASP) funded by the Knut and Alice Wallenberg Foundataion.

\end{document}